\documentclass[sn-mathphys]{sn-jnl}

\jyear{2021}%

\theoremstyle{thmstyleone}%
\usepackage{amssymb}
\usepackage{amsfonts}
\usepackage{mathrsfs}
\usepackage{amsmath}
\usepackage{graphicx}
\usepackage{dsfont}
\usepackage{float}
\usepackage{graphics}
\usepackage{titlesec}
\usepackage{braket}
\usepackage{xcolor,graphicx}
\usepackage{verbatim}
\usepackage{mathrsfs}
\usepackage{extarrows}
\usepackage{mathrsfs}
\usepackage{mathtools}
\usepackage{subfig}
\usepackage{hyperref}
\usepackage[percent]{overpic}

\theoremstyle{thmstyletwo}%

\theoremstyle{thmstylethree}%

\newcommand{\bbR}{\mathbb R}
\newcommand{\bA}{\mathbf A}

\newcommand{\bI}{\mathbf I}

\newcommand{\bP}{\mathbf P}

\newcommand{\bV}{\mathbf V}

\newcommand{\bg}{\mathbf g}

\newcommand{\bn}{\mathbf n}
\newcommand{\be}{\mathbf e}

\newcommand{\bu}{{\bf u}}

\newcommand{\bv}{\mathbf v}
\newcommand{\bw}{\mathbf w}

\newcommand{\bx}{\mathbf x}

\newcommand{\bbf}{\mathbf f}

\newcommand{\mT}{\mathcal T}

\newcommand{\divG}{{\mathop{\,\rm div}}_{\Gamma}}
\newcommand{\gradG}{\nabla_{\Gamma}}
\newcommand{\nablaG}{\nabla_{\Gamma}}
\newcommand{\laplG}{\Delta_{\Gamma}}

\newcommand{\OGamma}{\Omega^\Gamma_h}

\newcommand{\tr}{{\rm tr}}

\DeclareGraphicsExtensions{.pdf,.eps,.ps,.eps.gz,.ps.gz,.eps.Y}

\def\cl {\nonumber \\}

\begin{document}

\title[A comparison of CH and NSCH models on manifolds]{A comparison of Cahn--Hilliard and Navier--Stokes--Cahn--Hilliard models on manifolds}


\author[1]{\fnm{Maxim} \sur{Olshanskii}}\email{maolshanskiy@uh.edu}

\author[1]{\fnm{Yerbol} \sur{Palzhanov}}\email{ypalzhanov@uh.edu}

\author*[1]{\fnm{Annalisa} \sur{Quaini}}\email{aquaini@uh.edu}

\affil[1]{\orgdiv{Department of Mathematics}, \orgname{University of Houston}, \orgaddress{\street{3551 Cullen Blvd}, \city{Houston}, \state{Texas}, \postcode{77204} \country{USA}}}


\abstract{We consider phase-field models with and without lateral flow 
for the numerical simulation of lateral phase separation and coarsening in lipid
membranes. For the numerical solution of these models, we apply an unfitted finite element method that 
is flexible in handling complex and possibly evolving shapes in the absence of an explicit surface
parametrization. 
Through several numerical tests, we investigate the effect of the presence of lateral flow 
on the evolution of phases. In particular, we focus on understanding how
variable line tension, viscosity, membrane composition, and surface shape affect the 
pattern formation. 
}




\maketitle

\begin{center}
\emph{Dedicated to Prof.~Alfio Quarteroni on the occasion of his 70$^\text{th}$ birthday}
\end{center}

\section{Introduction}\label{sec:intro}
Recent years have witnessed an increased interest in studying phase separation in biological membranes
\cite{Levental2020,BENNETT20131765}.
This is due to the fact that lateral phase separation
has been recognized as a critical mechanism for dynamic control of the spatial organization of membrane components
\cite{Kahya2003,Niemela2007}. The lipid bilayer in biological membranes may be organized into one 
of two phases: liquid disordered and liquid ordered \cite{Balint2017}. The liquid ordered domains, also known as lipid rafts,
have been linked to a wide range of cellular functions, from membrane trafficking
to inter- and intracellular signaling \cite{Heberle2011}. In addition, domain formation on membranes has also been 
utilized to create novel membrane-based materials with heterogenous surfaces \cite{Bandekar2013,Sempkowski2016}.

Phase separation and pattern formation in lipid bilayers has been studied theoretically (see, e.g., \cite{Andelman_et_al1992,Seifert1993,Kawakatsu_et_al1993,Harden_MacKintosh1994}), experimentally
(see, e.g., \cite{Baumgart_et_al2003,VEATCH20033074}), and numerically (see, e.g., \cite{Marrink_Mark2003,Laradji_Sunil2004,Wang2008,Lowengrub2009,sohn2010dynamics,nitschke_voigt_wensch_2012,Li_et_al2012,Funkhouser_et_al2014}). 
Computational studies are particularly useful to observe the dynamics of phases, which is hard to address
theoretically, and to gain insights that are too expensive (or even impossible) to obtain experimentally. 
In this paper, we choose a continuum-based computational approach that relies on a phase-field description.
Emerged as a powerful computational approach to modeling and predicting phase separation
in materials and fluids, the phase-field method describes the system using a set of field variables
that are continuous across the interfacial regions separating the phases.

In our previous work, we have developed a computationally efficient method based on
the surface Cahn--Hilliard (CH) phase-field model to predict the phase behavior and domain formation on heterogeneous 
membranes \cite{Yushutin_IJNMBE2019,Yushutin2019}. More recently, such method has been 
validated against laboratory experiments \cite{zhiliakov2021experimental}. While a good agreement was achieved
between numerical results and experimental data, the CH model does not account for viscous and fluidic phenomena 
that are recognized to be important in lipid membranes \cite{Kawano2009}. In fact,
it has been demonstrated that membrane fluidity within the liquid ordered 
domains can be substantially lower than that in the liquid disordered phase \cite{SEZGIN20121777}, 
affecting the coarsening dynamics of rafts on membranes \cite{Stanich2013}. In order to capture these phenomena, in \cite{Palzhanov2021} we have considered the more complex surface Navier--Stokes--Cahn--Hilliard (NSCH) model and a numerical method for it.

Although the importance of viscous dissipation and fluidity in lipid membranes is acknowledged, it remains to be 
understood is how lateral flow affects pattern formation. Thus, 
in this paper we compare the evolution of phases as predicted by the CH
model (i.e., without lateral flow) and NSCH model (i.e., with lateral flow)
through a series of numerical tests. For the numerical solution of both models, 
we apply an unfitted finite element method called the trace finite element method (TraceFEM)~\cite{ORG09,olshanskii2017trace}. 
We opted for an unfitted finite element method because of its flexibility in handling complex shapes, as we will show in this
paper, and possibly evolving surfaces, as shown in \cite{Yushutin2019} for the CH model. Although the surfaces treated in 
this paper are steady, our interest in evolving surfaces is associated with our long term goal of simulating 
membrane-based drug carriers that used phase-separated patterns to facilitate fusion with the target cell
\cite{zhiliakov2021experimental}. Among all unfitted finite element methods,
TraceFEM has several advantages that make it appealing: i) it employs a
sharp surface representation, ii) surfaces can be defined implicitly and no surface parametrization is required, 
iii) the number of active degrees of freedom is asymptotically optimal, and iv) the order of convergence is optimal. 

The paper outline is as follows. In Sec.~\ref{sec:model}, we state the two phase-field models and their
variational formulations. The application of TraceFEM to both models is described in Sec.~\ref{sec:method}.
In Sec.~\ref{sec:num_res}, we report several numerical results obtained with both models on the 
surface of a sphere and an asymmetric torus. Sec.~\ref{sec:concl} provides concluding remarks.

\section{Mathematical model}\label{sec:model}

In order to formulate the surface the CH and NSCH equations, 
we need some notation. Let $\Gamma$ be an arbitrary-shaped closed, smooth, and stationary surface, with
the outward pointing unit normal $\bn$.
Let  $\bP=\bP(\bx)\coloneqq\bI -\bn(\bx)\bn(\bx)^T$ for $\bx \in \Gamma$ be the orthogonal projection onto the tangent plane.  
For a scalar function $p:\, \Gamma \to \mathbb{R}$ or a vector function $\bu:\, \Gamma \to \mathbb{R}^3$  we define
$p^e\,:\,\mathcal{O}(\Gamma)\to\mathbb{R}$, $\bu^e\,:\,\mathcal{O}(\Gamma)\to\mathbb{R}^3$ as suitable
extensions of $p$ and $\bu$ from $\Gamma$ to its neighborhood $\mathcal{O}(\Gamma)$.
The surface gradient and covariant derivatives on $\Gamma$ are then defined as
$\nablaG p=\bP\nabla p^e$ and  $\nabla_\Gamma \bu\coloneqq \bP \nabla \bu^e \bP$. These  definitions are  independent of a particular smooth extension of $p$ and $\bu$ off $\Gamma$.
On $\Gamma$ we consider the surface rate-of-strain tensor \cite{GurtinMurdoch75} given by
\begin{equation} \label{strain}
E_s(\bu) \coloneqq \frac12(\nabla_\Gamma \bu + (\nabla_\Gamma \bu)^T).
\end{equation}

The surface divergence operators for a vector $\bg: \Gamma \to \bbR^3$ and
a tensor $\bA: \Gamma \to \mathbb{R}^{3\times 3}$ are defined as:
\[
 \divG \bg  \coloneqq \tr (\gradG \bg), \qquad
 \divG \bA  \coloneqq \left( \divG (\be_1^T \bA),\,
               \divG (\be_2^T \bA),\,
               \divG (\be_3^T \bA)\right)^T,
               \]
with $\be_i$ the $i$th standard basis vector in $\bbR^3$ and $\tr(\cdot)$ is the trace of a matrix. The Laplace-Beltrami operator for a sufficiently smooth function $g$ in a neighborhood of $\Gamma$ is:
\[
\laplG g := \divG (\gradG g).
 \]
Further $L^2(\Gamma)$ is the Lebesgue space of square-integrable functions on $\Gamma$ and $H^1(\Gamma)$ is the Sobolev space of all functions $g\in L^2(\Gamma)$ such that $\nabla_\Gamma g\in L^2(\Gamma)^3$.

On $\Gamma$ we consider a heterogeneous mixture of two species with
surface fractions $c_i = S_i/S$, $i = 1, 2$, where $S_i$
are the surface area occupied by the components and $S$ is the surface area of $\Gamma$.
Since $S = S_1 + S_2$, we have $c_1 + c_2 = 1$. Let $c_1$ be the representative surface fraction, i.e.~$c = c_1$.

\subsection{The Cahn--Hilliard problem}\label{sec:CH}

A well established model for the process of spinodal decomposition and phase separation alone (i.e., 
in the absence of surface fluid flow) is the CH phase-field model \cite{Cahn_Hilliard1958,CAHN1961}. The surface CH equation 
governs the evolution in time $t$ of $c=c(t,\bx)$, $\bx \in\Gamma\subset\mathbb R^3$: 
\begin{equation}\label{surfaceCH}
\frac{\partial c}{\partial t} = \divG \left(M\nabla_\Gamma\left(\frac{1}{\epsilon}f_0'(c) - \epsilon\,\Delta_\Gamma c\right)\right)\quad \text{on}~\Gamma,~\text{for}~t \in (0, T],
\end{equation}
suitably endowed with an initial condition $ c(0, \bx)= c_0$. 
In \eqref{surfaceCH}, $f_0(c) = \frac{1}{4}\,c^2\,(1 - c)^2$ 
is the specific free energy of a homogeneous phase, parameter~$\epsilon > 0$ defines the width of the (diffuse) interface between 
the phases, $M$ is the so-called mobility (see \cite{Landau_Lifshitz_1958}), and
$T$ is the end of a time interval of interest. 
We consider the degenerate mobility of the form
\begin{equation}\label{mobility}
M_c = D c\,(1-c)
\end{equation}
with diffusivity constant~$D > 0$. Mobility \eqref{mobility} is a popular choice for numerical studies.
Eq.~\eqref{surfaceCH} is obtained from
minimizing the total specific free energy $\int_\Gamma \frac{1}{\epsilon} f_0(c) + \frac{1}{2} \epsilon \lvert \nabla_\Gamma c \rvert^2 ds$
subject to the conservation of surface area $\int_\Gamma c\,ds$.

Eq.~\eqref{surfaceCH} is a fourth-order equation. In order to avoid higher order spatial derivatives, which would 
need careful numerical treatment, it is common to rewrite eq.~\eqref{surfaceCH} as
two coupled second-order equations:
\begin{align}
\frac{\partial c}{\partial t} &=  \divG \left(M \gradG \mu \right) \quad \text{on}~\Gamma \label{eq:sys_CH1}, \\
\mu &= \frac{1}{\epsilon} f_0' - \epsilon \laplG c \quad \text{on}~\Gamma. \label{eq:sys_CH2}
\end{align}
In \eqref{eq:sys_CH2}, $\mu$ represents the chemical potential. 

For the numerical method, we need a weak (integral) formulation. In order to
devise it, one multiplies \eqref{eq:sys_CH1} by $v\in H^1(\Gamma)$ and \eqref{eq:sys_CH2} by
$q\in H^1(\Gamma)$, integrates over $\Gamma$ and employs the integration by parts identity.
The weak (variational) formulation of problem \eqref{eq:sys_CH1}-\eqref{eq:sys_CH2} reads: Find $(c,\mu) \in H^1(\Gamma) \times H^1(\Gamma)$ 
such that
\begin{align}
&\int_\Gamma \rho \frac{\partial c}{\partial t} \,v \, ds + \int_\Gamma M \gradG \mu \, \gradG v \, ds = 0, \label{eq:sys_CH1_weak} \\
&\int_\Gamma  \mu \,q \, ds - \int_\Gamma \frac{1}{\epsilon} f_0'(c) \,q \, ds - \int_\Gamma \epsilon \gradG c \, \gradG q \, ds = 0, \label{eq:sys_CH2_weak}
\end{align}
for all $ (v,q) \in H^1(\Gamma) \times H^1(\Gamma)$.

\subsection{The Navier--Stokes--Cahn--Hilliard problem}

Let us now consider the case of phase separation occurring together with lateral flow. 
The classical phase-field model for the flow of two immiscible, incompressible, and Newtonian fluids with the same density
is the so-called \emph{Model H} \cite{RevModPhys.49.435}. 
To be able to account for non-matching densities, here
we focus on a thermodynamically consistent generalization of Model H first presented in \cite{Palzhanov2021}.

In order to state the NSCH model in \cite{Palzhanov2021}, let $m_i$ be the mass of component $i$ and $m$ is the total mass.
The density of the mixture can be expressed as
$\rho= \frac{m}{S} = \frac{m_1}{S_1} \frac{S_1}{S} + \frac{m_2}{S_2} \frac{S_2}{S}$. Thus, $\rho= \rho(c) = \rho_1 c+ \rho_2 (1-c)$,
where densities $ \rho_1,  \rho_2 > 0$ are given constants. Similarly, for the dynamic viscosity of the mixture we can write
$\eta=\eta(c)=\eta_1 c+\eta_2(1-c)$,
where $\eta_1 >0$ and $\eta_2 >0$ are the constant dynamic viscosities of the two species. 
Then, the model in \cite{Palzhanov2021} reads:
\begin{align}
&  \rho\partial_t \bu + \rho(\nabla_\Gamma\bu)\bu - \bP\divG(2\eta E_s(\bu))+ \nabla_\Gamma p =  -\sigma_\gamma c  \nablaG \mu \cl
& \quad\quad+ {M \theta(\nabla_\Gamma(\theta\bu)\,)\gradG \mu} \label{grache-1m} \\
&\divG \bu  =0 , \label{gracke-2}\\
& \partial_t c +\divG(c\bu)-  \divG \left(M \gradG \mu \right)  = 0  \label{gracke-3}, \\
& \mu = \frac{1}{\epsilon} f_0' - \epsilon \laplG c,  \label{gracke-4}
\end{align}
on $\Gamma \times (0, T]$.
Here, $\bu$ is the surface averaged tangential velocity $\bu = c \bu_1 + (1 - c) \bu_2$, $p$ is pressure, $\sigma_\gamma$ is line tension, 
and $\theta^2 = \frac{d\rho}{dc}$. All other variables and parameters
are the same as defined in Sec.~\ref{sec:CH}. Without loss of generality we let $\rho_1\ge\rho_2$. Then the model \eqref{grache-1m}--\eqref{gracke-4} assumes that $\rho$ is a smooth monotonic function of $c$, i.e.~$\frac{d\rho}{dc}\ge 0$. 
Notice that since $\Gamma$ is stationary all terms in \eqref{grache-1m}--\eqref{gracke-4} are tangential. 
The energy balance delivered by the model reads \cite{Palzhanov2021}:
\begin{multline}\label{energy}
	\frac{d}{dt}
	\int_\Gamma\left(
	\frac{\rho}2\vert\bu\vert^2 
+\sigma_\gamma\left(\frac{1}{\epsilon}f_0 + \frac{\epsilon}2 \vert\nabla_\Gamma c\vert^2
\right)
\right) ds \\
	+ \int_\Gamma 2\eta \vert E_s(\bu)\vert^2 ds + \int_\Gamma \sigma_\gamma M \vert\nabla_\Gamma \mu\vert^2 ds= 0.
\end{multline}

The only difference between model \eqref{grache-1m}--\eqref{gracke-4} and Model H is the last term in eq.~\eqref{grache-1m}, which
can be interpreted as an additional momentum flux due to diffusion of
the components driven by the gradient of the chemical potential. This term vanishes 
for matching densities since $\theta = \sqrt{\frac{d\rho}{dc}} = 0$, thereby recovering Model H. For other thermodynamic consistent extensions of 
Model H that involve a generic smooth $\rho(c)$ (no monotonicity assumption), the reader is referred to
\cite{abels2016weak,abels2019existence}.

For the purpose of writing the weak formulation of problem \eqref{grache-1m}--\eqref{gracke-4}, 
we define the spaces
\begin{equation}   \label{defVT}
 \bV_T\coloneqq \{\, \bu \in H^1(\Gamma)^3~\vert~ \bu\cdot \bn =0\,\},\quad E\coloneqq \{\, \bu \in \bV_T~\vert~ E_s(\bu)=\mathbf{0}\,\}.
\end{equation}
Moreover, we define the Hilbert space $\bV_T^0$ as an orthogonal complement of  $E$ in $\bV_T$
(hence $\bV_T^0 \sim \bV_T/E$) and $L_0^2(\Gamma)\coloneqq\{\, p \in L^2(\Gamma)~\vert~ \int_\Gamma p\,ds=0\,\}$.

The weak formulation of the surface NSCH problem \eqref{grache-1m}-\eqref{gracke-4} reads:
Find $(\bu,p,c,\mu) \in \bV_T \times L_0^2(\Gamma) \times H^1(\Gamma) \times H^1(\Gamma)$ 
such that
\begin{align}
&\int_\Gamma \left( \rho \partial_t \bu\cdot\bv + \rho(\nabla_\Gamma\bu)\bu \cdot \bv + 2\eta  E_s(\bu):E_s(\bv)\right) \, ds - \int_\Gamma p\,\divG \bv \, ds =  \cl
& \quad \quad -\int_\Gamma \sigma_\gamma  c  \nablaG \mu \cdot\bv \, ds  + \int_\Gamma M (\nabla_\Gamma(\theta\bu)) (\gradG \mu) \cdot(\theta\bv)  \, ds, \label{gracke-w1} \\
& \int_\Gamma q\,\divG \bu \, ds = 0, \label{gracke-w2} \\
&\int_\Gamma \partial_t c \,v \, ds - \int_\Gamma c \bu \cdot  \gradG v \, ds + \int_\Gamma M \gradG \mu \cdot \gradG v \, ds = 0, \label{gracke-w3} \\
&\int_\Gamma  \mu \,g \, ds =  \int_\Gamma \frac{1}{\epsilon} f_0'(c) \,g \, ds + \int_\Gamma \epsilon \gradG c \cdot \gradG g \, ds, \label{gracke-w4}
\end{align}
for all $ (\bv,q,v,g) \in \bV_T \times L^2(\Gamma) \times H^1(\Gamma) \times H^1(\Gamma)$. More details on the derivation of
\eqref{gracke-w1}-\eqref{gracke-w4} can be found in \cite{Palzhanov2021}.

\section{Numerical method}\label{sec:method}

For the numerical solution of the problems presented in Sec.~\ref{sec:model}, we apply the trace finite element method (TraceFEM)~\cite{ORG09,olshanskii2017trace}. 
TraceFEM relies on a tessellation of a 3D bulk computational domain
$\Omega$ ($\Gamma\subset\Omega$ holds) into shape-regular tetrahedra untangled to the position of $\Gamma$.

Surface $\Gamma$ is defined as the zero level set of a function $\phi$ (where $\phi$ is at least Lipschitz continuous), i.e.~$\Gamma=\{\bx\in\Omega\,:\, \phi(\bx)=0\}$,
such that $\vert \nabla\phi \vert \ge c_0>0$ in a 3D neighborhood $U(\Gamma)$ of the surface.
The vector field $\bn=\nabla\phi/\vert \nabla\phi \vert$ is normal on $\Gamma$ and  defines  quasi-normal directions in $U(\Gamma)$.
Let $\mT_h$ be  the collection of all tetrahedra such that $\overline{\Omega}=\cup_{T\in\mT_h}\overline{T}$, with 
$h$ denoting the characteristic tetrahedra size.
The subset of tetrahedra that have a nonzero intersection with $\Gamma$ is denoted by $\mT_h^\Gamma$.
We allow for local refinement of the grid towards $\Gamma$.
The domain formed by all tetrahedra in $\mT_h^\Gamma$ is denoted by $\OGamma$.

In order to state the fully discretized CH and NSCH problems, we introduce some finite element spaces.
Let $V_h^k$ denote the bulk (volumetric) finite element space of continuous functions that are polynomials
of degree $k$ on each $T\in\mT_h^{\Gamma}$:
\[
V_h^k=\{v\in C(\OGamma)\,:\, v\in P_k(T)~\text{for any}~T\in\mT_h^{\Gamma}\}.
\]
The traces of functions from $V_h^1$ on $\Gamma$ will be used to approximate the surface fraction and the chemical potential.
Our bulk velocity and pressure finite element spaces are the Taylor--Hood elements on $\OGamma$:
\begin{equation}\label{TH}
 \bV_h = (V_h^{2})^3, \quad Q_h = V_h^1 \cap L^2_0(\Gamma). 
 \end{equation}
Higher order approximations are possible (see, e.g., 
\cite{Palzhanov2021,grande2018analysis}) but will not be considered here. 

For the purpose of numerical integration, we approximate $\Gamma$ by a ``discrete'' surface $\Gamma_h$
so that integrals over $\Gamma_h$ can be computed accurately and efficiently.
 For first order finite elements, a straightforward  polygonal approximation of $\Gamma$
 ensures that the geometric approximation error is consistent with the finite element
 interpolation error. See, e.g., \cite{ORG09}. 
 
 Next, we introduce two finite element bilinear forms that are common to the discrete versions of both
 NS and NSCH problems:
\begin{align}
& a_\mu(\mu, v) \coloneqq \int_{\Gamma} M \nabla_{\Gamma} \mu \cdot \nabla_{\Gamma}v\,ds +  \tau_\mu \int_{\OGamma} (\bn\cdot\nabla \mu) (\bn\cdot\nabla v)\,d\bx, \label{eq:a_mu} \\
& a_c(c, g) \coloneqq  \epsilon\int_{\Gamma} \gradG c \cdot \gradG g \, ds + \tau_c \int_{\OGamma} ( \bn\cdot\nabla c) (\bn\cdot\nabla g) \,d\bx. \label{eq:a_c}
\end{align}
Forms \eqref{eq:a_mu}--\eqref{eq:a_c} are well defined for $\mu,v, c, g \in H^1(\OGamma)$.
The volumetric terms in \eqref{eq:a_mu} and \eqref{eq:a_c} are there to recover algebraic stability
as possible small cuts of tetrahedra from $\mT_h^\Gamma$ by $\Gamma$ may lead to poorly conditioned algebraic systems.
Notice that these terms are consistent up to geometric errors related to the approximation of $\Gamma$ by $\Gamma_h$ and $\bn$ by $\bn_h$.
We set the stabilization parameters as follows: 
$
\tau_\mu=h,\quad \tau_c=\epsilon\,h^{-1}.
$

For the time discretization, let $\Delta t=\frac{T}{N}$ be a time step. 
At time instance $t^n=n\Delta t$, $\zeta^n$ denotes the approximation of generic
variable $\zeta(t^n, \bx)$. To approximate the time derivatives in problems \eqref{eq:sys_CH1_weak}-\eqref{eq:sys_CH2_weak}
and \eqref{gracke-w1}-\eqref{gracke-w4}, we use the backward differentiation formula of order 1 (BDF1):
 \begin{equation*}
\left[\zeta\right]_t^{n} =\frac{\zeta^{n}- \zeta^{n-1}}{\Delta t}.
\end{equation*}

Once fully discretized, CH problem \eqref{eq:sys_CH1_weak}-\eqref{eq:sys_CH2_weak} reads: 
Given $c^n_h\in V_h^1$, find $(c^{n+1}_h, \mu^{n+1}_h) \in V^1 _h \times V^1 _h$
such that:
\begin{align}
&\left(\left[c_h\right]_t^{n+1}, v_h\right) + a_\mu(\mu_h^{n+1}, v_h) = 0,  \cl
&\left(  \mu_h^{n+1} - \frac{\gamma_c\Delta t}{\epsilon}\left[c_h\right]_{t}^{n+1}  - \frac{1}{\epsilon} f'_0(c_h^{n}),\,g_h\right)
- a_c(c_h^{n+1}, g_h) = 0,  \label{eq:CH_FE2_0}
  \end{align}
for all  $(v_h, g_h) \in V^1 _h \times V^1 _h$. Following~\cite{Shen_Yang2010}, the second term in \eqref{eq:CH_FE2_0} stabilizes the explicit treatment  of non-linear part of the free energy variation.

Next, we turn to the NSCH problem \eqref{gracke-w1}-\eqref{gracke-w4}. For its numerical solution, 
we adopt a decoupled linear finite element method introduced in \cite{Palzhanov2021}. In order to described such method, 
we need to introduce some additional forms related to the Navier--Stokes part of the problem
and the decomposition of a vector field on $\Gamma$ into its tangential and normal components:
$\bu=\overline{\bu}+(\bu\cdot\bn)\bn$. The additional forms are defined as follows:
\begin{align}
& a(\eta; \bu,\bv) \coloneqq \int_\Gamma 2\eta E_s( \overline{\bu}):  E_s( \overline{\bv})\, ds+\tau \int_{\Gamma}(\bn\cdot\bu) (\bn\cdot\bv) \, ds \cl
&\hspace{2cm} + \beta_u \int_{\OGamma} [(\bn\cdot\nabla) \bu] \cdot [(\bn\cdot\nabla) \bv] + \widehat{\gamma}\,\int_\Gamma \divG \bu \, \divG \bv\,ds  \, d\bx, \label{eq:a} \\
& c(\rho; \bw, \bu,\bv)\coloneqq \int_\Gamma\rho\bv^T(\nabla_\Gamma\overline{\bu})\bw\, ds 
 +\frac12 \int_\Gamma\widehat{\rho}(\divG \overline{\bw})\overline{\bu}\cdot\overline{\bv}\, ds, \label{eq:c} \\
&b(\bu,q) =  \int_\Gamma \bu\cdot\nabla_\Gamma q \, ds, \label{eq:b} \\
&s(p,q) \coloneqq \beta_p  \int_{\OGamma} \nabla p\cdot\nabla q \, d\bx, \label{eq:s}
\end{align}
where $\widehat{\gamma}$ is the grad-div stabilization parameter \cite{olshanskii2002low} (set equal to 1) and $\widehat{\rho}=\rho-\frac{d\rho}{d\,c}c$.
Forms \eqref{eq:a}--\eqref{eq:s} are well defined for $p,q \in H^1(\OGamma)\cap H^1(\Gamma)$, $\bu,\bv,\bw \in H^1(\OGamma)^3\cap H^1(\Gamma)^3$. In \eqref{eq:a}, $\tau>0$ is a penalty parameter to enforce the tangential constraint (i.e., condition 
$\bu_h \cdot \bn = 0$ on $\Gamma$ for $\bu_h \in  \bV_h$),
while $\beta_u\ge0$ in \eqref{eq:a} and $\beta_p\ge0$ in \eqref{eq:s} are stabilization parameters
to deal with possible small cuts. They are set according to \cite{Jankuhn2020}:
$
\tau=h^{-2},\quad \beta_p=h, \quad \beta_u=h^{-1}.
$

The decoupled finite element method from \cite{Palzhanov2021} requires the solution of 
one linear problem Chan--Hilliard type system (step 1) and one  linearized Navier--Stokes system (step 2)
per time step $t^{n+1}$, thereby ensuring low computational costs. This scheme, which is 
provably stable under relatively mild restrictions \cite{Palzhanov2021}, reads:
\begin{itemize}
\item[-] \underline{Step 1}: Given $\bu^n_h\in \bV_h$ and $c^n_h\in V_h^1$, find $(c^{n+1}_h, \mu^{n+1}_h) \in V^1 _h \times V^1 _h$
such that:
\begin{align}
&\left(\left[c_h\right]_t^{n+1}, v_h\right) - \left(\bu^{n}_h  c^{n+1}_h, \nablaG v_h\right) + a_\mu(\mu_h^{n+1}, v_h) = 0,  \label{eq:CH_FE1} \\
&\left(  \mu_h^{n+1} - \frac{\gamma_c\Delta t}{\epsilon}\left[c_h\right]_{t}^{n+1}  - \frac{1}{\epsilon} f'_0(c_h^{n}),\,g_h\right)
- a_c(c_h^{n+1}, g_h) = 0,  \label{eq:CH_FE2}
  \end{align}
for all  $(v_h, g_h) \in V^1 _h \times V^1 _h$.
\item[-] \underline{Step 2}: Set 
$\theta^{n+1} =\sqrt{\frac{d\rho}{d c}(c^{n+1}_h)}$.
Find $(\bu_h^{n+1}, p_h^{n+1}) \in \bV_h \times Q_h$ such that
 \begin{align}
&  (\rho^n\left[\overline{\bu}_h\right]_t^{n+1},\bv_h)+ c(\rho^{n+1}; \bu^{n}_h, \bu^{n+1}_h,\bv_h) + a(\eta^{n+1};\bu_h^{n+1},\bv_h)  + b(\bv_h,p_h^{n+1}) \cl
& \quad \quad = -(\sigma_\gamma c^{n+1}_h \nablaG \mu^{n+1}_h, \bv_h) + M \left((\nabla_\Gamma(\theta^{n+1}\overline{\bu}_h^{n+1}))\gradG \mu^{n+1}_h, \theta^{n+1}\bv_h\right) \cl
& \hspace{1.1cm} + (\bbf_h^{n+1}, \bv_h)  \label{NSEh1} \\
 & b(\bu_h^{n+1},q_h)-s(p_h^{n+1},q_h)  = 0  \label{NSEh2}
 \end{align}
for all  $(\bv_h,q_h) \in \bV_h \times Q_h$.
\end{itemize}

\section{Numerical results}\label{sec:num_res}

We present a series of numerical results aimed at understanding the difference in the evolution of 
phases when modeled by the Cahn--Hilliard or Navier--Stokes--Cahn--Hilliard equations posed on a closed smooth surface.
For the latter model, we experiment with different settings for the physical parameters. 

We start by comparing the numerical results produced by the two models
on a sphere in Sec.~\ref{sec:sphere}. Then, we consider an asymmetric torus in Sec.~\ref{sec:torus}
to see the effects of a different geometry on the evolution of phases. For all the simulations we
set $\epsilon = 0.02$ and $D = 0.02$ in \eqref{mobility}. In order to model an initially homogenous mix of components, 
the initial are fraction $c_0$ is defined  as a realization of Bernoulli random variable~$c_\text{rand} \sim \text{Bernoulli}(a)$
with mean value $a$, i.e. we set:
\begin{equation}\label{raftIC}
	c_0 \coloneqq c_\text{rand}(\bx)\quad\text{for active mesh nodes $\bx$}.
\end{equation}
We set $a=0.5$ for the 50\%-50\% composition (meaning that 50\% of the surface is covered by one phase and 
the remaining 50\% by the other phase) and $a=0.3$ for the 30\%-70\% composition. 
The other physical parameters will be specified for each case.
We run all the simulations till $T = 100$ and an adaptive time stepping
technique~\cite{gomez2008isogeometric}.

For all the simulations, we will study the evolution of the discrete Lyapunov energy:
\begin{align}
E^L_h(c_h) = \int_{\Gamma_h} f(c_h) ds = \int_{\Gamma_h} \left( \frac{1}{\epsilon}f_0(c_h) + \frac{1}{2} \epsilon \vert \gradG c_h \vert^2 \right)ds \label{eq:Lyapunov_E}
\end{align}
and we will visually compare the evolutions of phases. In addition, for the
NSCH model we will compare the flow in a qualitative way. 

\subsection{Phase separation on a sphere} \label{sec:sphere}

The surface of the sphere is appealing for its simplicity and for its relevance in practical applications. In fact,
lipid vesicles used as drug carriers have a spherical shape \cite{zhiliakov2021experimental}. We characterize $\Gamma$
as the zero level set of function $\phi(\bx) = \|\bx\|_2 -1$ and we embed it in an outer cubic domain $\Omega=[-5/3,5/3]^3$.

We experimented with different meshes to find one with an appropriate level of refinement
for the given value of $\epsilon$. 
The initial triangulation $\mathcal{T}_{h_\ell}$ of $\Omega$ we considered consists of eight sub-cubes,
where each of the sub-cubes is further subdivided into six tetrahedra. 
We applied several level of refinement $\ell\in\Bbb{N}$, with associated mesh size $h_\ell= \frac{10/3}{2^{\ell+1}}$.
Each mesh also features a refinement towards the surface.  
We found that $\ell = 5$ is a good compromise between accuracy and computational cost. See also \cite{Palzhanov2021}.
Thus, the results reported in this section refer to the mesh with $\ell = 5$. We note that for the NSCH model 
such mesh has 225822 active degrees of freedom (193086 for $\bu_h$
and 10912 for $p_h$, $c_h$, and $\mu_h$). 

\subsubsection{Variable line tension}\label{sec:variable_gamma}

In this section, we focus on composition 50\%-50\%. One initial condition \eqref{raftIC} is generated and 
used to compare phase separation given by the CH model and the NSCH model
with variable line tension. 
We assign density $\rho_1 = 3$ and viscosity $\eta_1 =0.1$ to species 1, 
while species 2 has $\rho_2 = 1$ and $\eta_2 =0.008$. We consider a low value of line tension
$\sigma_\gamma =  0.004$ and one high value $\sigma_\gamma =  0.4$.

\begin{figure}[htb]
\centering
\includegraphics[width=0.48\textwidth]{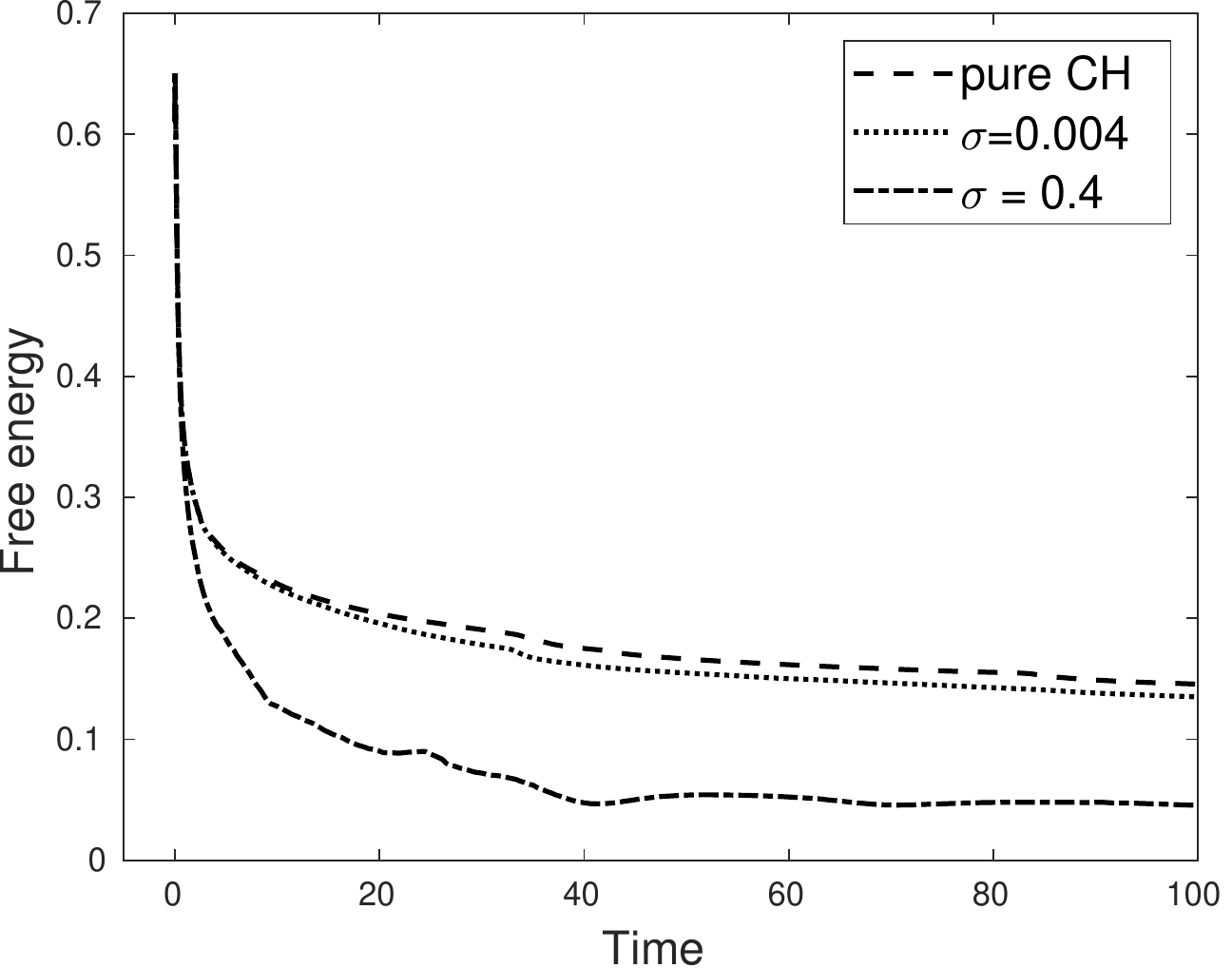}
	\caption{Discrete Lyapunov energy \eqref{eq:Lyapunov_E} given by the CH model, NSCH model with $\sigma_\gamma = 0.004$, and NSCH model with $\sigma_\gamma = 0.4$}
	\label{fig:free_energy_sigma}		
\end{figure}

Fig.~\ref{fig:free_energy_sigma}	 shows the Lyapunov energy \eqref{eq:Lyapunov_E}
over time computed by the CH model, NSCH model with low and high line tension. 
We observe that when switching form the CH model to NSCH model with low line tension
the Lyapunov energy decay is slightly faster. It becomes substantially faster when the value
 of $\sigma_\gamma$ is increased, which can be expected from the energy balance~\eqref{energy}
 since the last dissipative term scales with  $\sigma_\gamma$.

\begin{figure}
\begin{center}
\hskip .7cm
\begin{overpic}[width=.15\textwidth,grid=false]{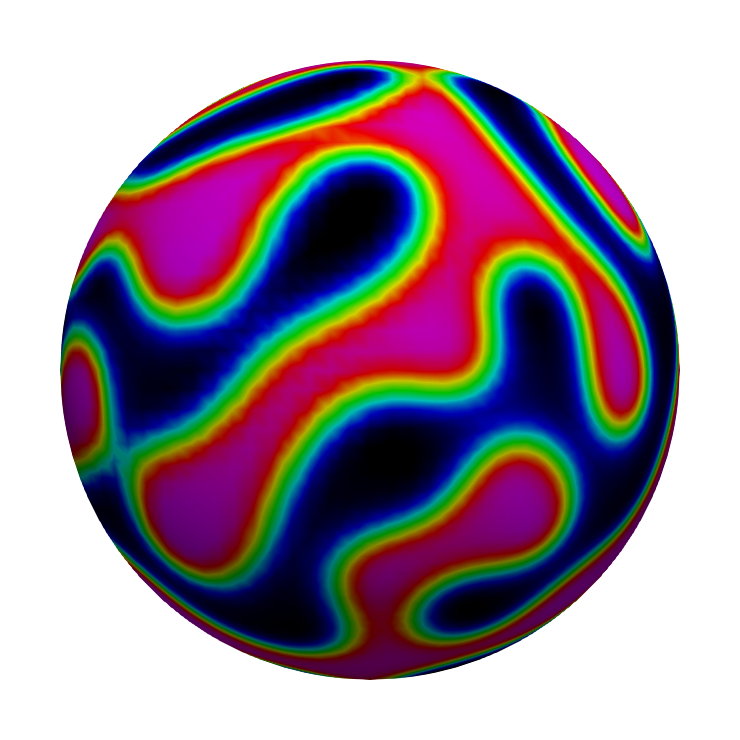}
\put(30,98){\small{$t = 2$}}
\put(-30,50){\small{CH}}
\end{overpic}
\begin{overpic}[width=.15\textwidth,grid=false]{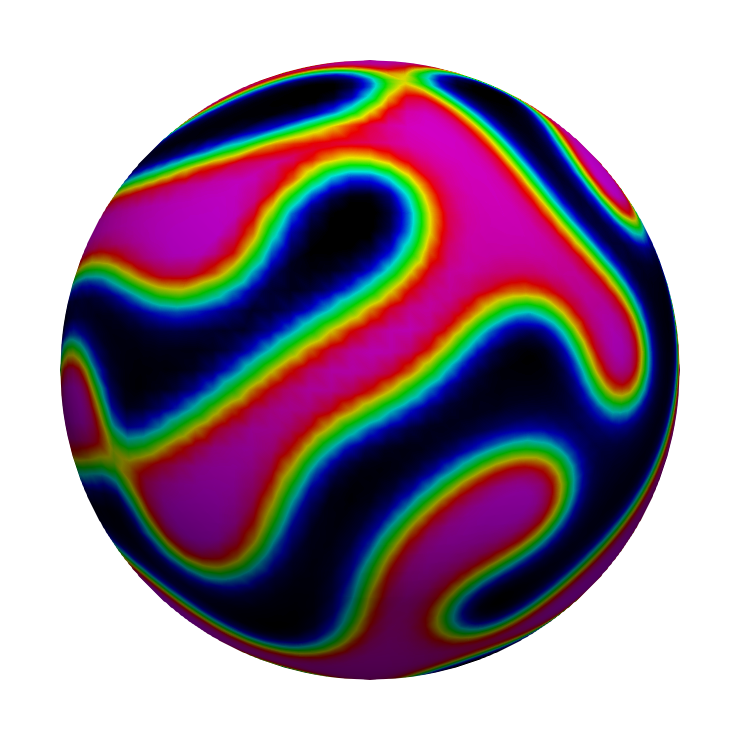}
\put(30,98){\small{$t = 5$}}
\end{overpic}
\begin{overpic}[width=.15\textwidth,grid=false]{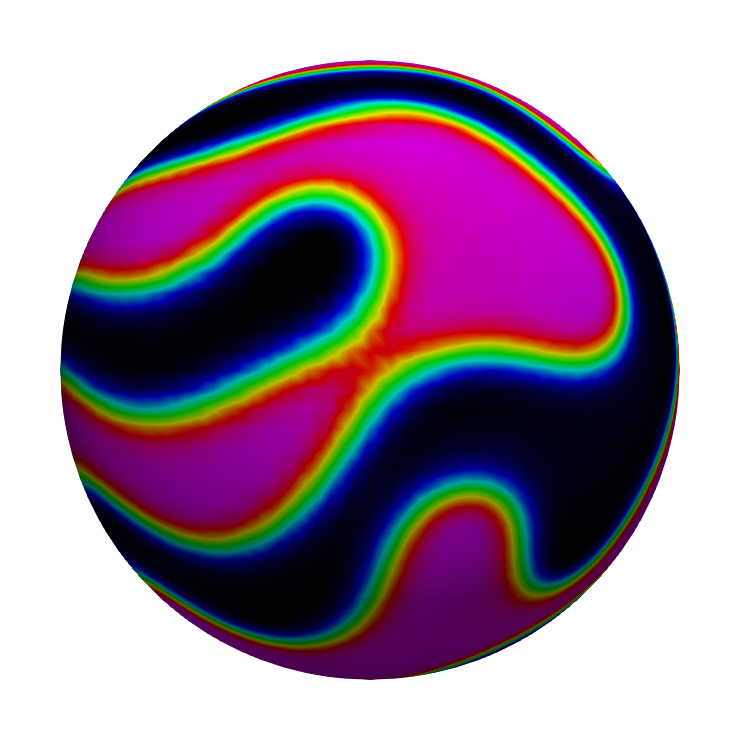}
\put(30,98){\small{$t = 15$}} 
\end{overpic}
\begin{overpic}[width=.15\textwidth,grid=false]{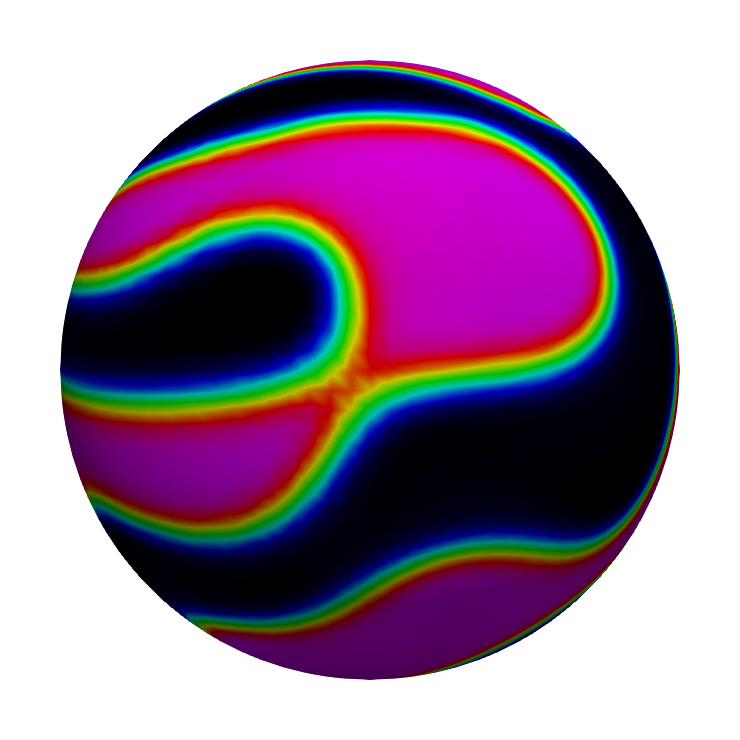}
\put(28,98){\small{$t = 30$}}
\end{overpic}
\begin{overpic}[width=.15\textwidth,grid=false]{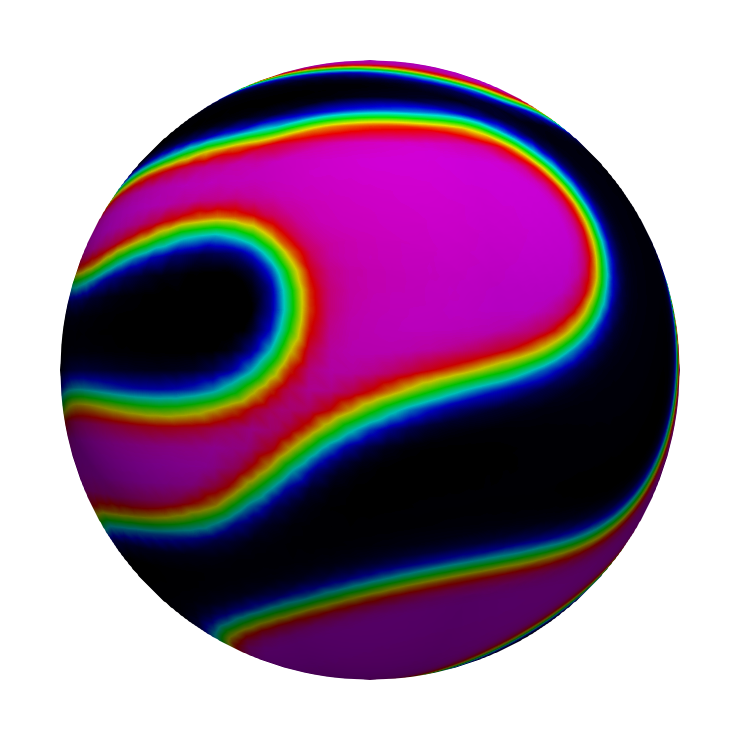}
\put(28,98){\small{$t = 50$}}
\end{overpic}
\begin{overpic}[width=.15\textwidth,grid=false]{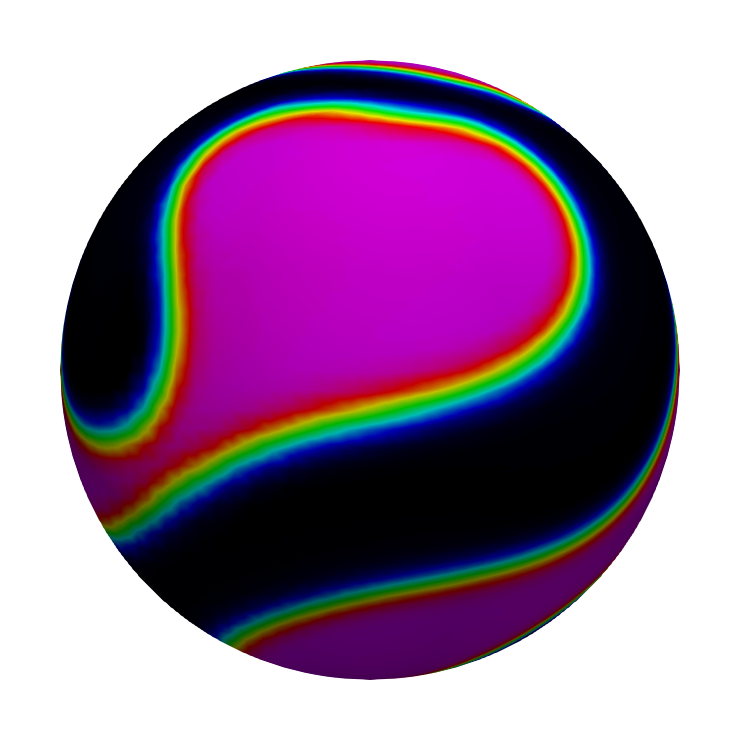}
\put(26,95){\small{$t = 100$}}
\end{overpic}
\\
\hskip .7cm
\begin{overpic}[width=.15\textwidth,grid=false]{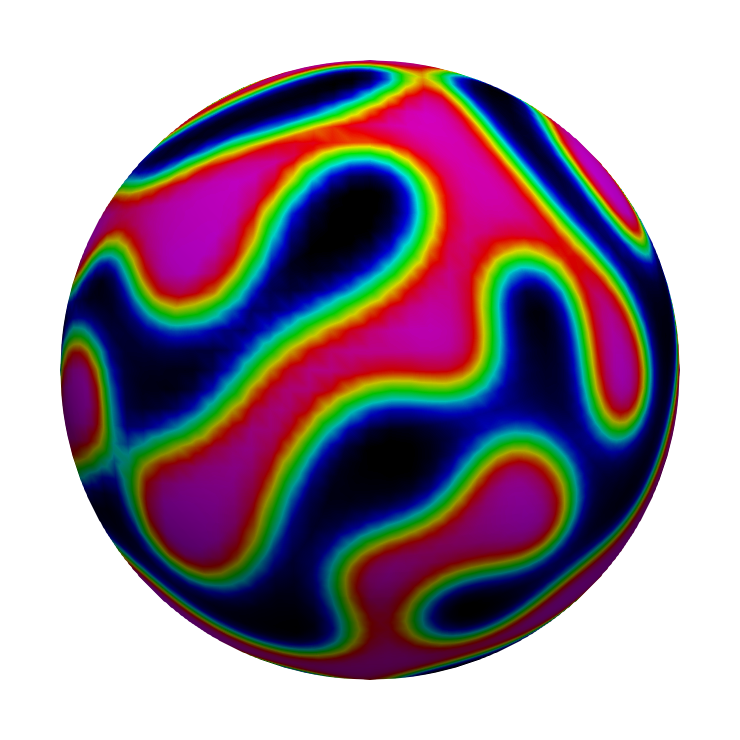}
\put(-50,50){\small{NSCH}}
\put(-50,33){\small{low $\sigma_\gamma$}}
\end{overpic}
\begin{overpic}[width=.15\textwidth,grid=false]{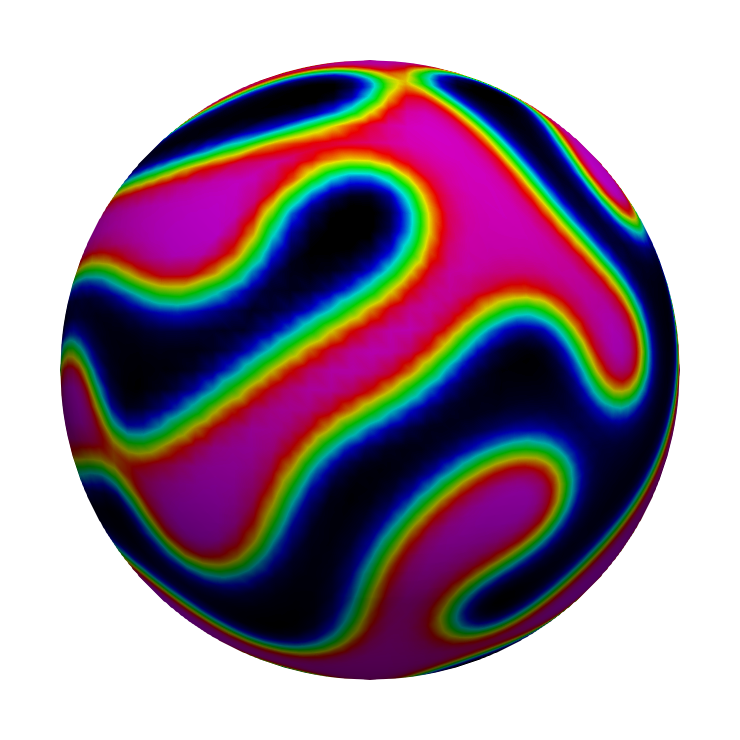}
\end{overpic}
\begin{overpic}[width=.15\textwidth,grid=false]{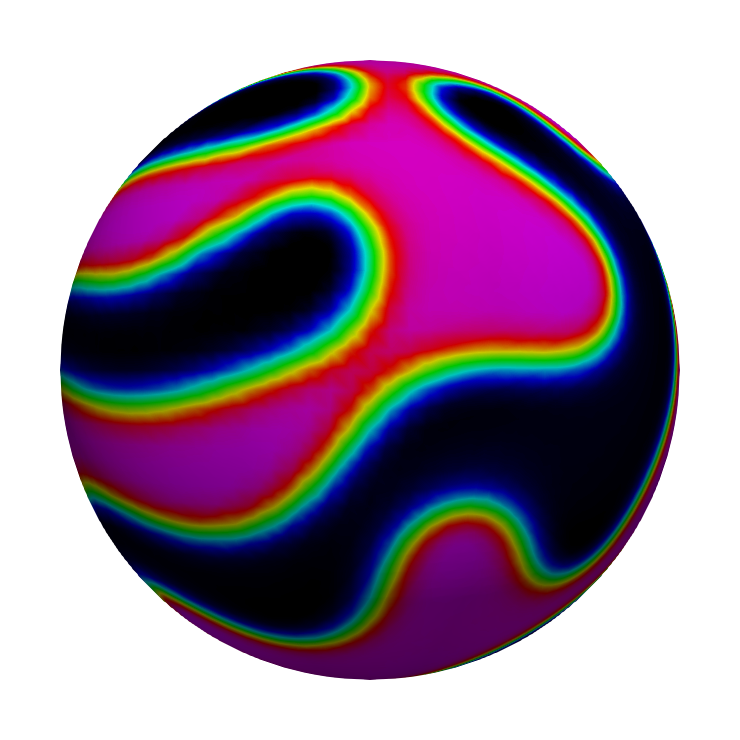}
\end{overpic}
\begin{overpic}[width=.15\textwidth,grid=false]{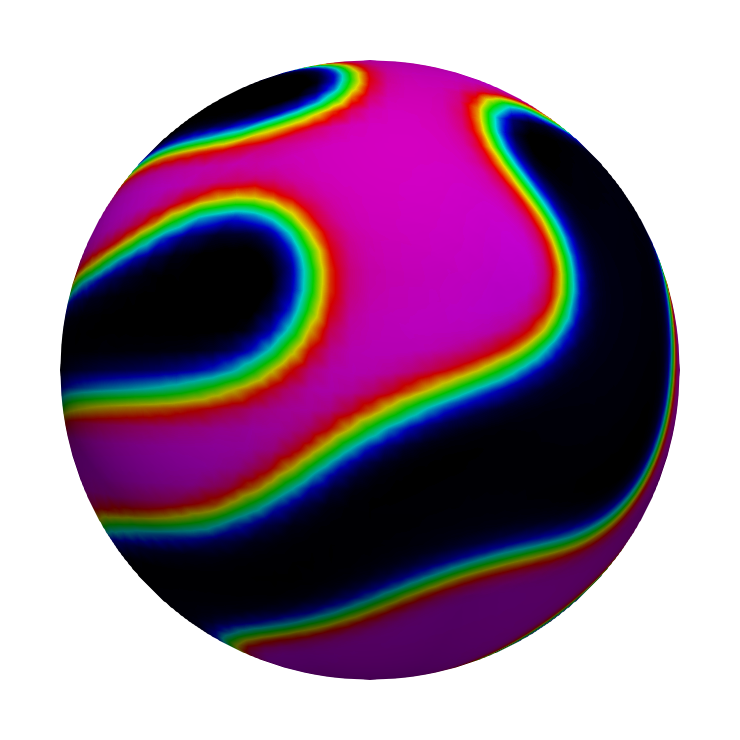}
\end{overpic}
\begin{overpic}[width=.15\textwidth,grid=false]{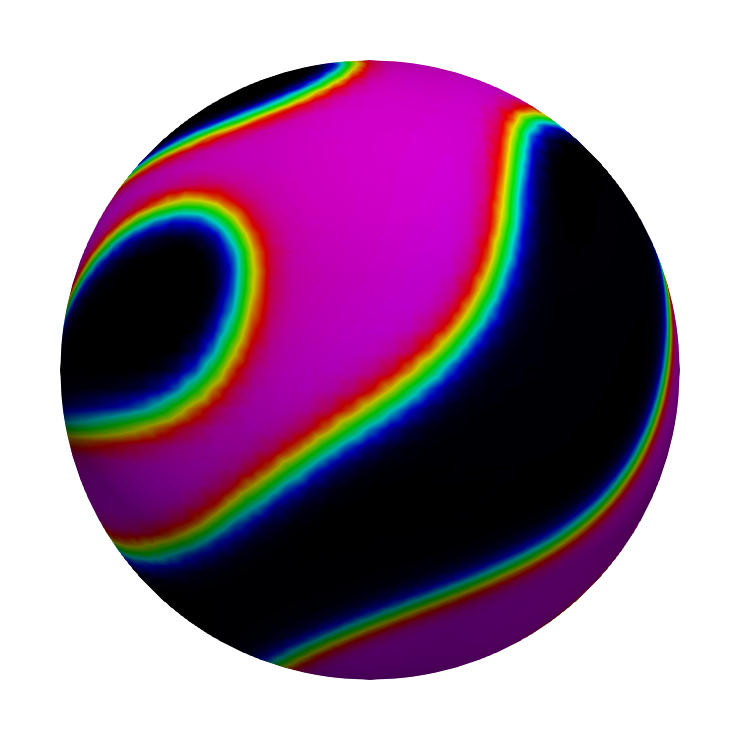}
\end{overpic}
\begin{overpic}[width=.15\textwidth,grid=false]{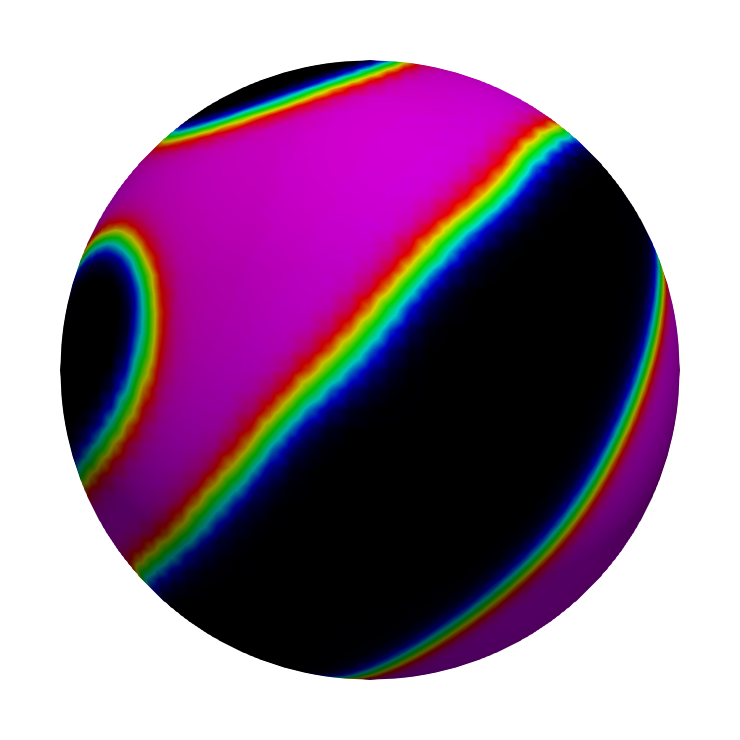}
\end{overpic}
\\
\hskip .7cm
\begin{overpic}[width=.15\textwidth,grid=false]{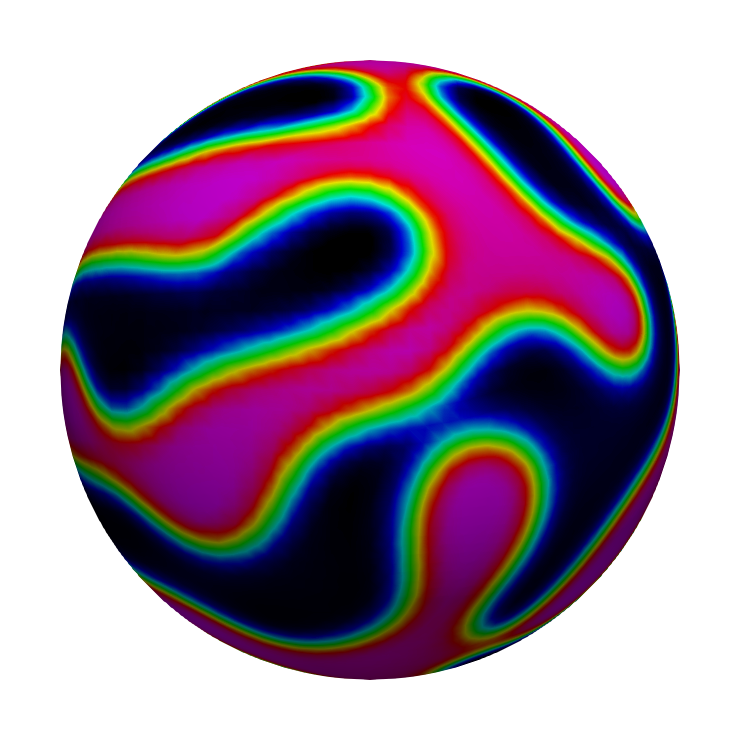}
\put(-50,50){\small{NSCH}}
\put(-50,33){\small{high $\sigma_\gamma$}}
\end{overpic}
\begin{overpic}[width=.15\textwidth,grid=false]{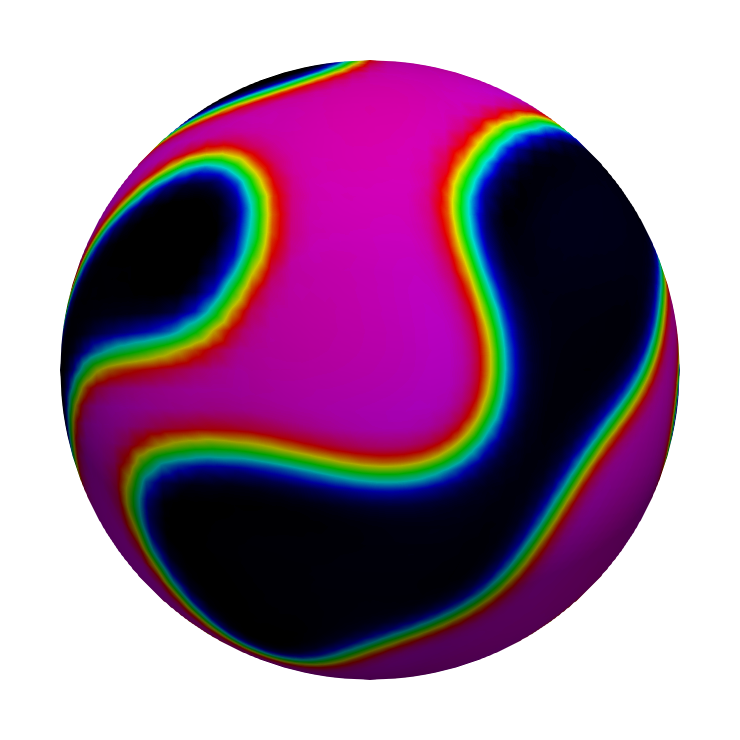}
\end{overpic}
\begin{overpic}[width=.15\textwidth,grid=false]{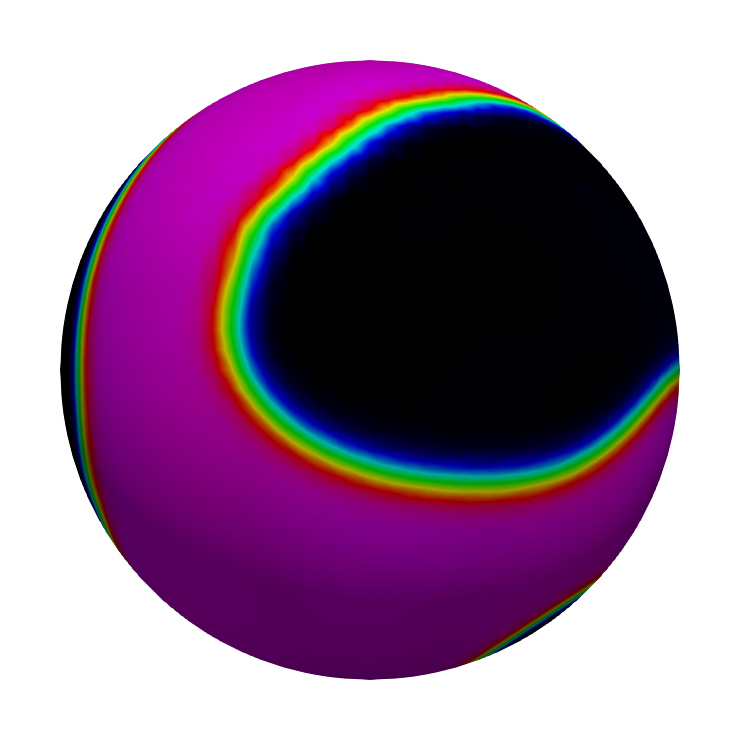}
\end{overpic}
\begin{overpic}[width=.15\textwidth,grid=false]{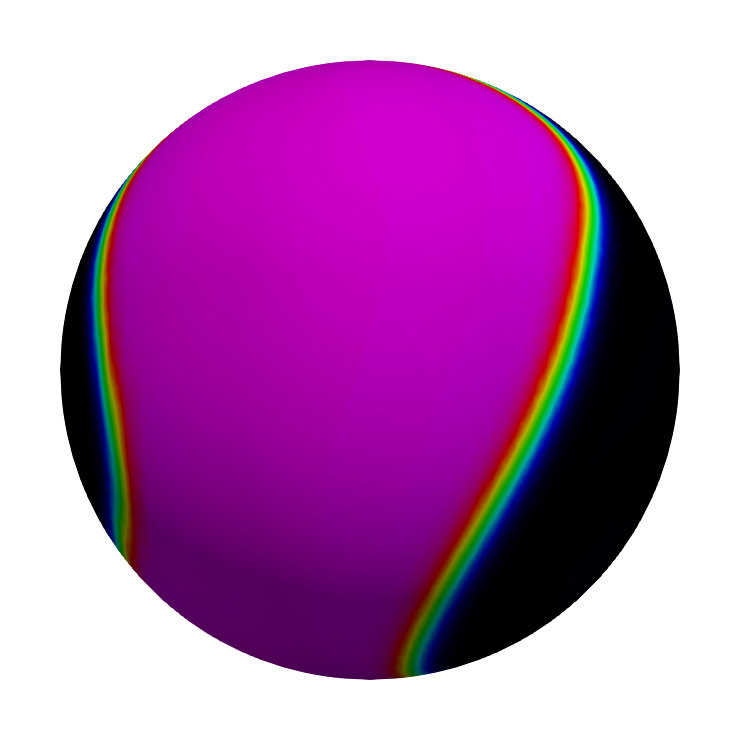}
\end{overpic}
\begin{overpic}[width=.15\textwidth,grid=false]{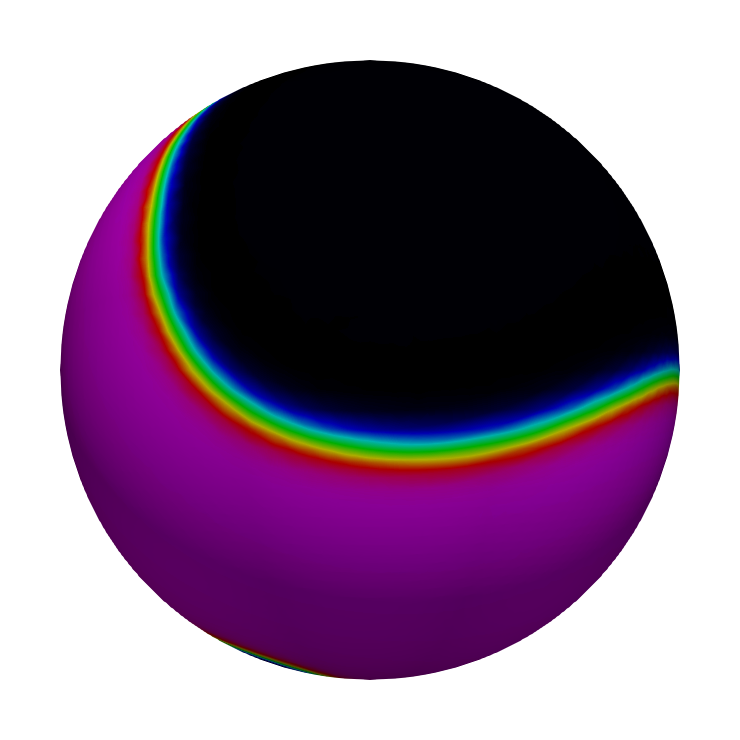}
\end{overpic}
\begin{overpic}[width=.15\textwidth,grid=false]{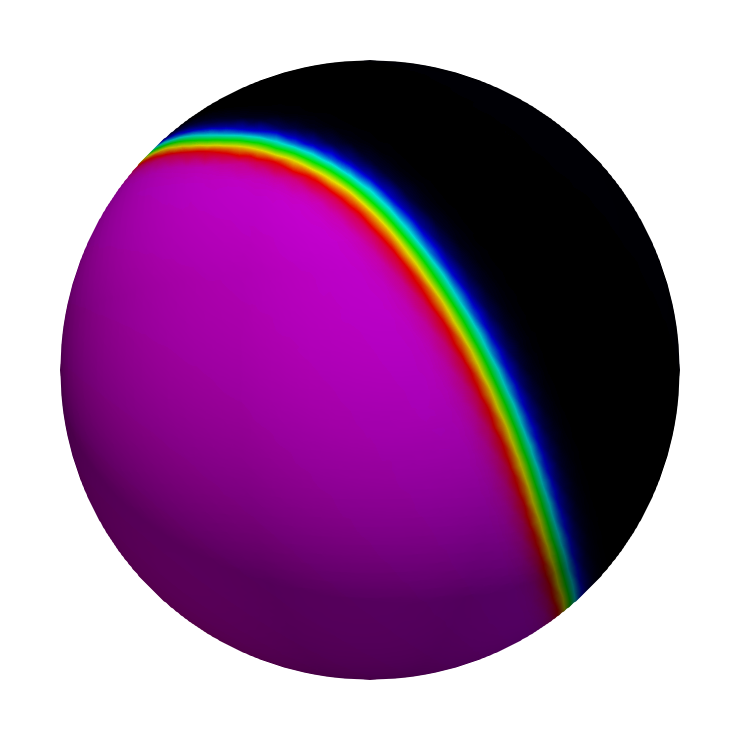}
\end{overpic}
\\
\vskip .2cm
\begin{overpic}[width=0.5\textwidth,grid=false,tics=10]{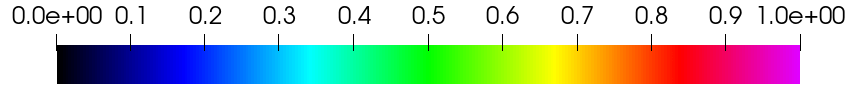}
\end{overpic}
\end{center}
\caption{Phase separation given by the CH model (top), NSCH model with $\sigma_\gamma = 0.004$
(center), and NSCH model with $\sigma_\gamma = 0.4$ (bottom). }\label{fig:sigma}
\end{figure}

These differences are reflected in the evolution of phases displayed in 
Fig.~\ref{fig:sigma}. The evolution of the surface fraction 
does not vary significantly when going from the CH model to 
the NSCH model with $\sigma_\gamma = 0.004$,
although some differences can be noticed from $t =30$ on. 
Changing to $\sigma_\gamma = 0.4$ produces
more evident differences, starting already from $t =5$. Moreover, by comparing the
center and bottom rows in Fig.~\ref{fig:sigma} it is clear that a larger value of $\sigma_\gamma$
accelerates the transition towards a steady state, i.e. one large black domain and one large
pink domain separated by a minimal length interface. 

\begin{figure}
\begin{center}
\hskip .7cm
\begin{overpic}[width=.15\textwidth,grid=false]{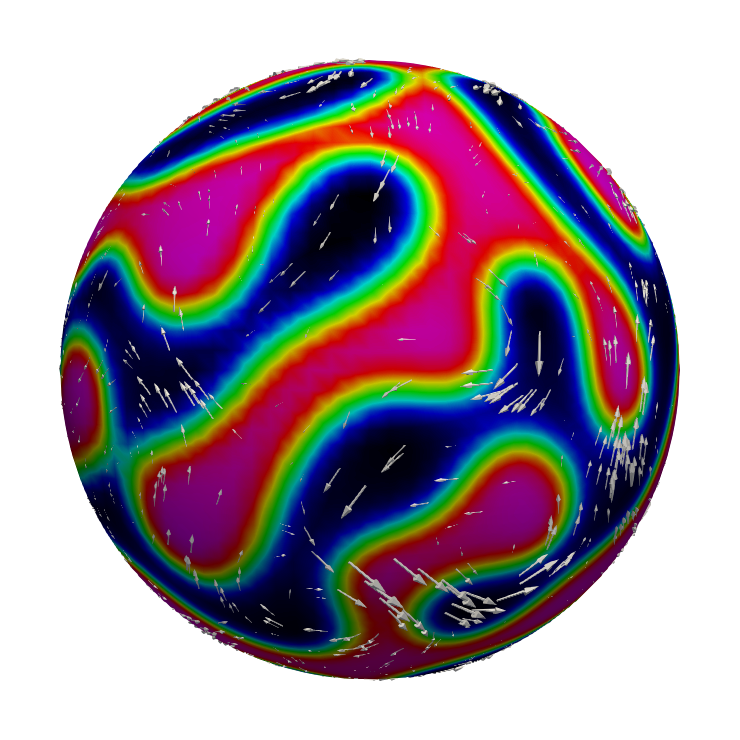}
\put(-53,45){\small{low $\sigma_\gamma$}}
\put(30,98){\small{$t = 2$}}
\end{overpic}
\begin{overpic}[width=.15\textwidth,grid=false]{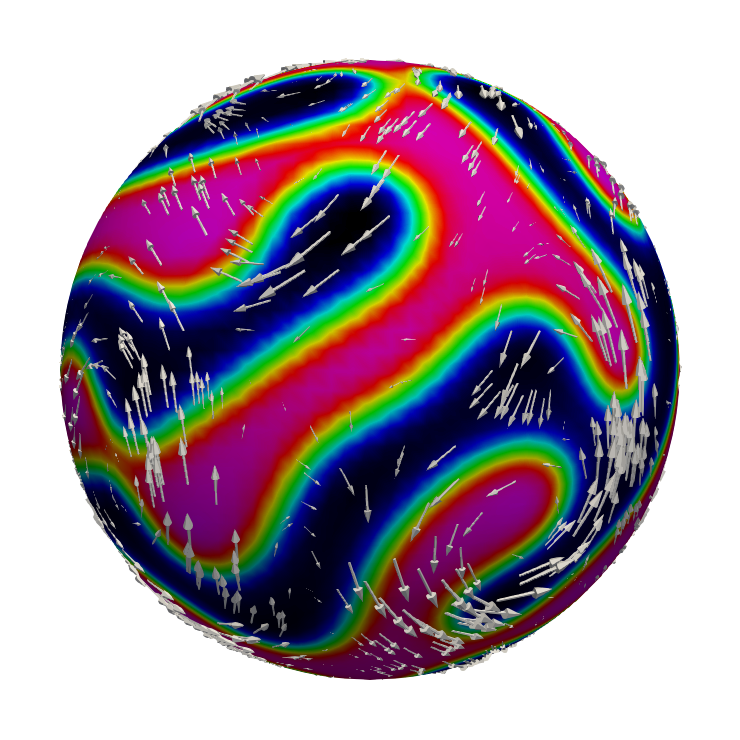}
\put(30,98){\small{$t = 5$}}
\end{overpic}
\begin{overpic}[width=.15\textwidth,grid=false]{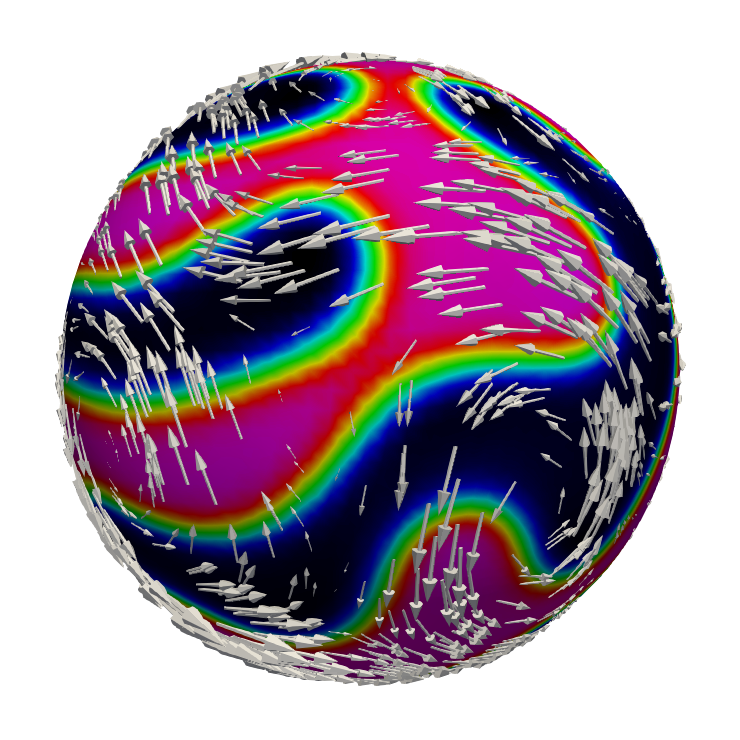}
\put(28,98){\small{$t = 15$}}
\end{overpic}
\begin{overpic}[width=.15\textwidth,grid=false]{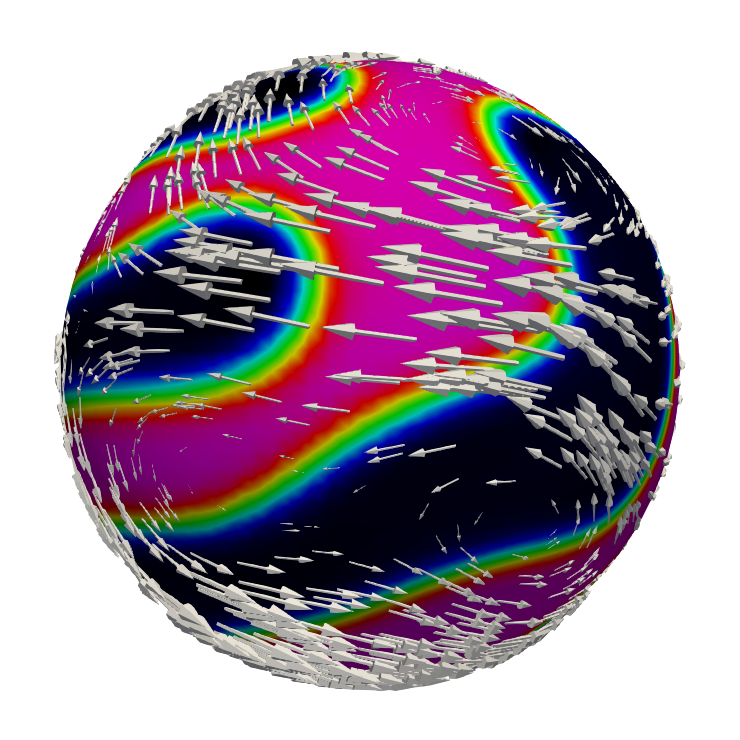}
\put(28,98){\small{$t = 30$}}
\end{overpic}
\begin{overpic}[width=.15\textwidth,grid=false]{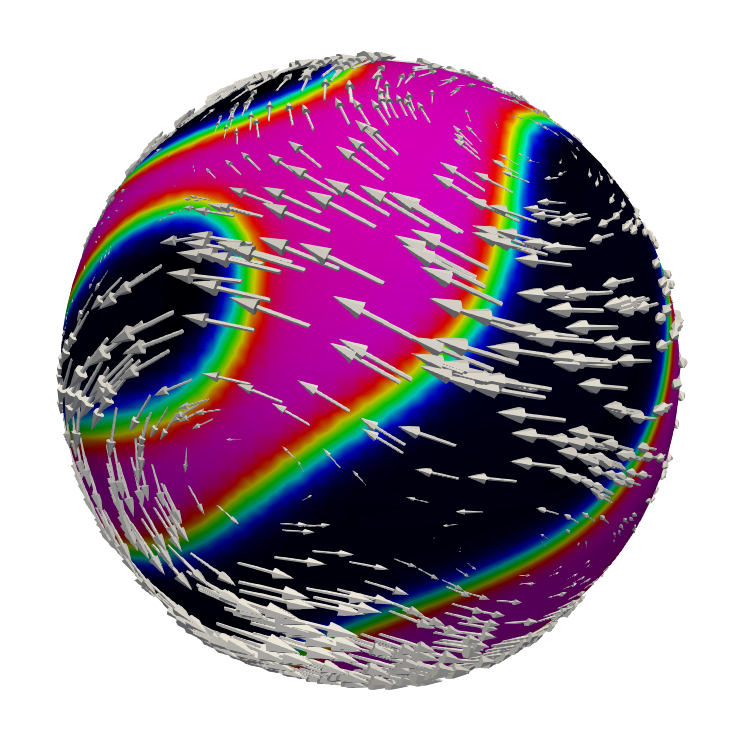}
\put(28,98){\small{$t = 50$}}
\end{overpic}
\begin{overpic}[width=.15\textwidth,grid=false]{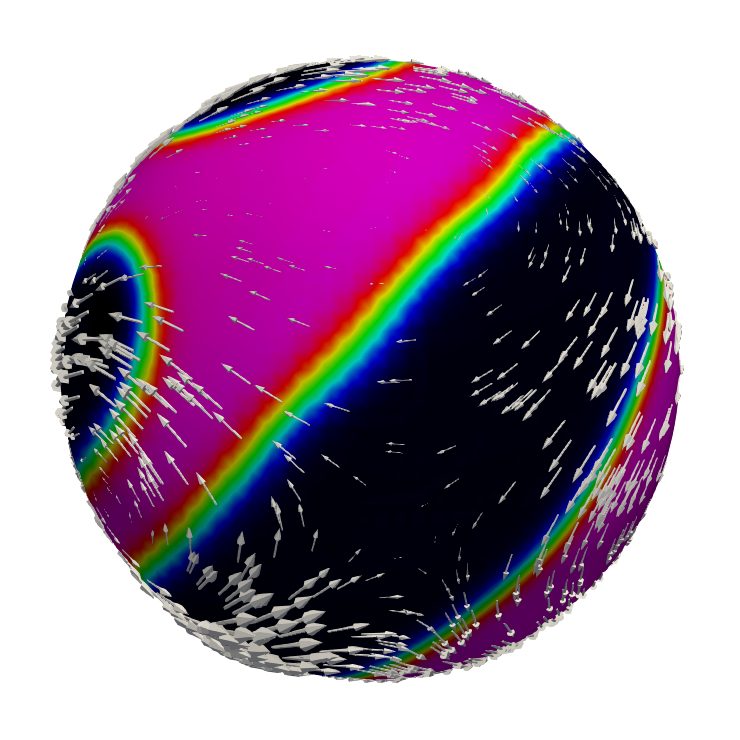}
\put(26,98){\small{$t = 100$}}
\end{overpic}
\\
\hskip .7cm
\begin{overpic}[width=.15\textwidth,grid=false]{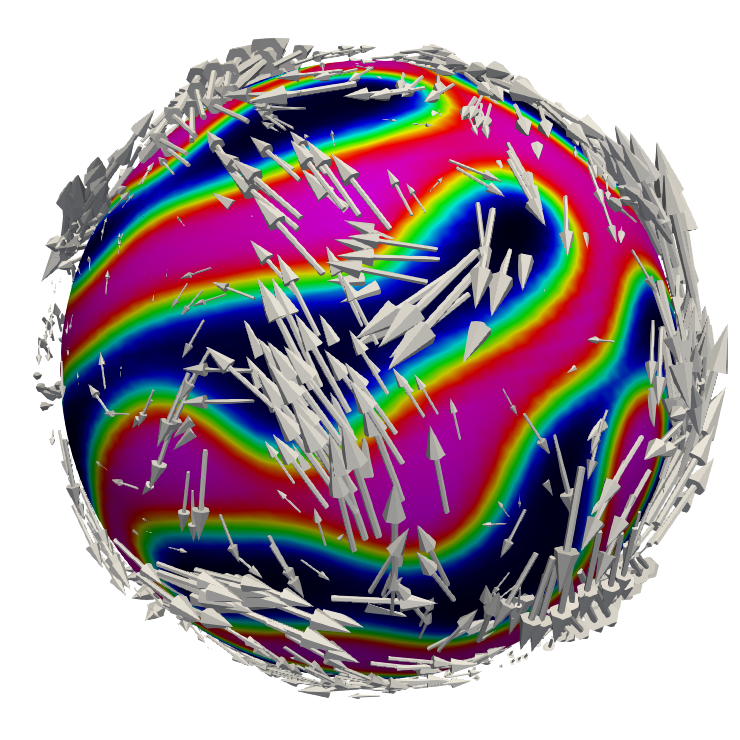}
\put(-53,45){\small{high $\sigma_\gamma$}}
\end{overpic}
\begin{overpic}[width=.15\textwidth,grid=false]{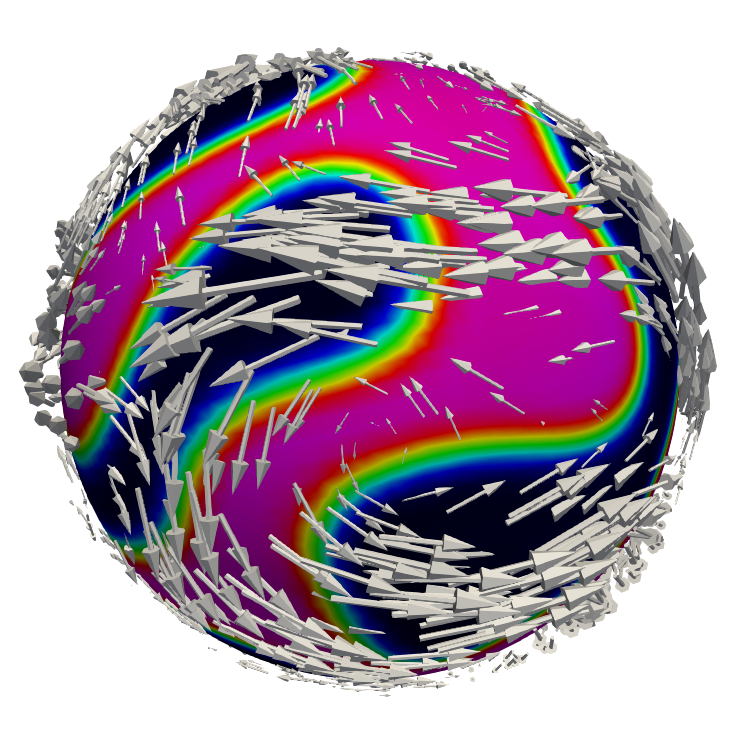}
\end{overpic}
\begin{overpic}[width=.15\textwidth,grid=false]{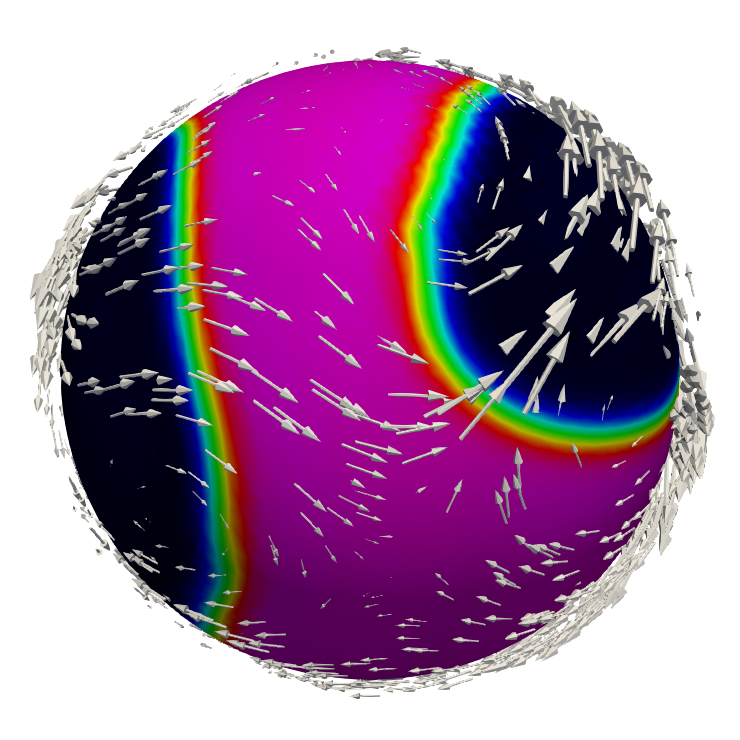}
\end{overpic}
\begin{overpic}[width=.15\textwidth,grid=false]{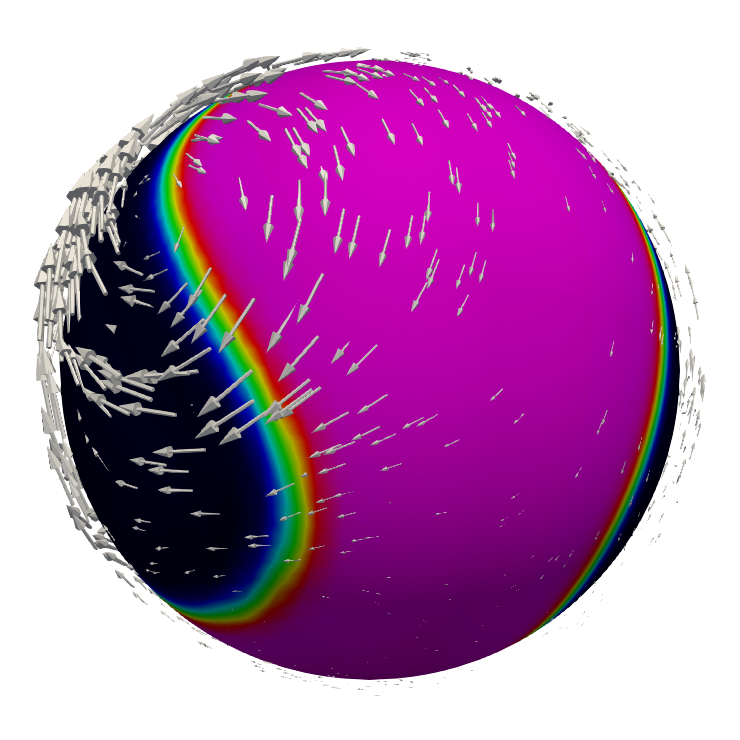}
\end{overpic}
\begin{overpic}[width=.15\textwidth,grid=false]{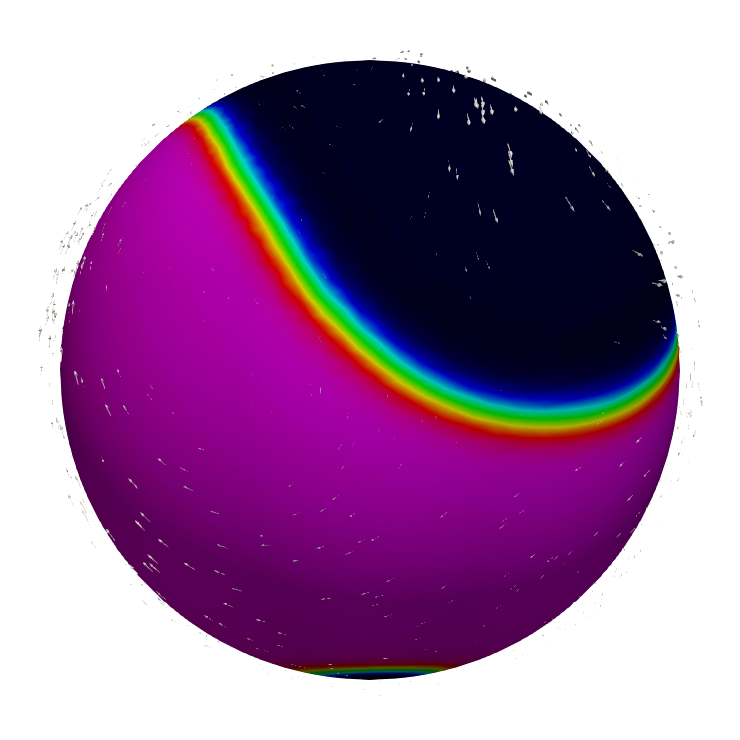}
\end{overpic}
\begin{overpic}[width=.15\textwidth,grid=false]{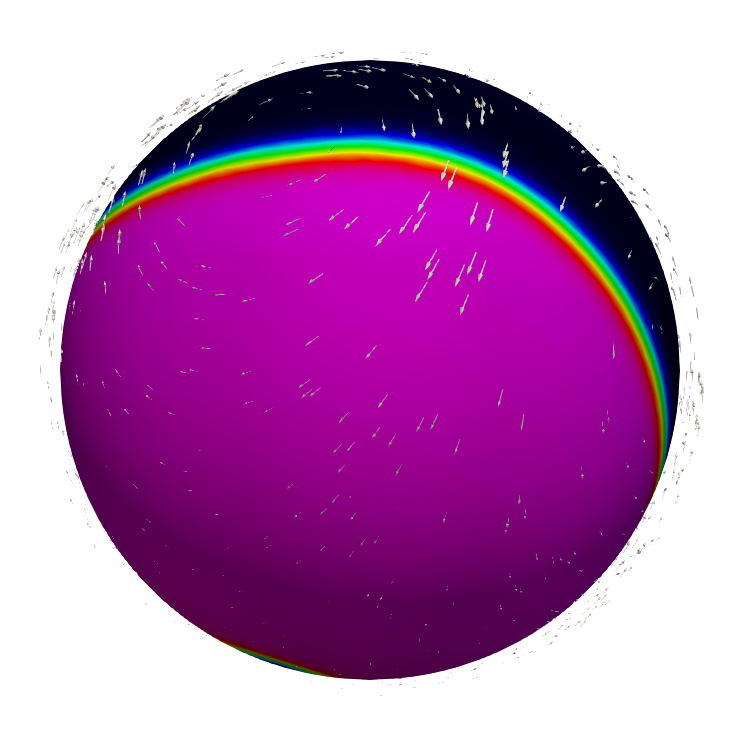}
\end{overpic}
\end{center}
\caption{Velocity vectors superimposed to the surface fraction for $\sigma_\gamma = 0.004$ (top)
and $\sigma_\gamma = 0.4$ (bottom). For visualization purposes, the  
velocity vectors are magnified by a factor 40 in the top row and 2 in the bottom row. 
}\label{fig:flow_sigma}
\end{figure}

Fig.~\ref{fig:flow_sigma} displays the velocity vectors superimposed to the surface fraction
for the bottom two cases in Fig.~\ref{fig:sigma}. Since for visualization purposes the arrows have 
been magnified with different factors, the velocity vectors cannot be compared across rows. 
In the NSCH model in \eqref{grache-1m}-\eqref{gracke-4}, the fluid flow is purely driven by the coupling
with the phase separation process. For $\sigma_\gamma = 0.4$ (bottom row in Fig.~\ref{fig:flow_sigma}), 
the larger surface tension forces initially produce more significant fluid motion, which however decays faster over time. 
This is due to the fact that the
system evolves more rapidly towards a steady state, as mentioned above. This is not the case for 
$\sigma_\gamma = 0.004$. See Fig.~\ref{fig:flow_sigma}, top row.

\subsubsection{Variable viscosity}\label{sec:variable_visc_sphere}

Now, we set line tension to $\sigma_\gamma = 0.04$ and vary the viscosity for composition
50\%-50\% and 30\%-70\% (which means that 30\% of the surface is in phase 1). 
We consider a high viscosity case ($\eta_1 = 0.01$, $\eta_2 = 0.0008$) and a low viscosity 
case ($\eta_1 = 0.0001$, $\eta_2 = 0.000008$). The densities are set like in Sec.~\ref{sec:variable_gamma}.

\begin{figure}[htb]
\centering
\includegraphics[width=0.48\textwidth]{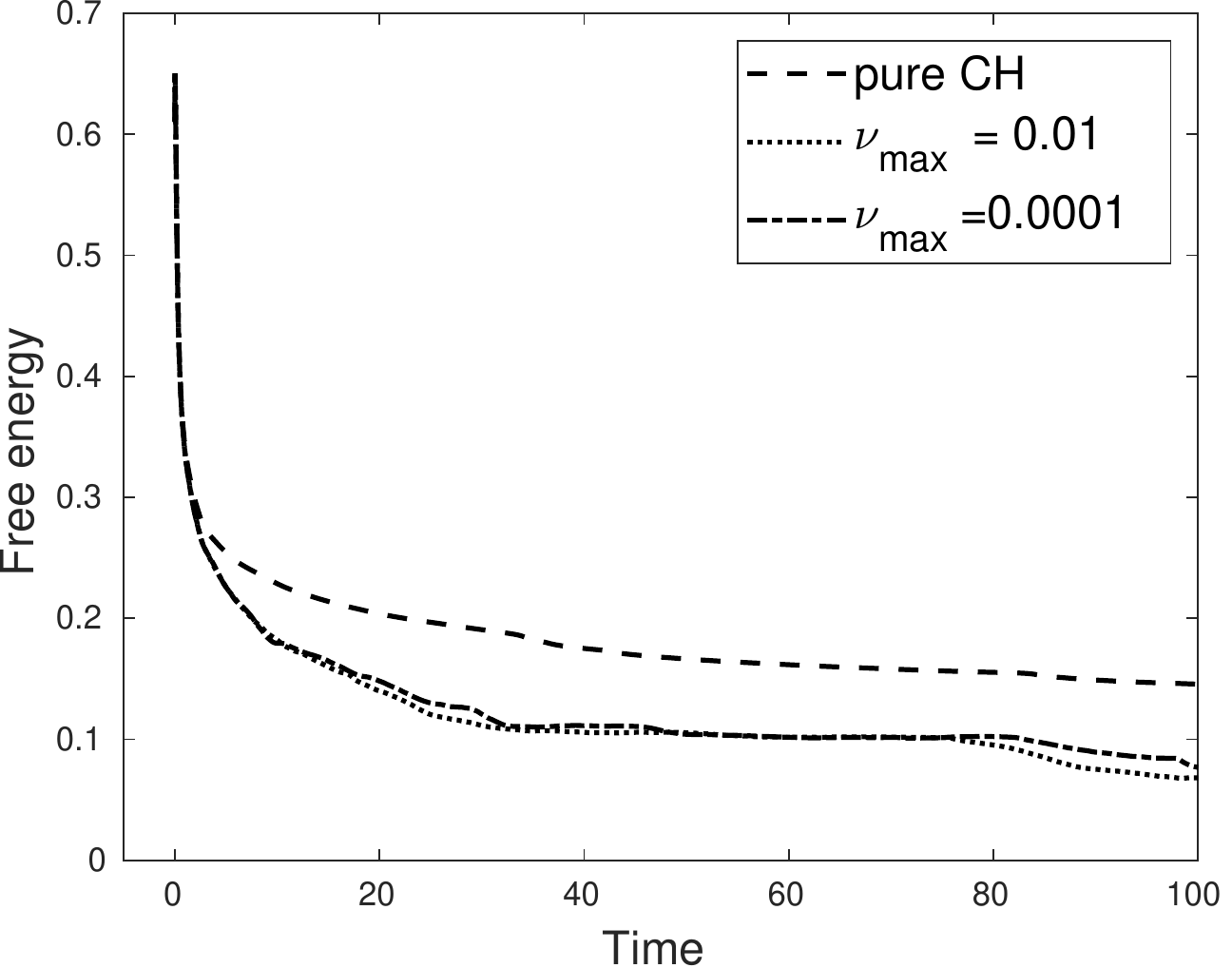}
\includegraphics[width=0.48\textwidth]{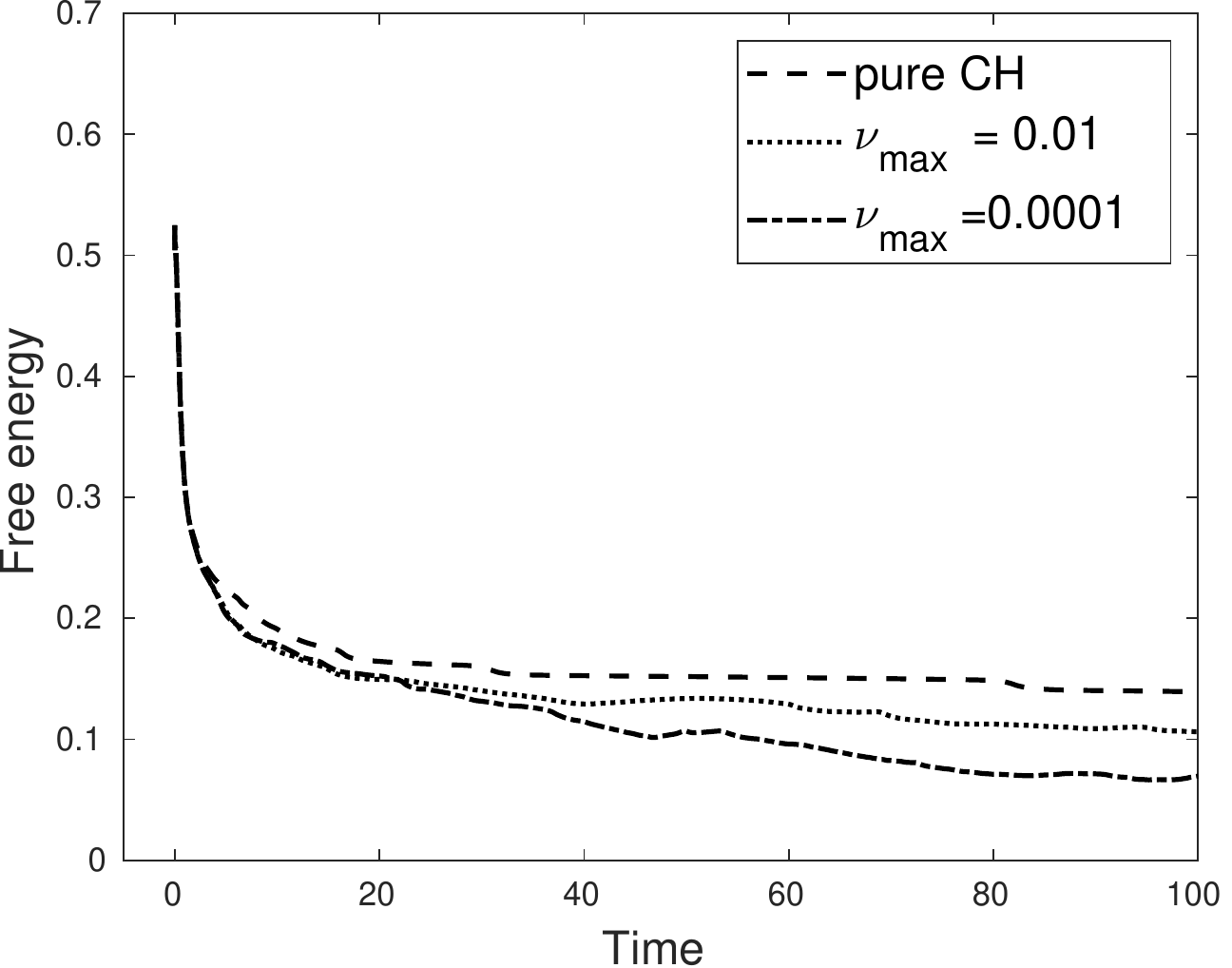}
	\caption{Discrete Lyapunov energy \eqref{eq:Lyapunov_E} given by the CH model, NSCH model with high values of viscosities
	($\eta_1 = 0.01$, $\eta_2 = 0.0008$), and NSCH model with low values of viscosities ($\eta_1 = 0.0001$, $\eta_2 = 0.000008$)
	for composition 50\%-50\% (left) and 30\%-70\% (right).
	}
	\label{fig:free_energy_eta}		
\end{figure}

\begin{figure}[htb]
\centering
\includegraphics[width=0.48\textwidth]{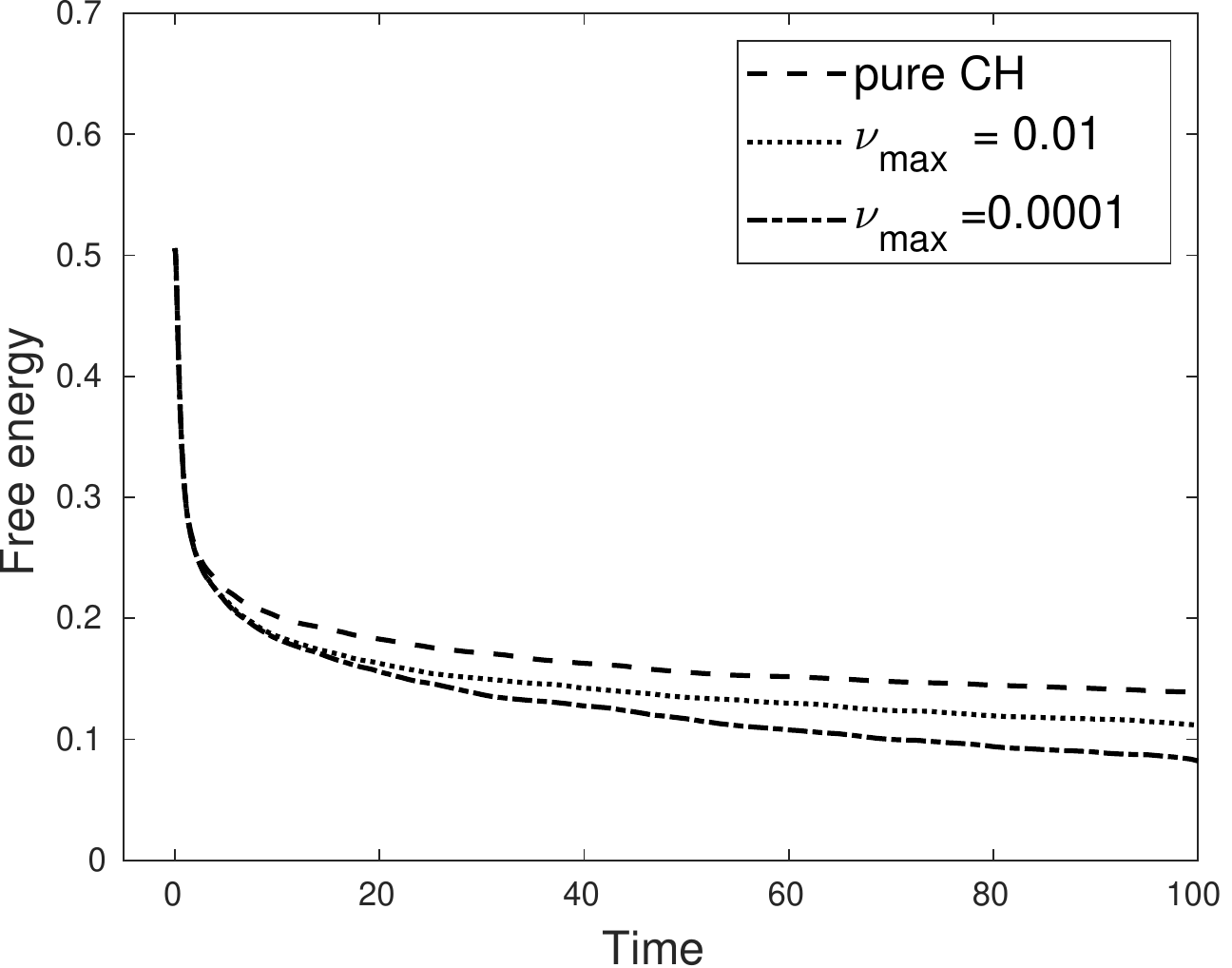}
	\caption{Average discrete Lyapunov energy \eqref{eq:Lyapunov_E} given by he CH model and
	 the NSCH model with high values of viscosities
	($\eta_1 = 0.01$, $\eta_2 = 0.0008$) and low values of viscosities ($\eta_1 = 0.0001$, $\eta_2 = 0.000008$). Average is taken over ten simulations
	for composition 30\%-70\% (right) with different realizations of random initial distribution.
	}
	\label{fig:free_energy_eta_30}		
\end{figure}

Fig.~\ref{fig:free_energy_eta} reports the discrete Lyapunov energy \eqref{eq:Lyapunov_E}
over time computed by the CH model, NSCH model with low and high viscosities for both compositions.
We observe that the presence of surface flow leads to a faster Lyapunov energy decay. In the case of composition
50\%-50\%, we see that switching from high to low values of the viscosity does not produce a significant
change in the energy decay. Instead, for composition 30\%-70\% Fig.~\ref{fig:free_energy_eta} shows that lower values of viscosity 
lead to a faster Lyapunov energy decay than higher values. An explanation of this phenomenon is not obvious from the energy balance~\eqref{energy}, since the first  dissipation term scales with viscosity. To ensure that this was not an accident, 
we repeated the numerical experiment for composition 30\%-70\% ten times with different realizations of initial condition
\eqref{raftIC}. The average Lyapunov energy computed by CH model and NSCH model with low and high viscosities is reported in
Fig.~\ref{fig:free_energy_eta_30}, which confirms the trend. 
It seems plausible to hypothesize the following: 
for lower viscosity the surface tension produces higher speed lateral flows (as illustrated in Fig.~\ref{fig:flow_sigma}), which increases the probability of small rafts coming together and merging, thereby releasing the free energy.        

\begin{figure}
\begin{center}
\hskip .7cm
\begin{overpic}[width=.15\textwidth,grid=false]{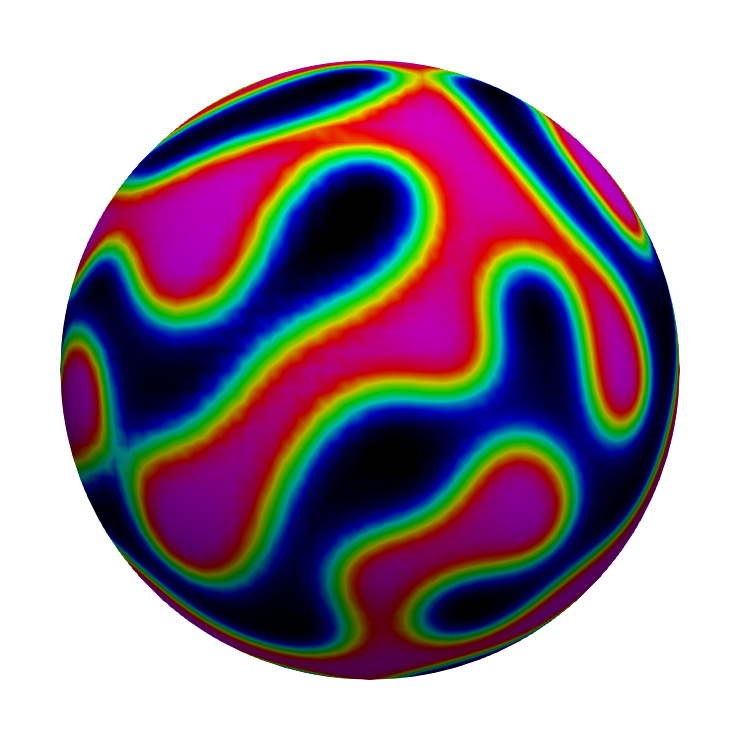}
\put(30,98){\small{$t = 2$}}
\put(-30,50){\small{CH}}
\end{overpic}
\begin{overpic}[width=.15\textwidth,grid=false]{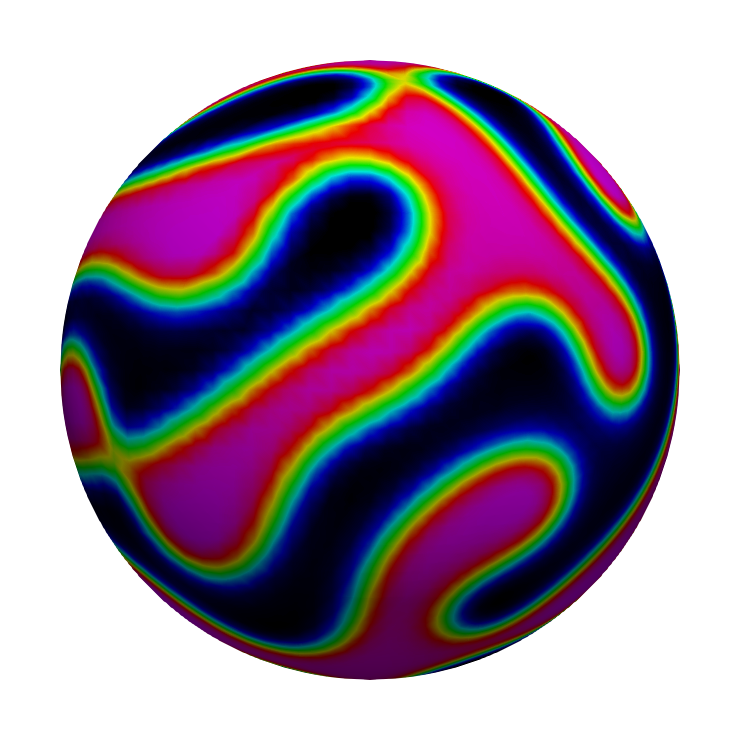}
\put(30,98){\small{$t = 5$}}
\end{overpic}
\begin{overpic}[width=.15\textwidth,grid=false]{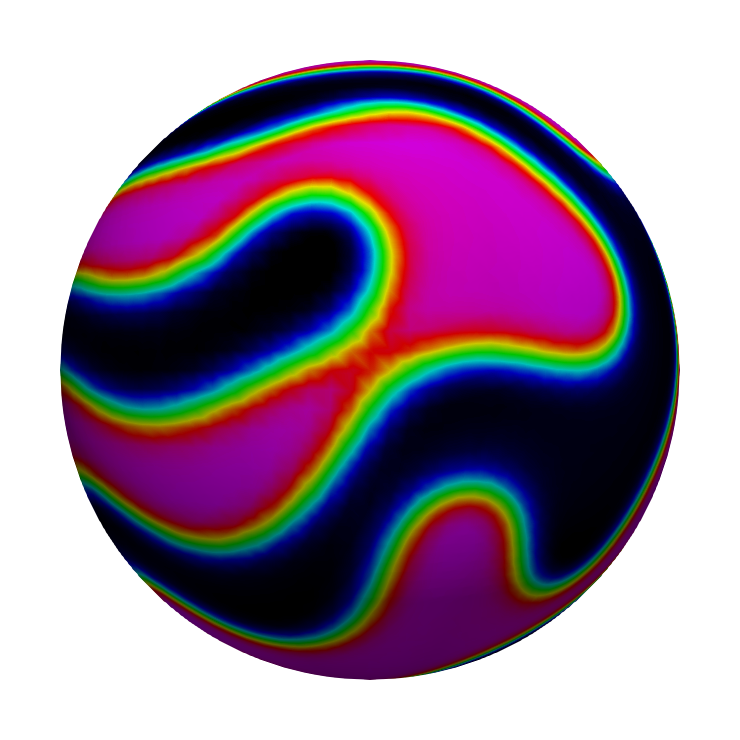}
\put(30,98){\small{$t = 15$}} 
\end{overpic}
\begin{overpic}[width=.15\textwidth,grid=false]{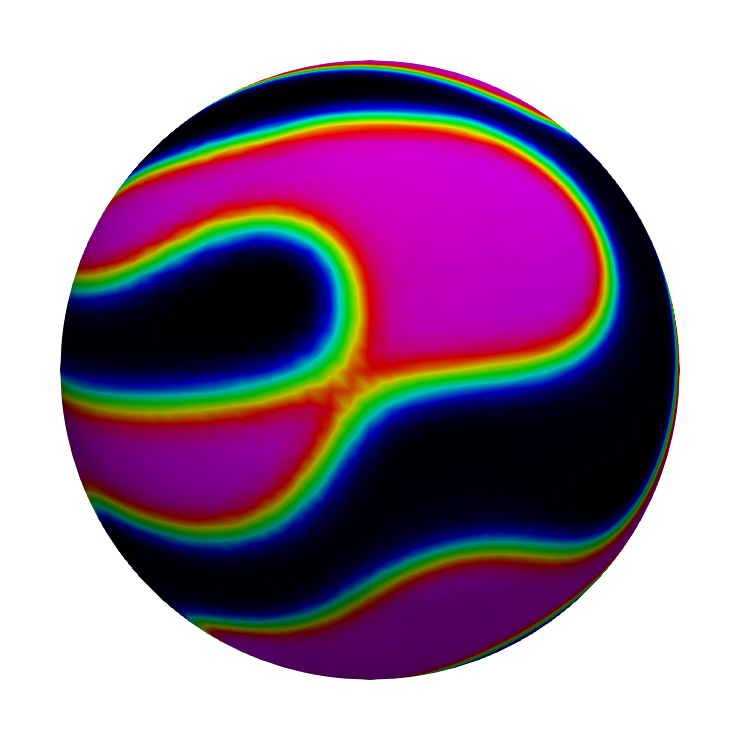}
\put(30,98){\small{$t = 30$}}
\end{overpic}
\begin{overpic}[width=.15\textwidth,grid=false]{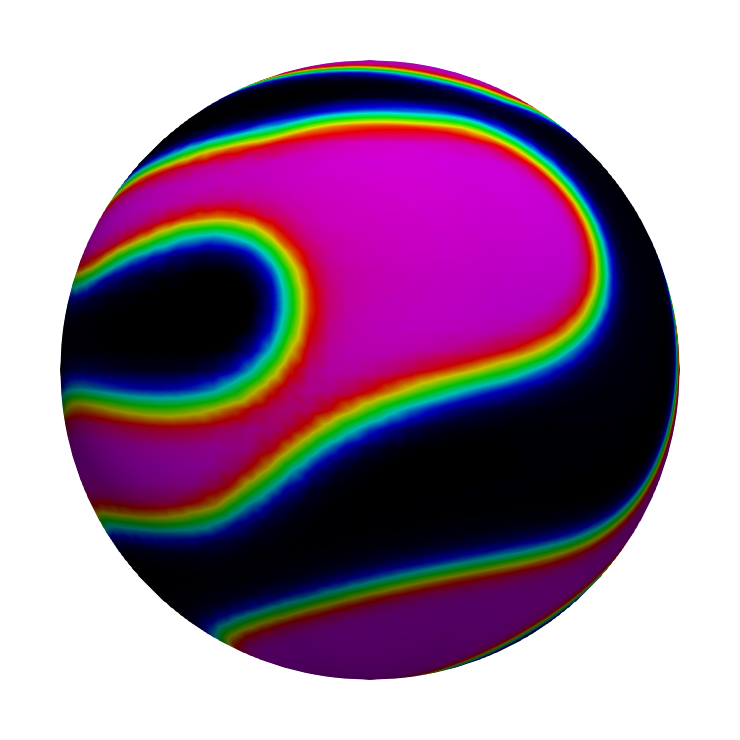}
\put(30,98){\small{$t = 50$}}
\end{overpic}
\begin{overpic}[width=.15\textwidth,grid=false]{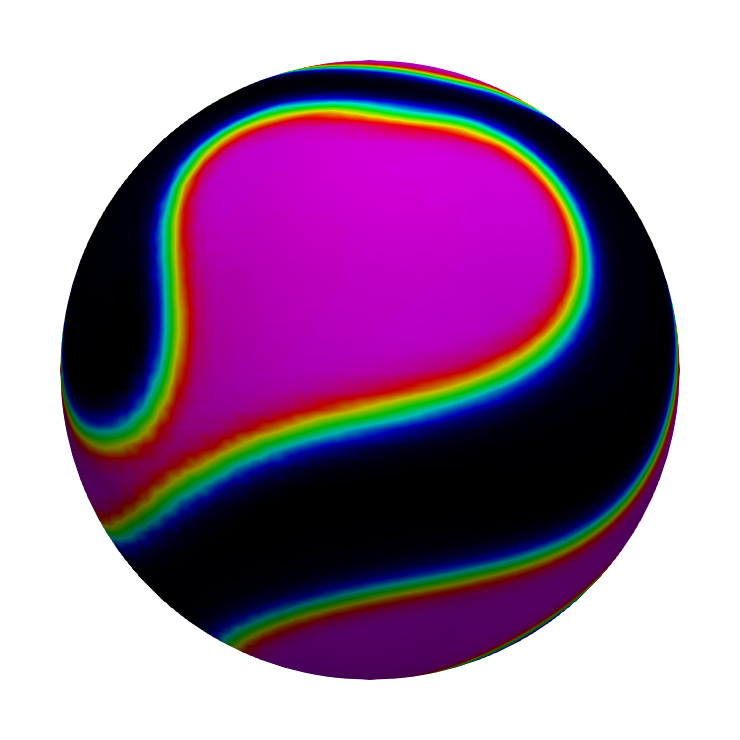}
\put(30,95){\small{$t = 100$}}
\end{overpic}
\\
\hskip .7cm
\begin{overpic}[width=.15\textwidth,grid=false]{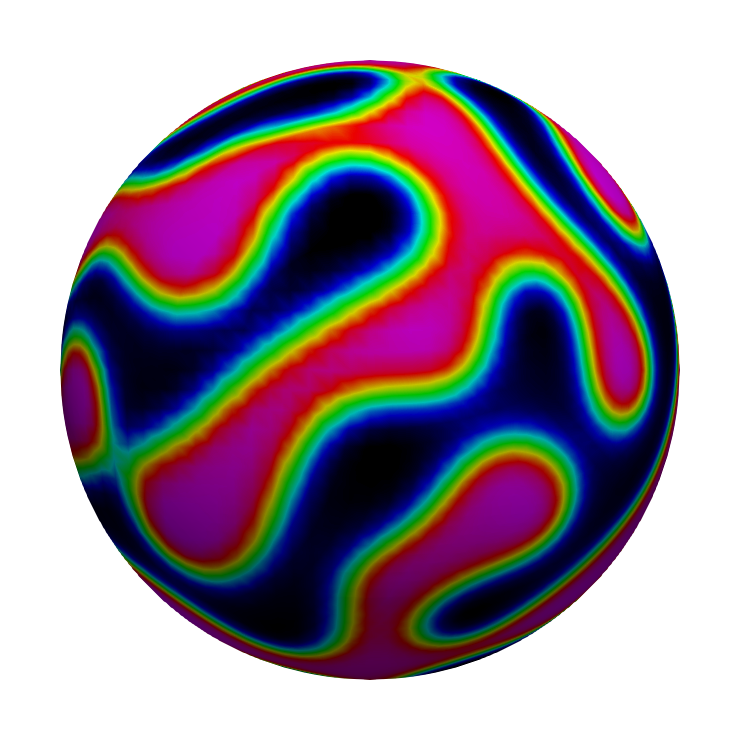}
\put(-50,50){\small{NSCH}}
\put(-50,33){\small{high $\eta$}}
\end{overpic}
\begin{overpic}[width=.15\textwidth,grid=false]{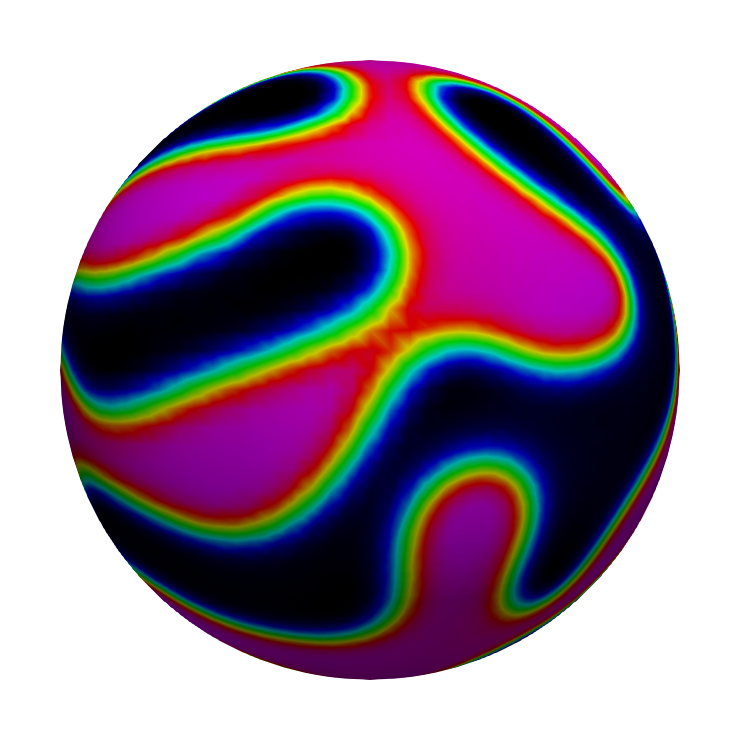}
\end{overpic}
\begin{overpic}[width=.15\textwidth,grid=false]{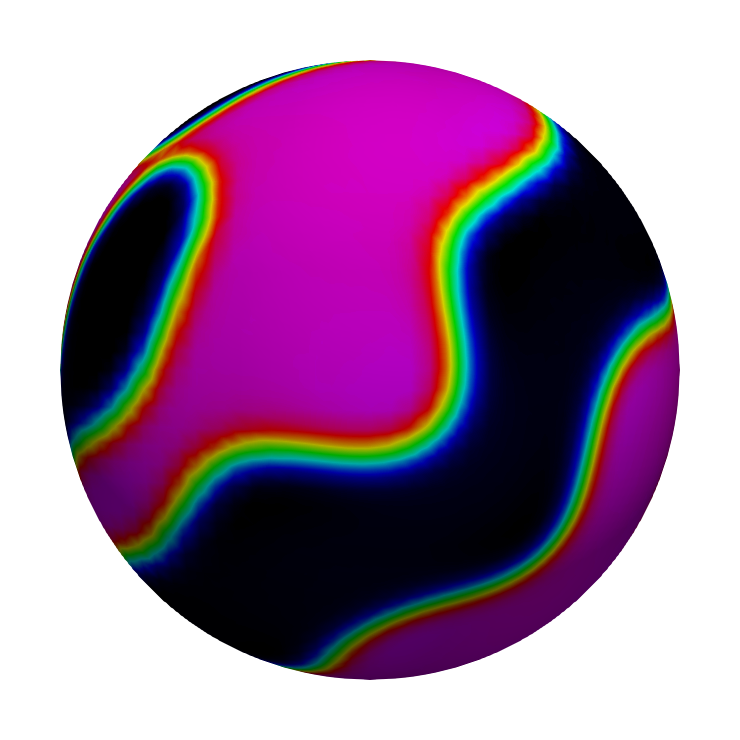}
\end{overpic}
\begin{overpic}[width=.15\textwidth,grid=false]{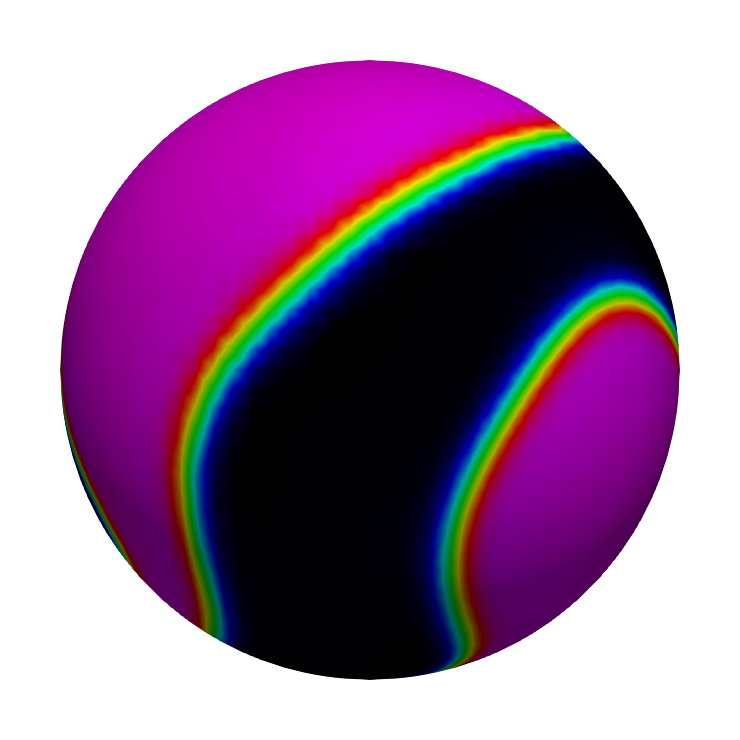}
\end{overpic}
\begin{overpic}[width=.15\textwidth,grid=false]{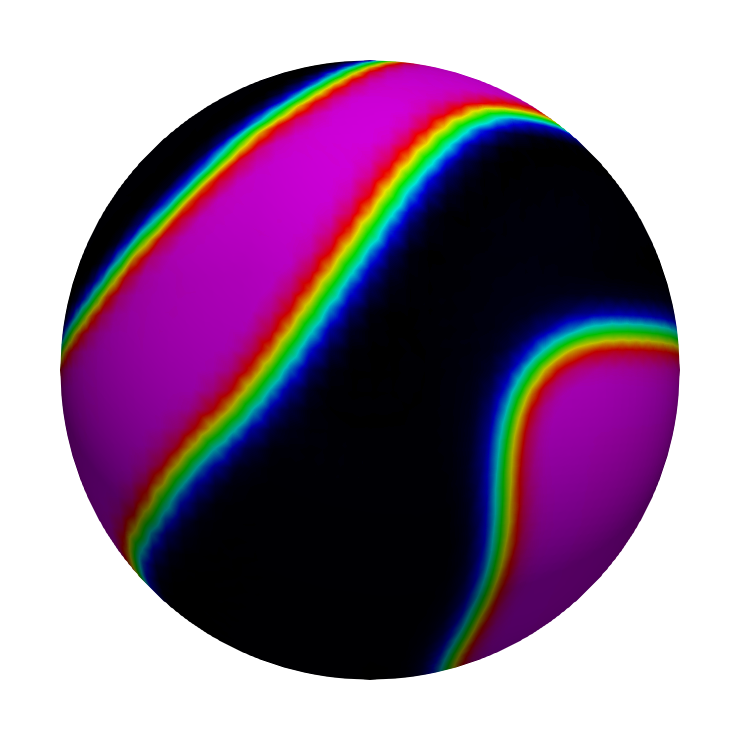}
\end{overpic}
\begin{overpic}[width=.15\textwidth,grid=false]{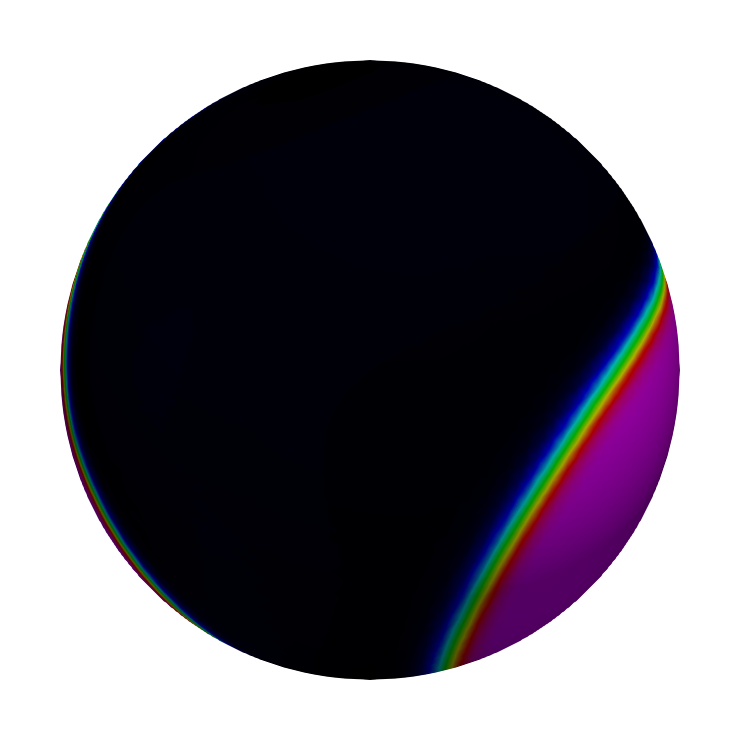}
\end{overpic}
\\
\hskip .7cm
\begin{overpic}[width=.15\textwidth,grid=false]{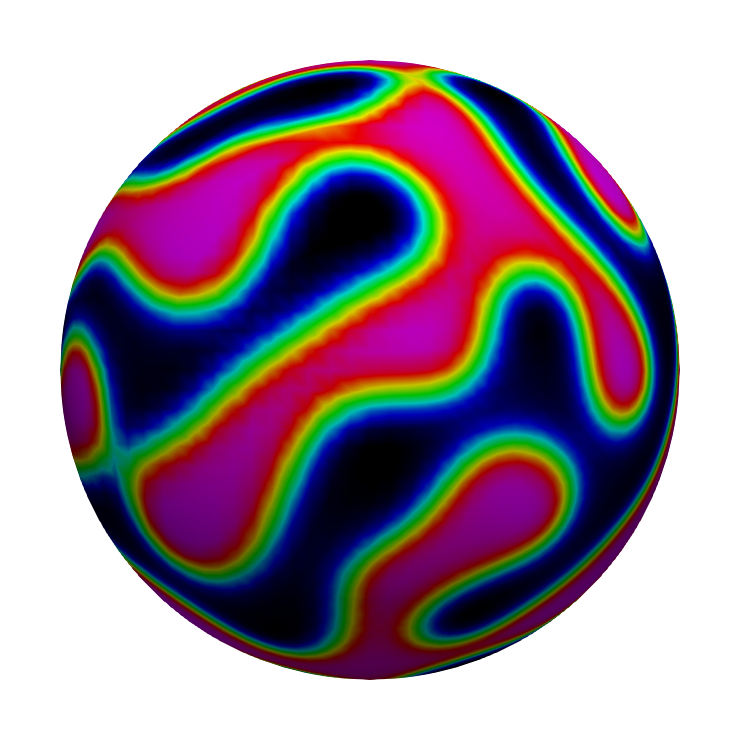}
\put(-50,50){\small{NSCH}}
\put(-50,33){\small{low $\eta$}}
\end{overpic}
\begin{overpic}[width=.15\textwidth,grid=false]{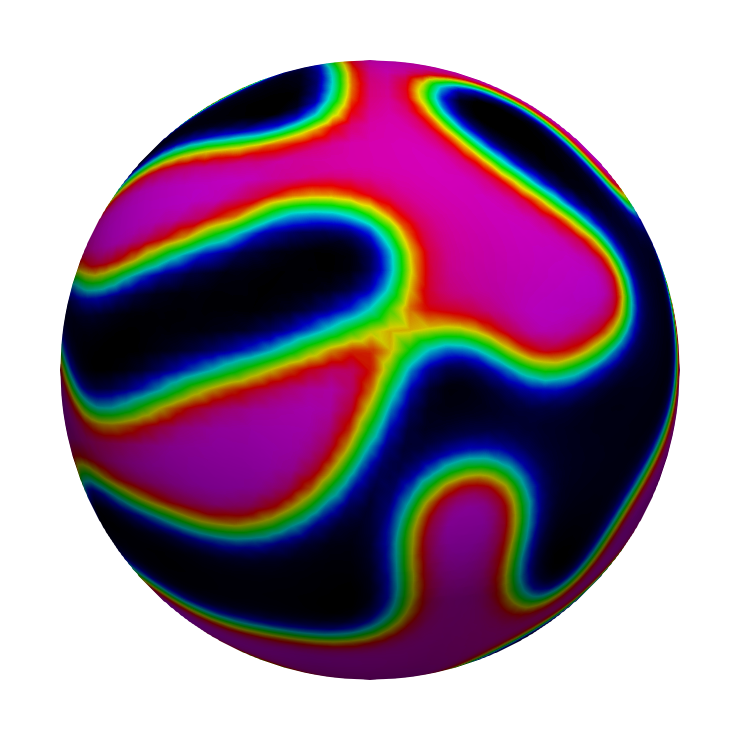}
\end{overpic}
\begin{overpic}[width=.15\textwidth,grid=false]{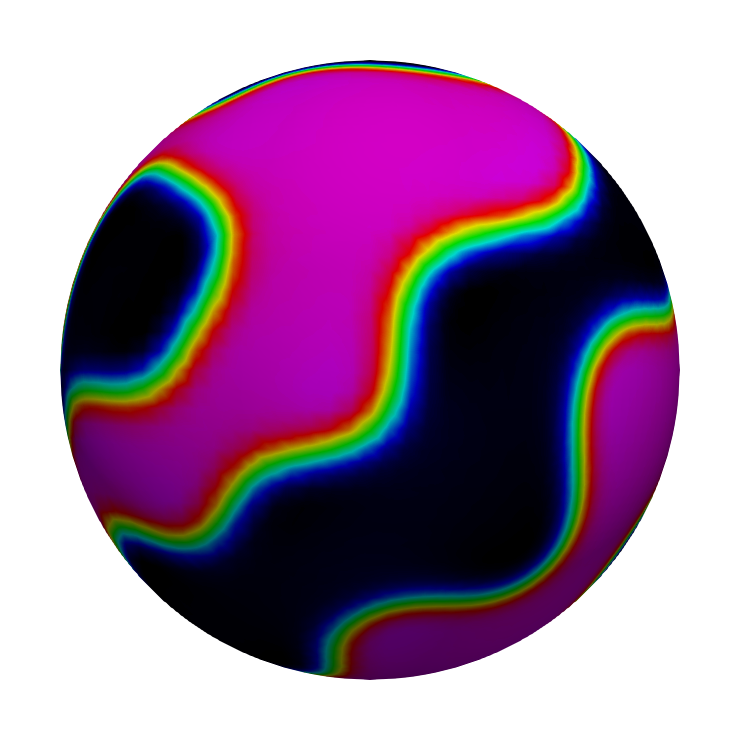}
\end{overpic}
\begin{overpic}[width=.15\textwidth,grid=false]{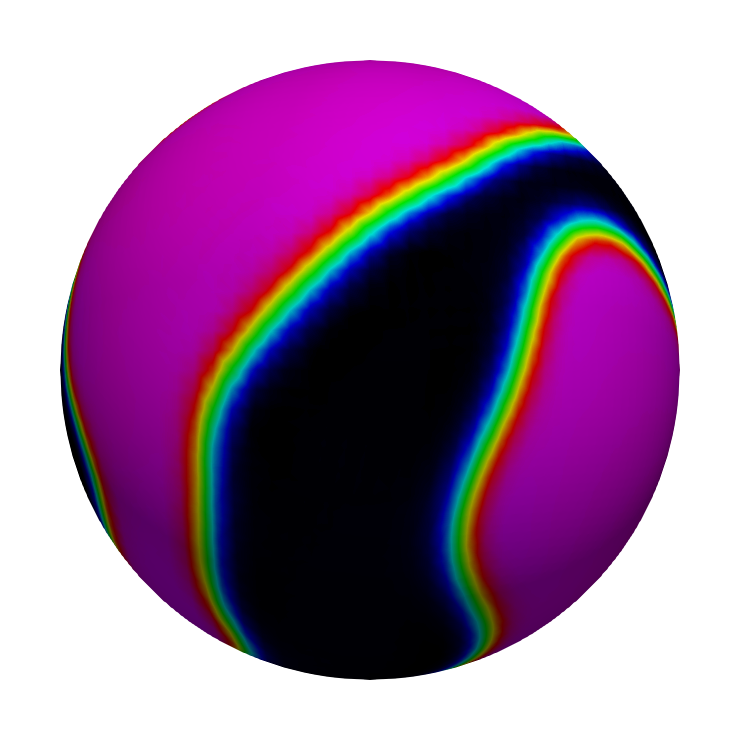}
\end{overpic}
\begin{overpic}[width=.15\textwidth,grid=false]{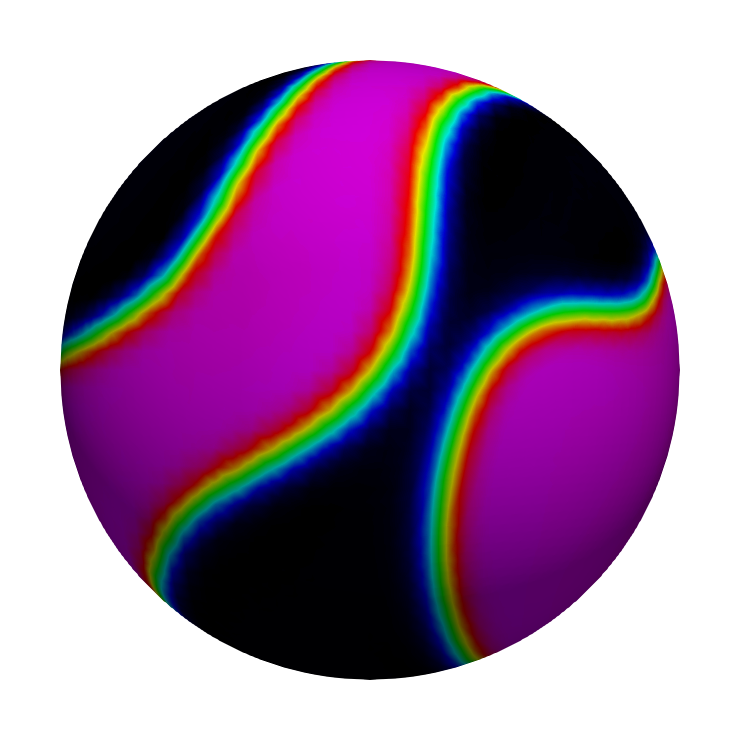}
\end{overpic}
\begin{overpic}[width=.15\textwidth,grid=false]{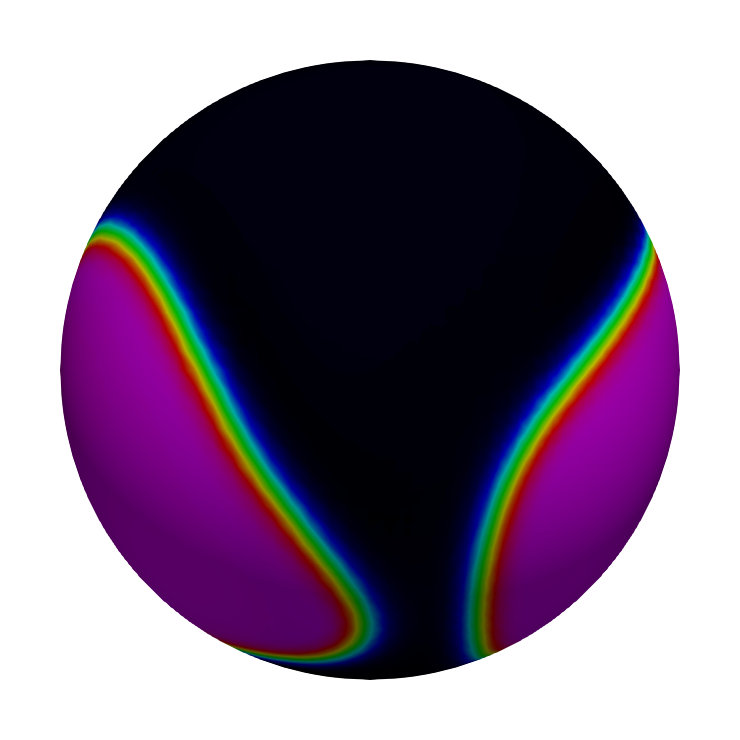}
\end{overpic}
\\
\vskip .2cm
\begin{overpic}[width=0.5\textwidth,grid=false,tics=10]{images/legend.png}
\end{overpic}
\end{center}
\caption{Phase separation given by the CH model (top), NSCH model with viscosities
$\eta_1 = 0.01$, $\eta_2 = 0.0008$ (center), and NSCH model with viscosities $\eta_1 = 0.0001$, $\eta_2 = 0.000008$ (bottom)
for composition 50\%-50\%.}\label{fig:sphere_Re_50}
\end{figure}

\begin{figure}
\begin{center}
\hskip .7cm
\begin{overpic}[width=.15\textwidth,grid=false]{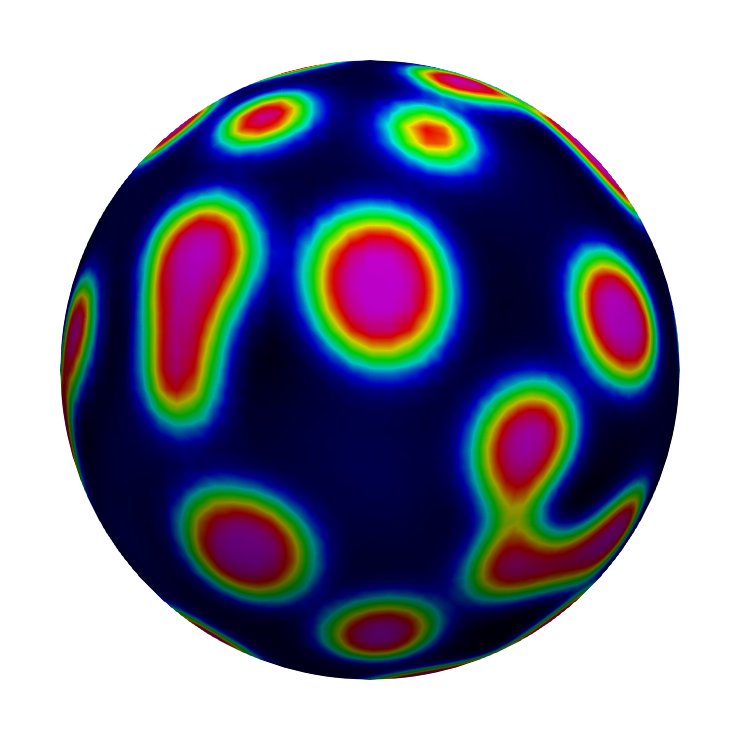}
\put(30,98){\small{$t = 2$}}
\put(-30,50){\small{CH}}
\end{overpic}
\begin{overpic}[width=.15\textwidth,grid=false]{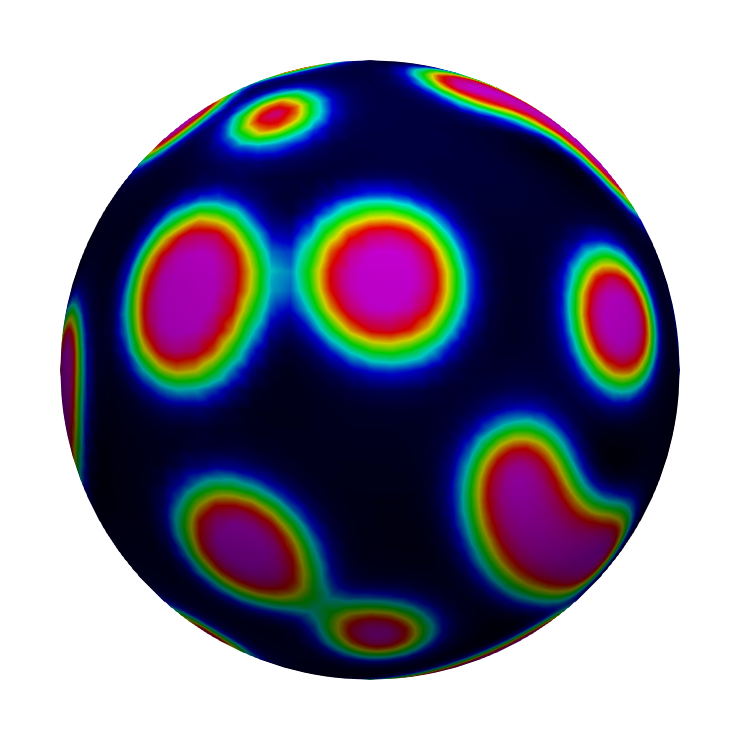}
\put(30,98){\small{$t = 5$}}
\end{overpic}
\begin{overpic}[width=.15\textwidth,grid=false]{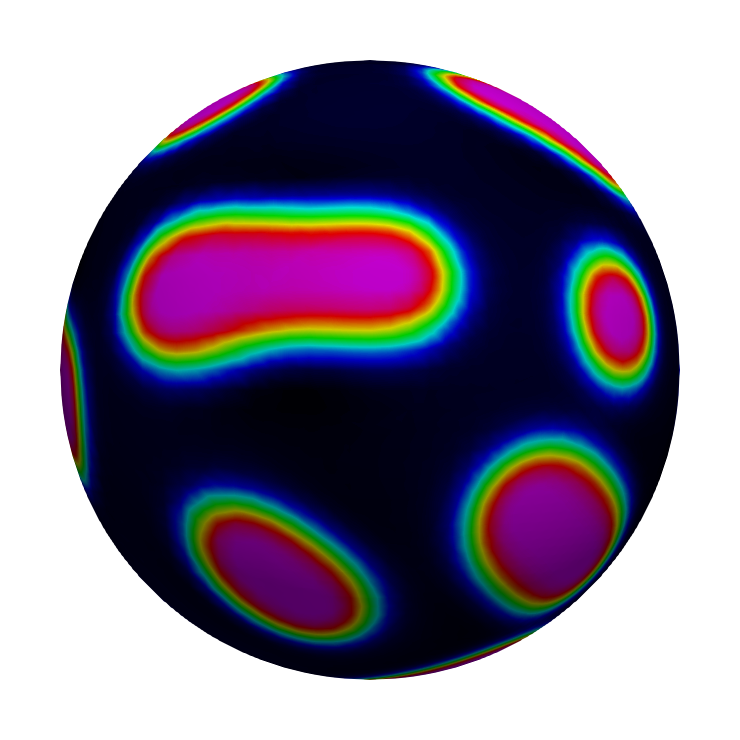}
\put(30,98){\small{$t = 10$}} 
\end{overpic}
\begin{overpic}[width=.15\textwidth,grid=false]{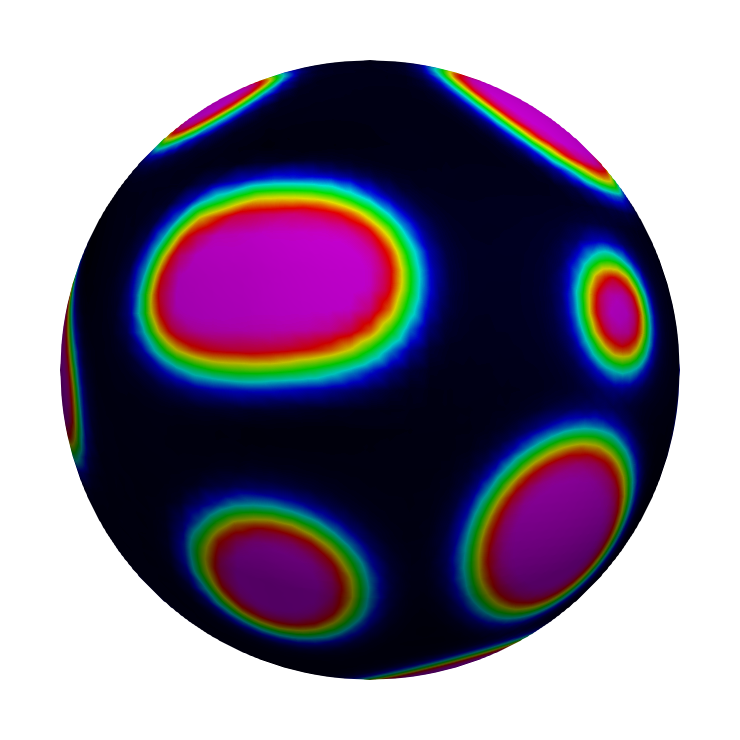}
\put(30,98){\small{$t = 20$}}
\end{overpic}
\begin{overpic}[width=.15\textwidth,grid=false]{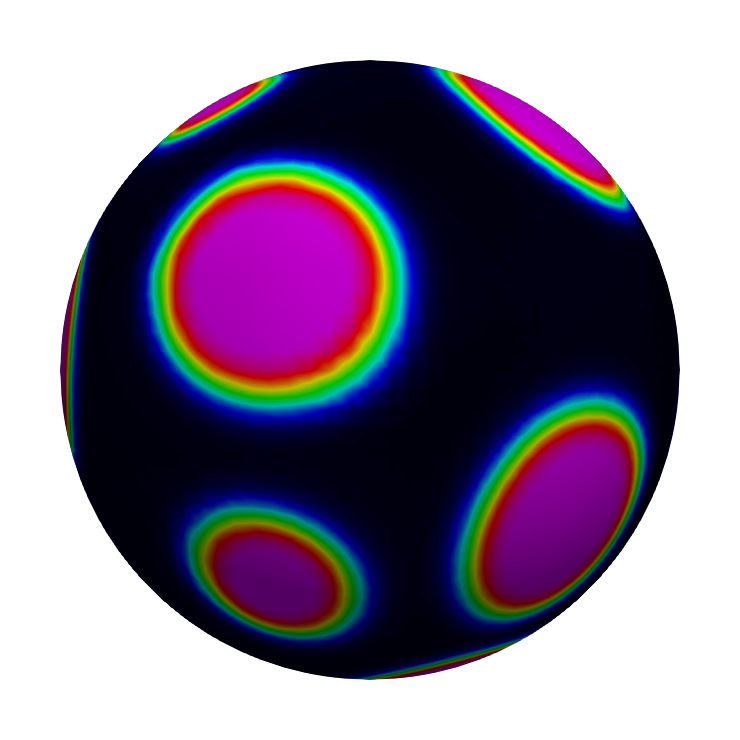}
\put(30,98){\small{$t = 40$}}
\end{overpic}
\begin{overpic}[width=.15\textwidth,grid=false]{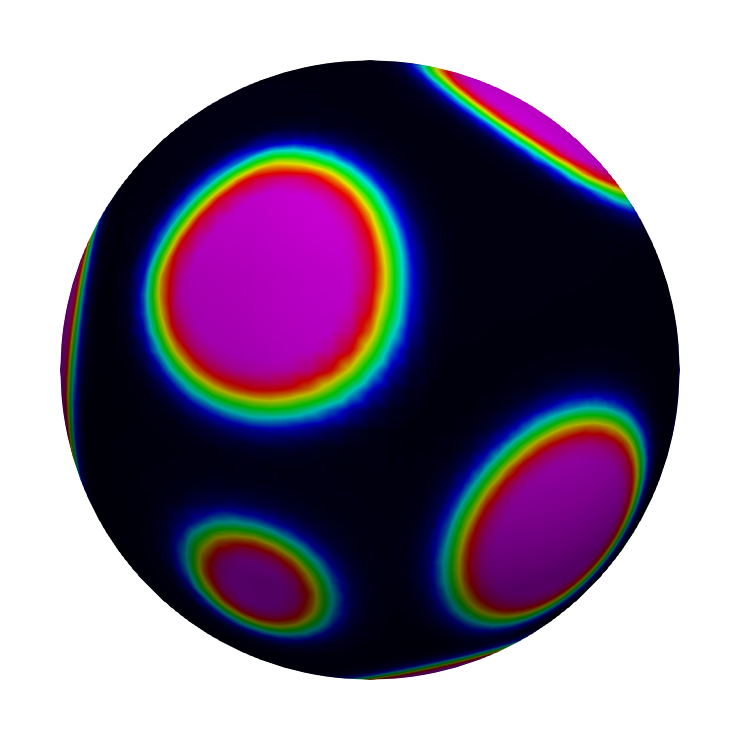}
\put(30,95){\small{$t = 100$}}
\end{overpic}
\\
\hskip .7cm
\begin{overpic}[width=.15\textwidth,grid=false]{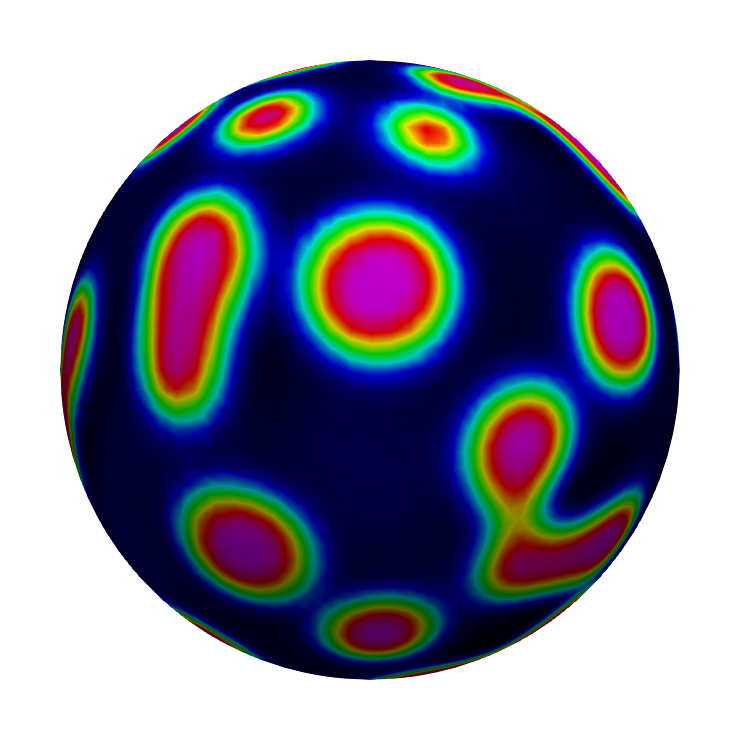}
\put(-50,50){\small{NSCH}}
\put(-50,33){\small{high $\eta$}}
\end{overpic}
\begin{overpic}[width=.15\textwidth,grid=false]{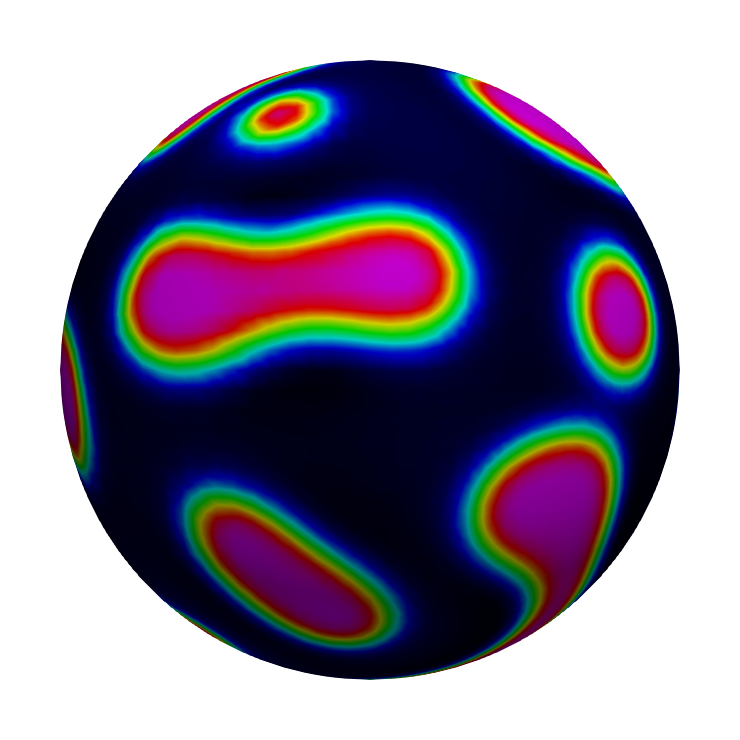}
\end{overpic}
\begin{overpic}[width=.15\textwidth,grid=false]{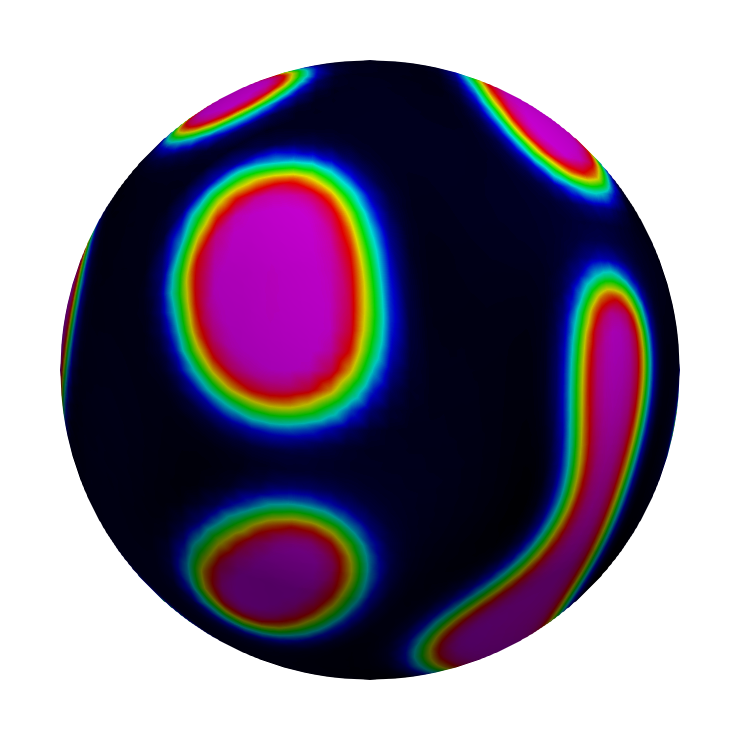}
\end{overpic}
\begin{overpic}[width=.15\textwidth,grid=false]{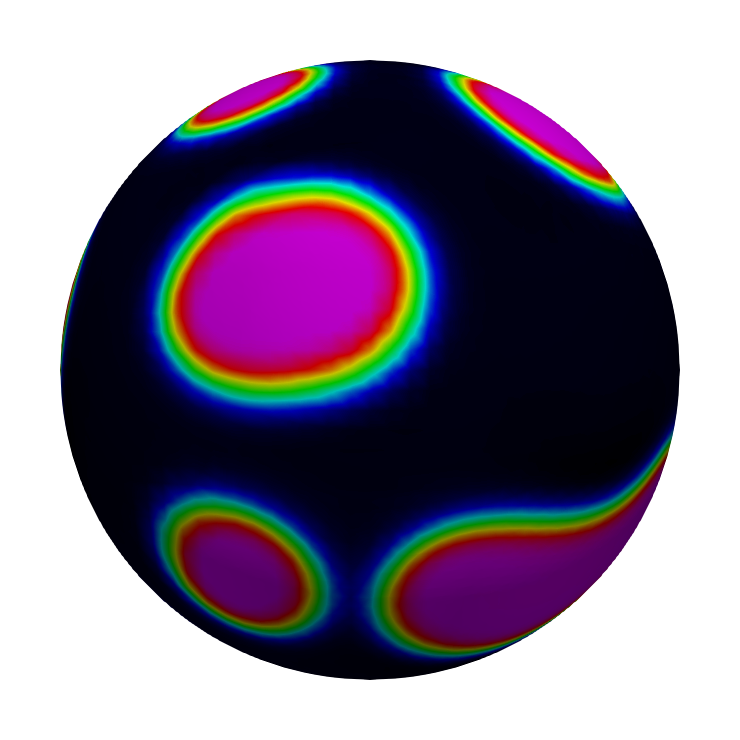}
\end{overpic}
\begin{overpic}[width=.15\textwidth,grid=false]{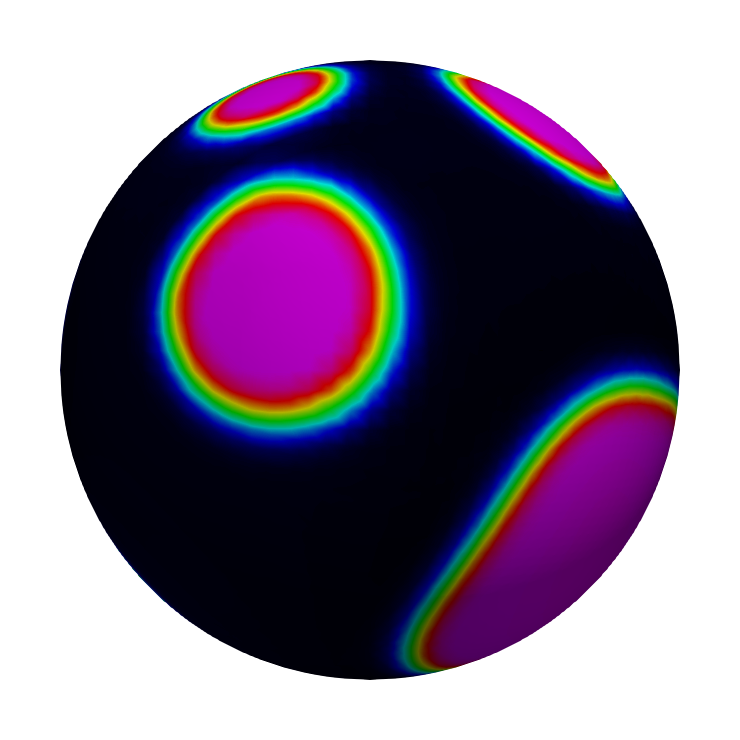}
\end{overpic}
\begin{overpic}[width=.15\textwidth,grid=false]{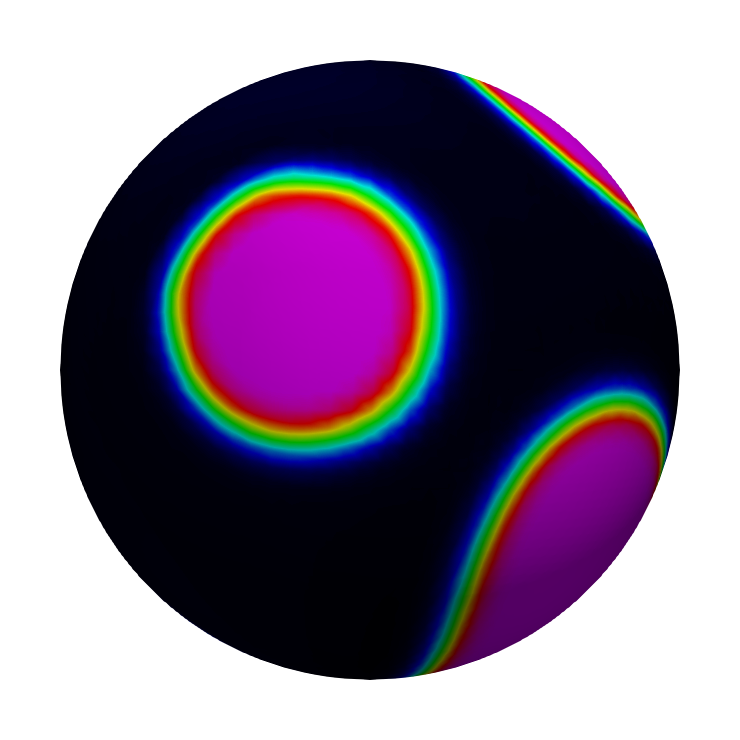}
\end{overpic}
\\
\hskip .7cm
\begin{overpic}[width=.15\textwidth,grid=false]{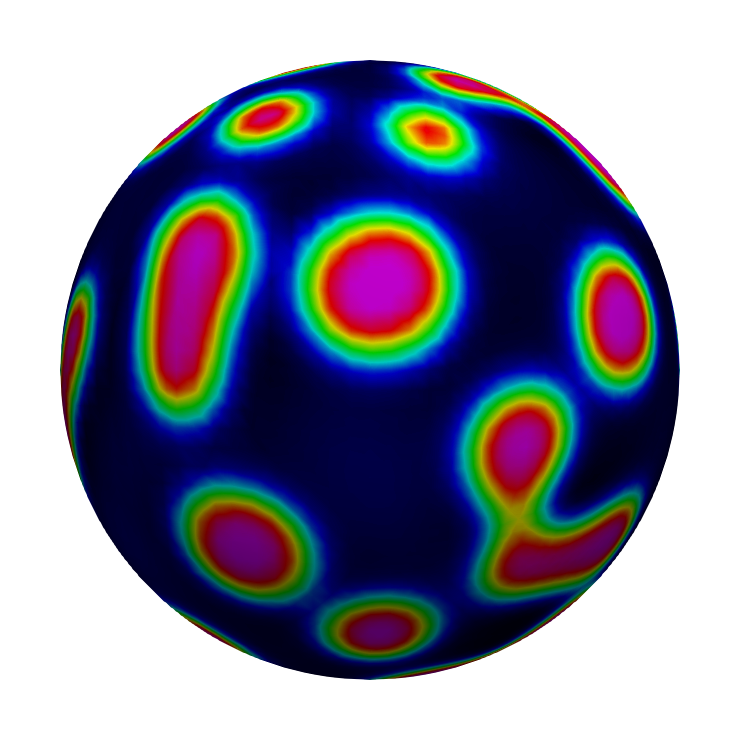}
\put(-50,50){\small{NSCH}}
\put(-50,33){\small{low $\eta$}}
\end{overpic}
\begin{overpic}[width=.15\textwidth,grid=false]{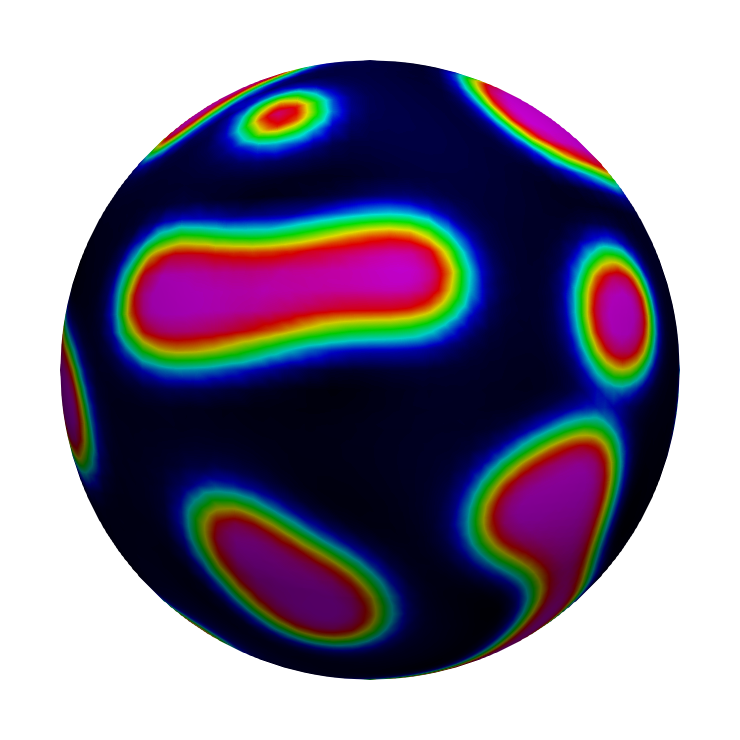}
\end{overpic}
\begin{overpic}[width=.15\textwidth,grid=false]{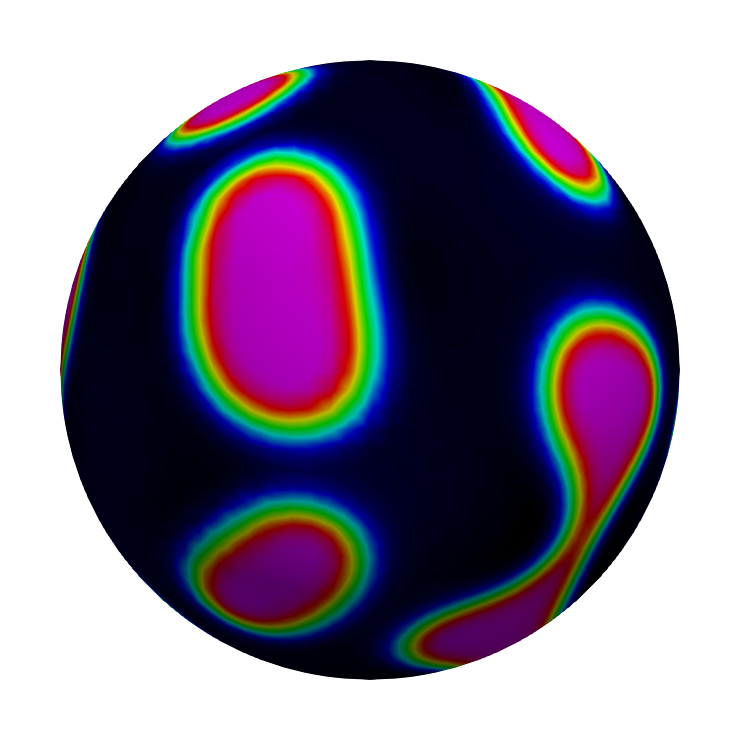}
\end{overpic}
\begin{overpic}[width=.15\textwidth,grid=false]{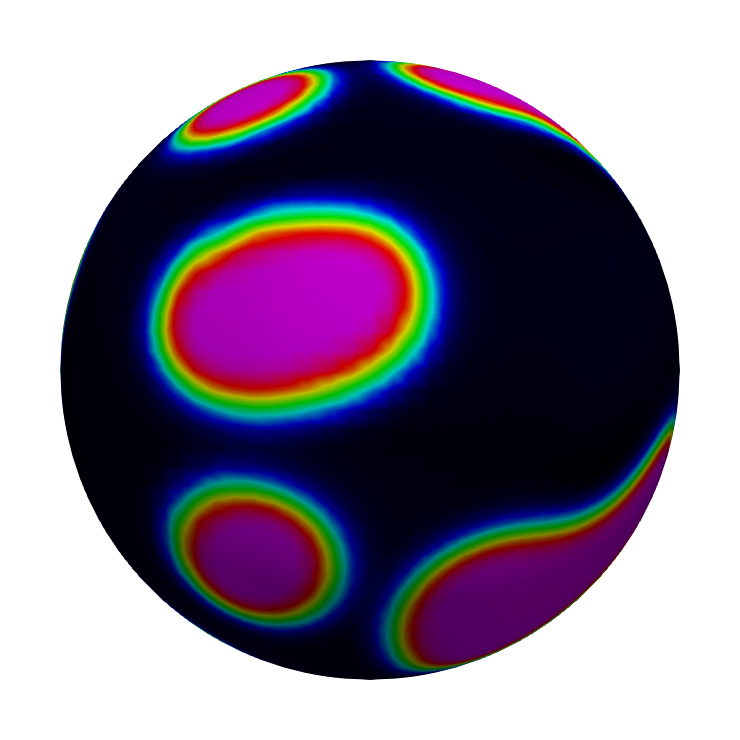}
\end{overpic}
\begin{overpic}[width=.15\textwidth,grid=false]{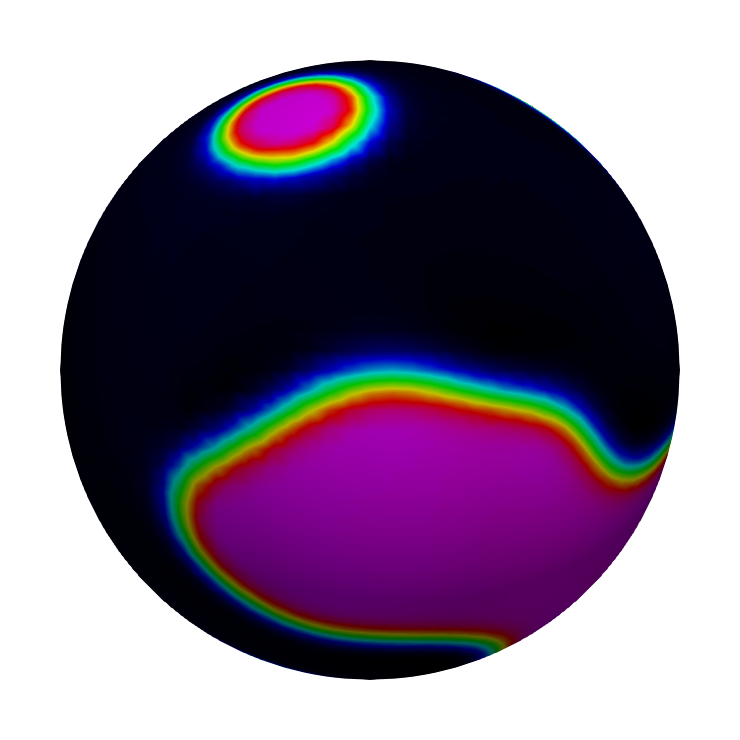}
\end{overpic}
\begin{overpic}[width=.15\textwidth,grid=false]{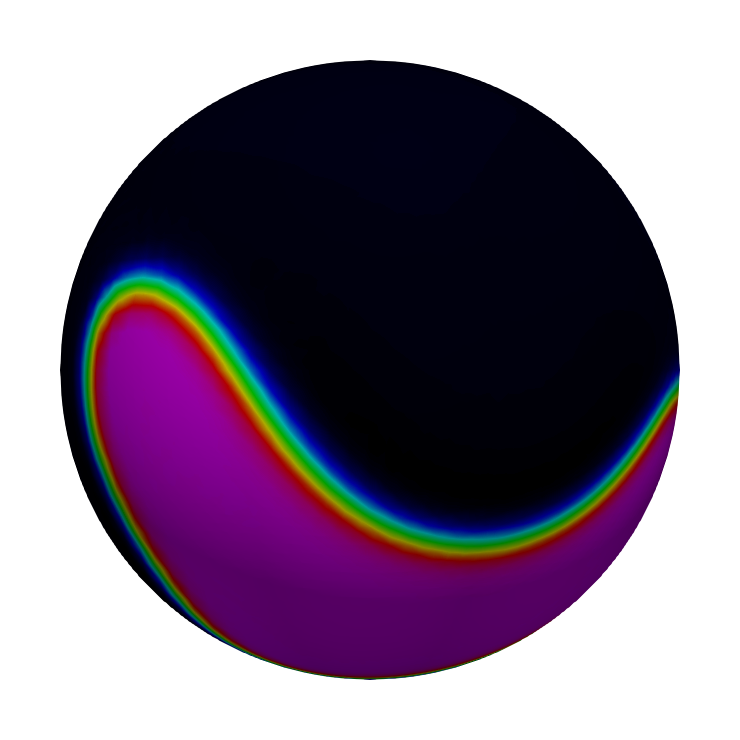}
\end{overpic}\\
\vskip .2cm
\begin{overpic}[width=0.5\textwidth,grid=false,tics=10]{images/legend.png}
\end{overpic}
\end{center}
\caption{Phase separation given by the CH model (top), NSCH model with viscosities
$\eta_1 = 0.01$, $\eta_2 = 0.0008$ (center), and NSCH model with viscosities $\eta_1 = 0.0001$, $\eta_2 = 0.000008$ (bottom)
for composition 30\%-70\%.}\label{fig:sphere_Re_30}
\end{figure}

Fig.~\ref{fig:sphere_Re_50} and \ref{fig:sphere_Re_30} show the evolution of phase
separation by the CH model and NSCH model for the high viscosity and low viscosity cases
for composition 50\%-50\% and 30\%-70\%, respectively. The patterns are very different
for the two compositions: composition 50\%-50\% gives rise to pink macrodomains with a tortuous
interface, while from composition 30\%-70\% one gets many small domains with a more or less elongated shape. 
Switching from high viscosities to low viscosities in the NSCH model does not produce significant differences
in the appearance of the domains until $t = 100$ for composition 50\%-50\%. Compare center and bottom row in 
Fig.~\ref{fig:sphere_Re_50}. This is somewhat expected from the fact that the lines corresponding to the two cases in 
Fig.~\ref{fig:free_energy_eta} (left) are almost superimposed until about $t = 80$. The change in domain appearance 
happens faster when going from high viscosities to low viscosities for composition 30\%-70\%. Indeed, by comparing
center and bottom row in Fig.~\ref{fig:sphere_Re_30} we observe remarkable differences already at $t = 40$. Again,
this can have been expected from looking at the two lines for the NSCH model in  Fig.~\ref{fig:free_energy_eta} (right).

\begin{figure}
\begin{center}
\hskip .7cm
\begin{overpic}[width=.15\textwidth,grid=false]{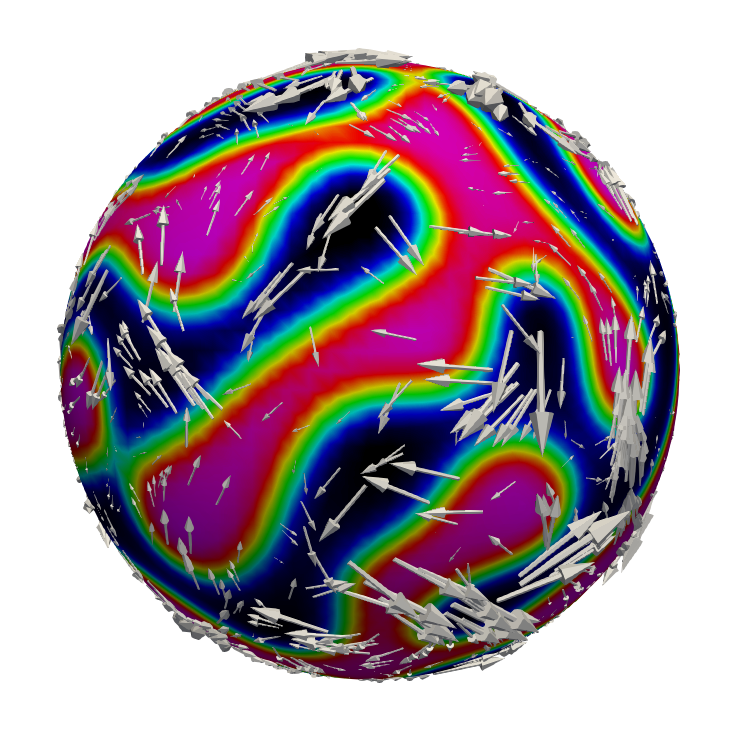}
\put(-50,45){\small{high $\eta$}}
\put(30,98){\small{$t = 2$}}
\end{overpic}
\begin{overpic}[width=.15\textwidth,grid=false]{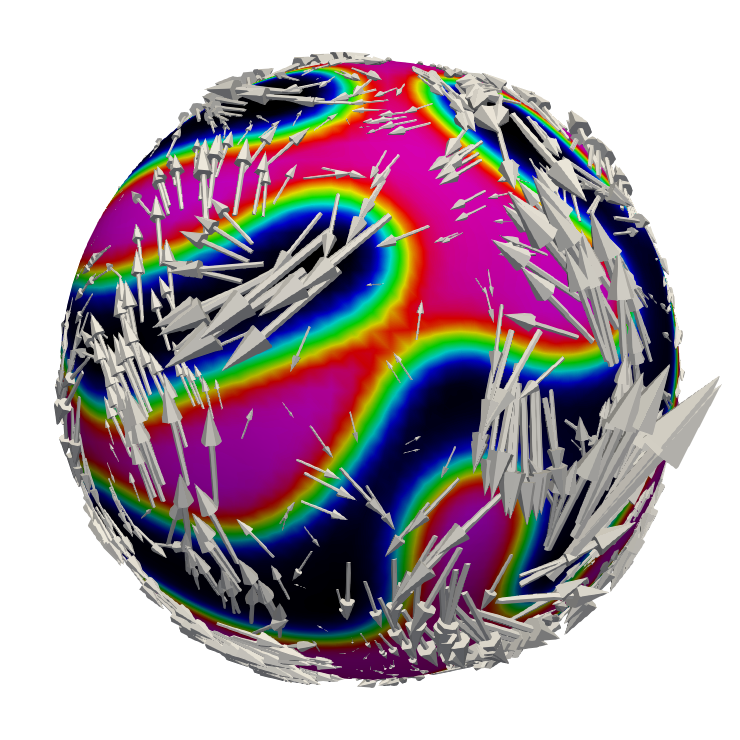}
\put(30,98){\small{$t = 5$}}
\end{overpic}
\begin{overpic}[width=.15\textwidth,grid=false]{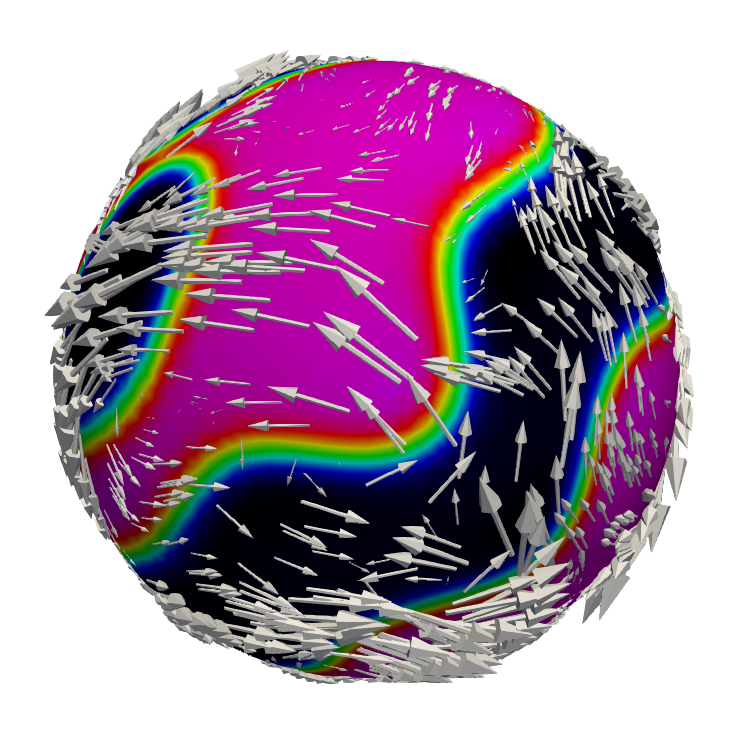}
\put(28,98){\small{$t = 15$}}
\end{overpic}
\begin{overpic}[width=.15\textwidth,grid=false]{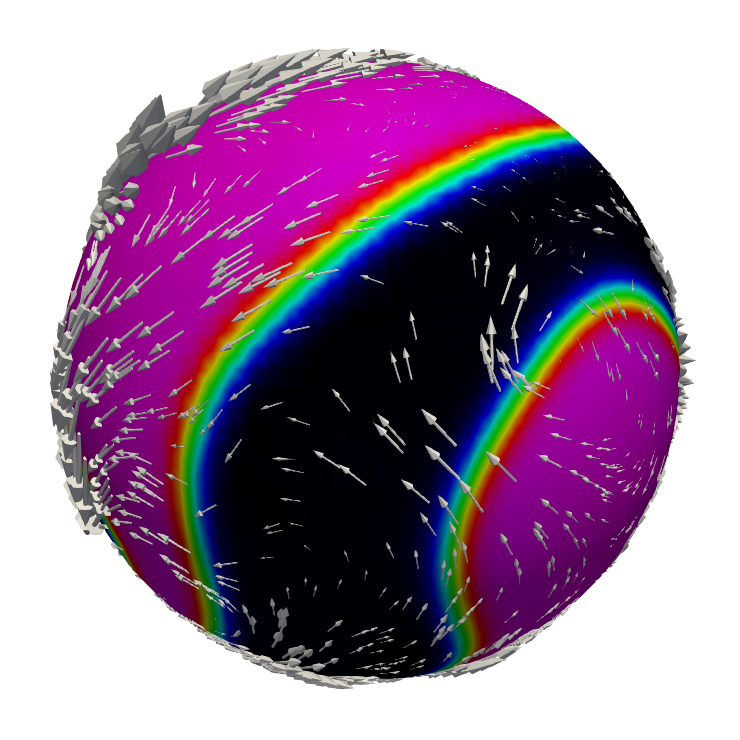}
\put(28,98){\small{$t = 30$}}
\end{overpic}
\begin{overpic}[width=.15\textwidth,grid=false]{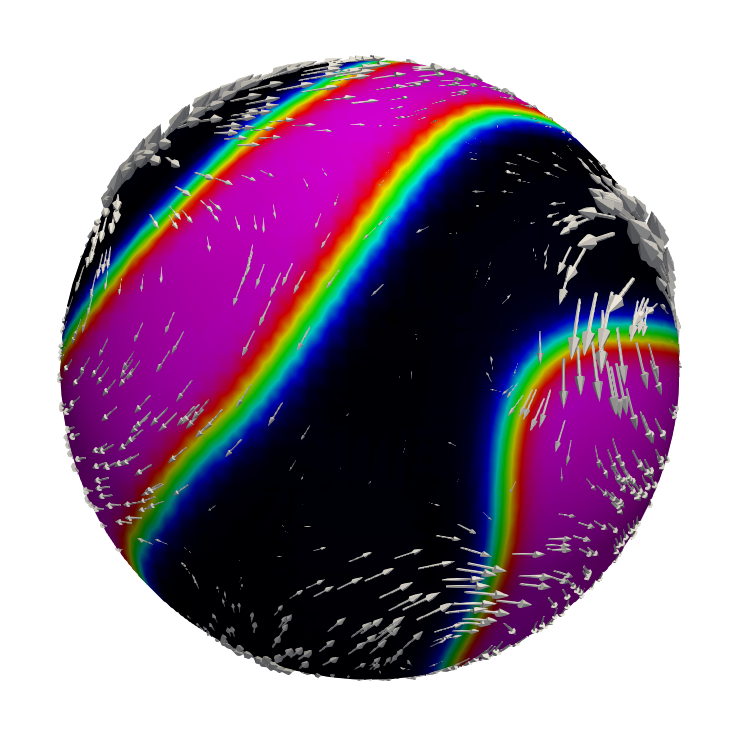}
\put(28,98){\small{$t = 50$}}
\end{overpic}
\begin{overpic}[width=.15\textwidth,grid=false]{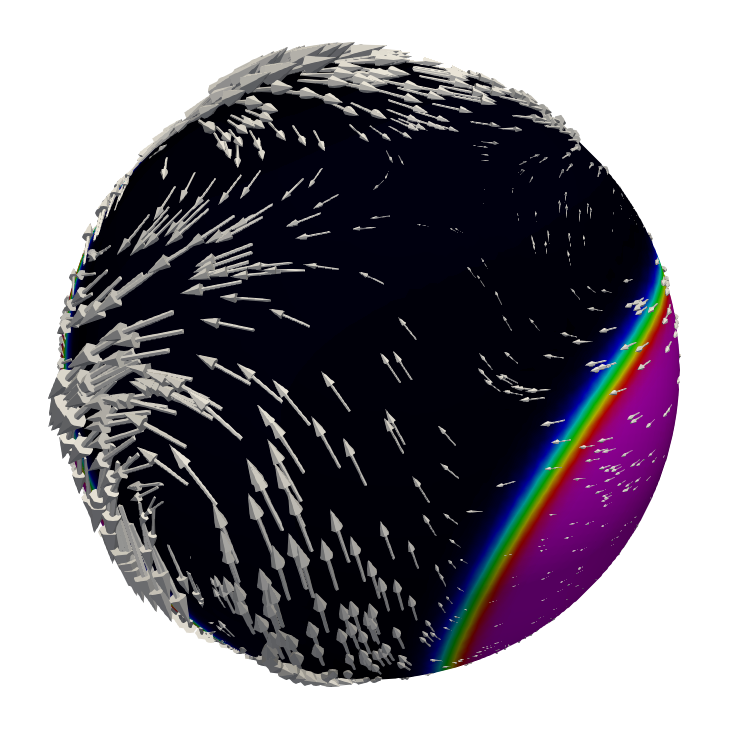}
\put(26,98){\small{$t = 100$}}
\end{overpic}
\\
\hskip .7cm
\begin{overpic}[width=.15\textwidth,grid=false]{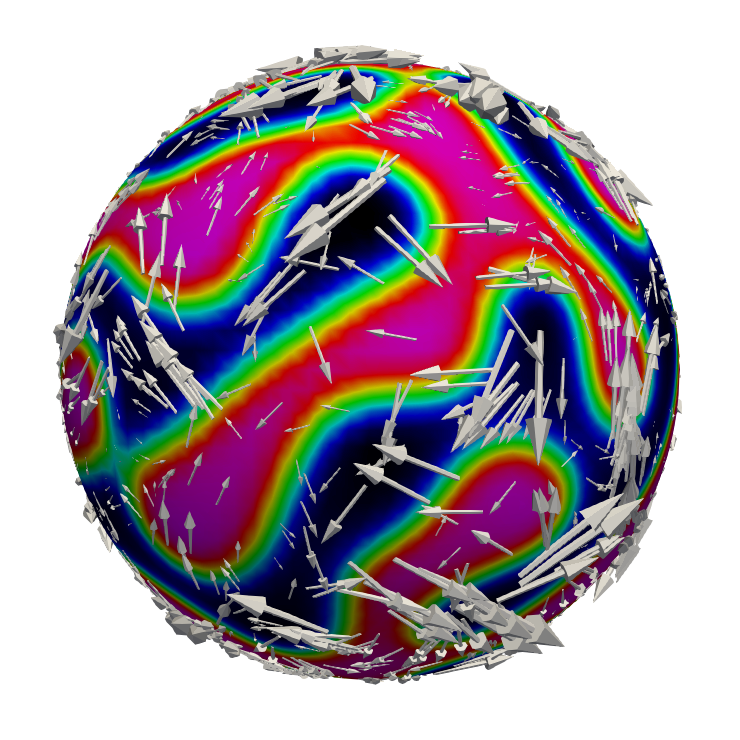}
\put(-50,45){\small{low $\eta$}}
\end{overpic}
\begin{overpic}[width=.15\textwidth,grid=false]{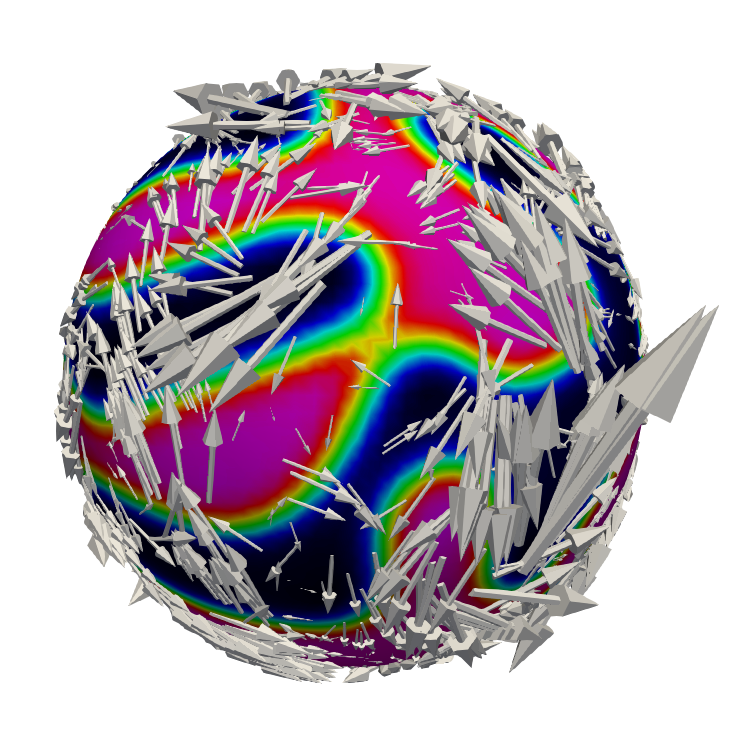}
\end{overpic}
\begin{overpic}[width=.15\textwidth,grid=false]{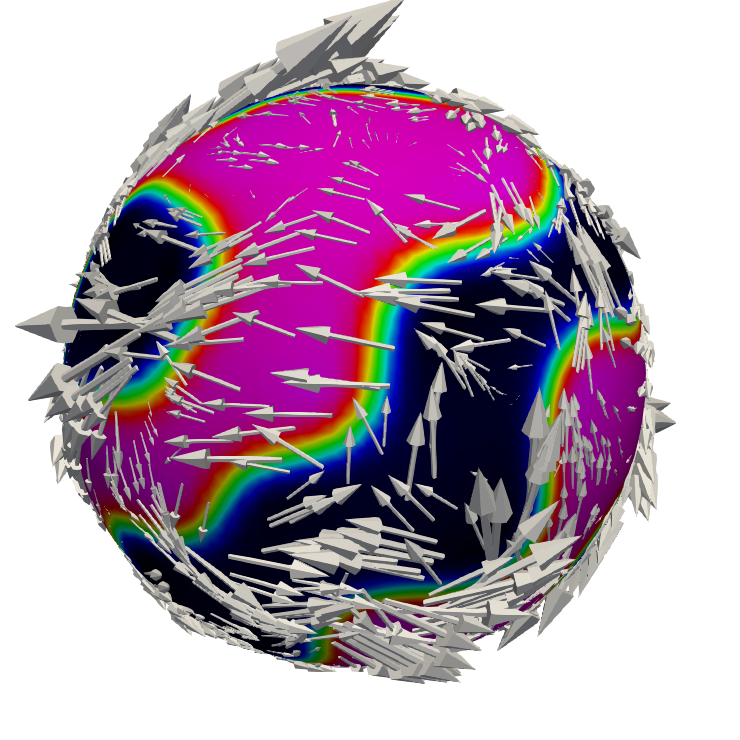}
\end{overpic}
\begin{overpic}[width=.15\textwidth,grid=false]{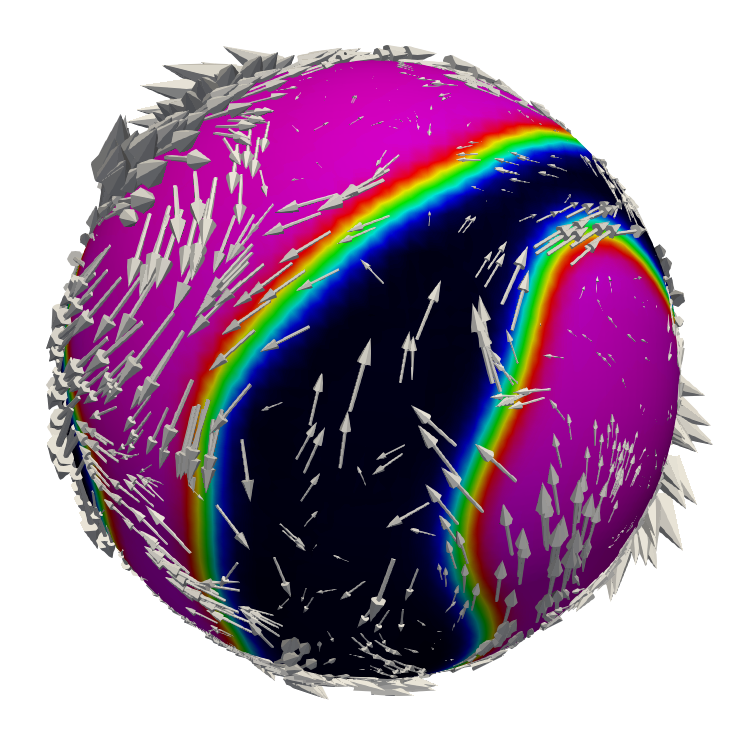}
\end{overpic}
\begin{overpic}[width=.15\textwidth,grid=false]{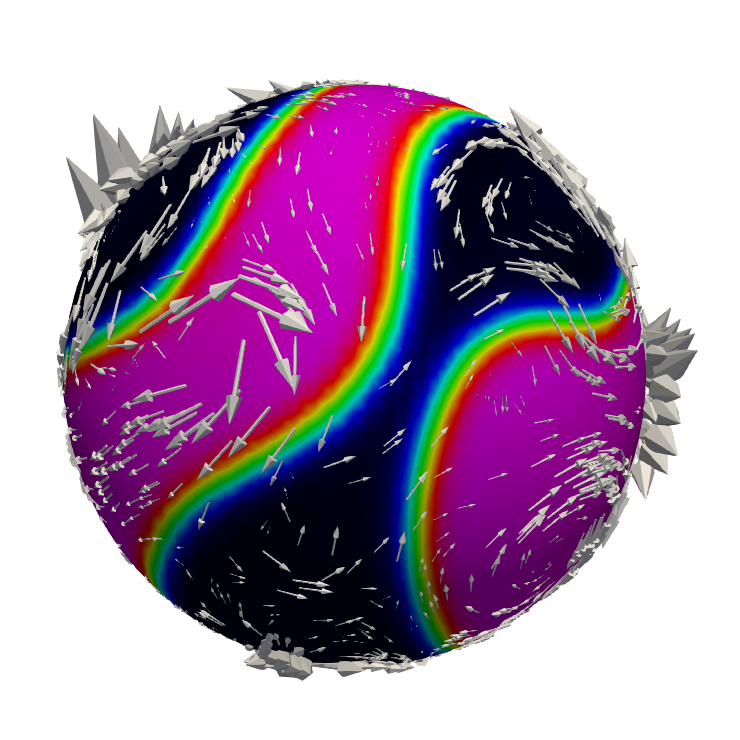}
\end{overpic}
\begin{overpic}[width=.15\textwidth,grid=false]{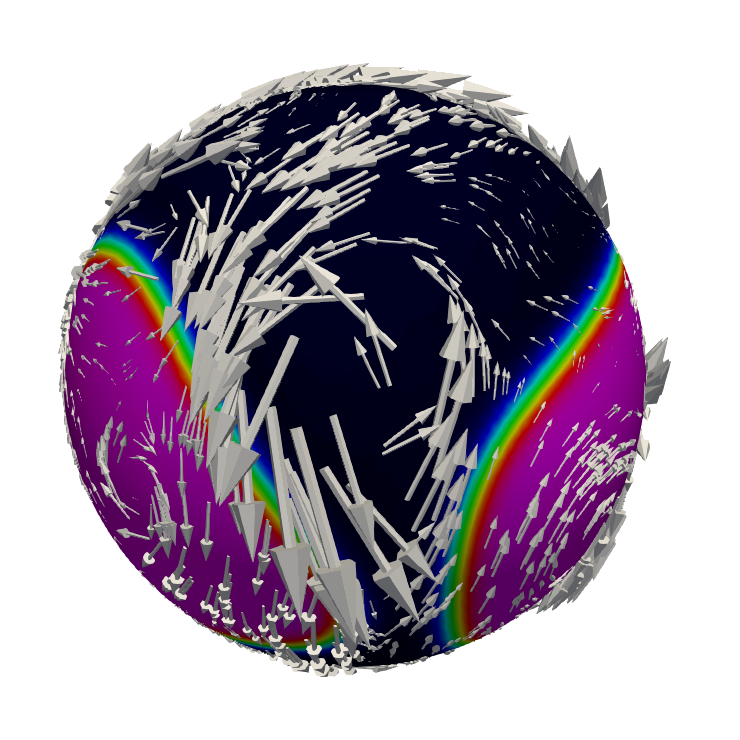}
\end{overpic}
\end{center}
\caption{Velocity vectors superimposed to the surface fraction for $\eta_1 = 10^{-2}, \eta_2 = 8 \cdot 10^{-4}$ 
(top row) and $\eta_1 = 10^{-4}, \eta_2 = 8 \cdot 10^{-6}$ (bottom row) for 
composition 50\%-50\%. For visualization purposes, the  
velocity vectors are magnified by a factor 5 in both rows. 
}\label{fig:flow_Re_50}
\end{figure}

\begin{figure}
\begin{center}
\hskip .7cm
\begin{overpic}[width=.15\textwidth,grid=false]{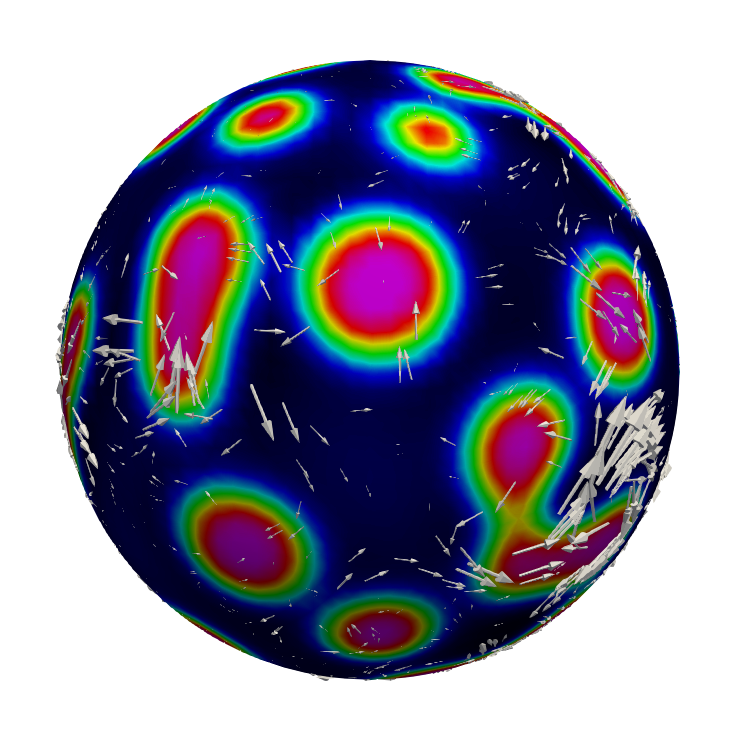}
\put(-50,45){\small{high $\eta$}}
\put(30,98){\small{$t = 2$}}
\end{overpic}
\begin{overpic}[width=.15\textwidth,grid=false]{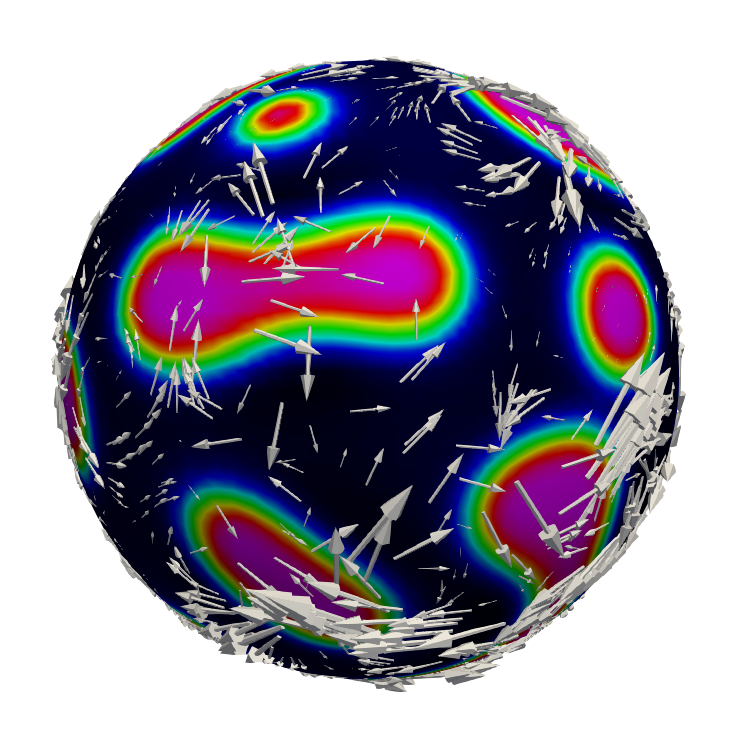}
\put(30,98){\small{$t = 5$}}
\end{overpic}
\begin{overpic}[width=.15\textwidth,grid=false]{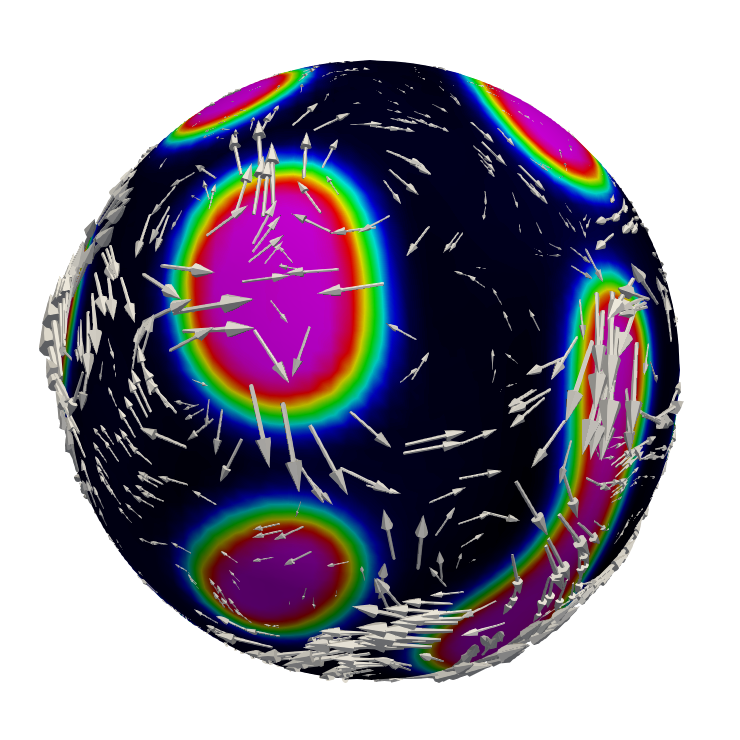}
\put(28,98){\small{$t = 10$}}
\end{overpic}
\begin{overpic}[width=.15\textwidth,grid=false]{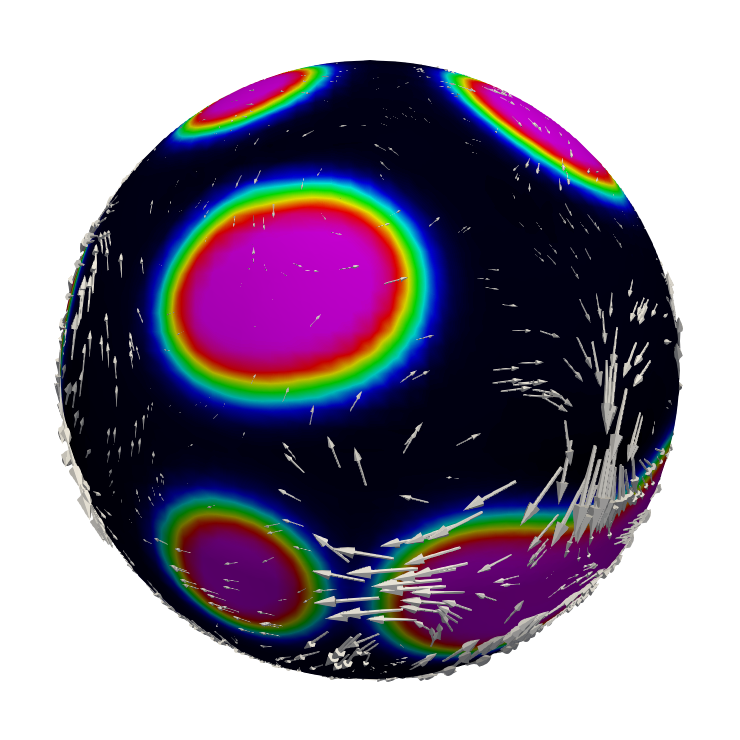}
\put(28,98){\small{$t = 20$}}
\end{overpic}
\begin{overpic}[width=.15\textwidth,grid=false]{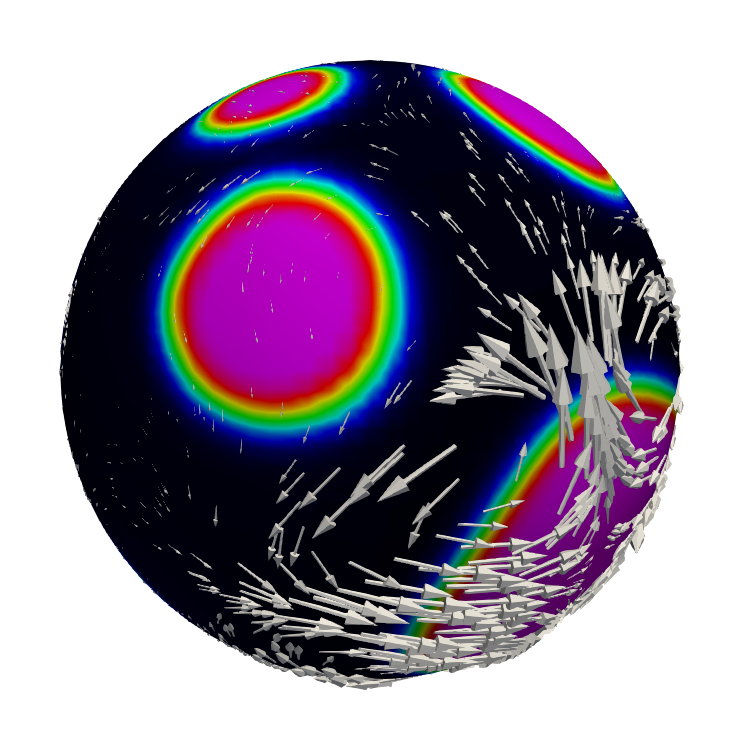}
\put(28,98){\small{$t = 40$}}
\end{overpic}
\begin{overpic}[width=.15\textwidth,grid=false]{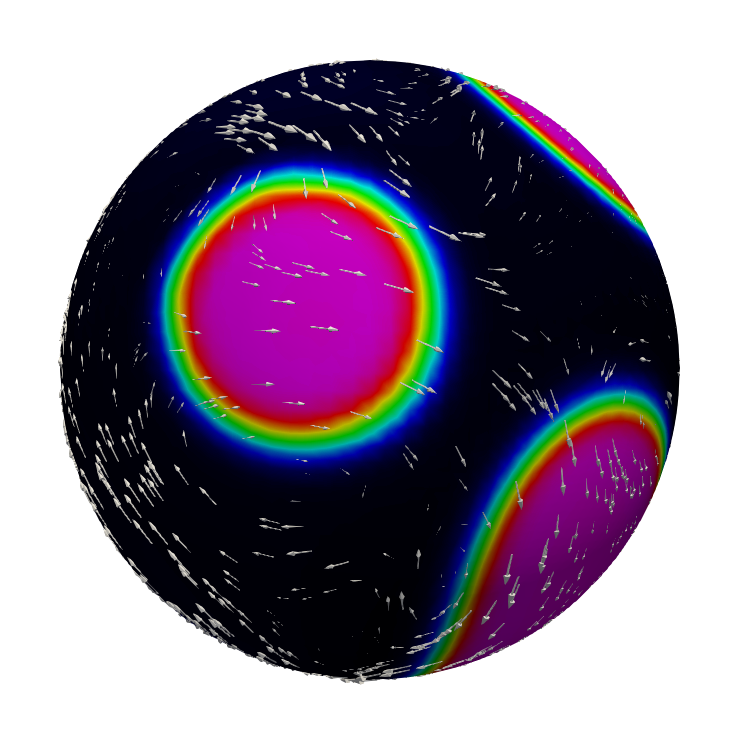}
\put(26,98){\small{$t = 100$}}
\end{overpic}
\\
\hskip .7cm
\begin{overpic}[width=.15\textwidth,grid=false]{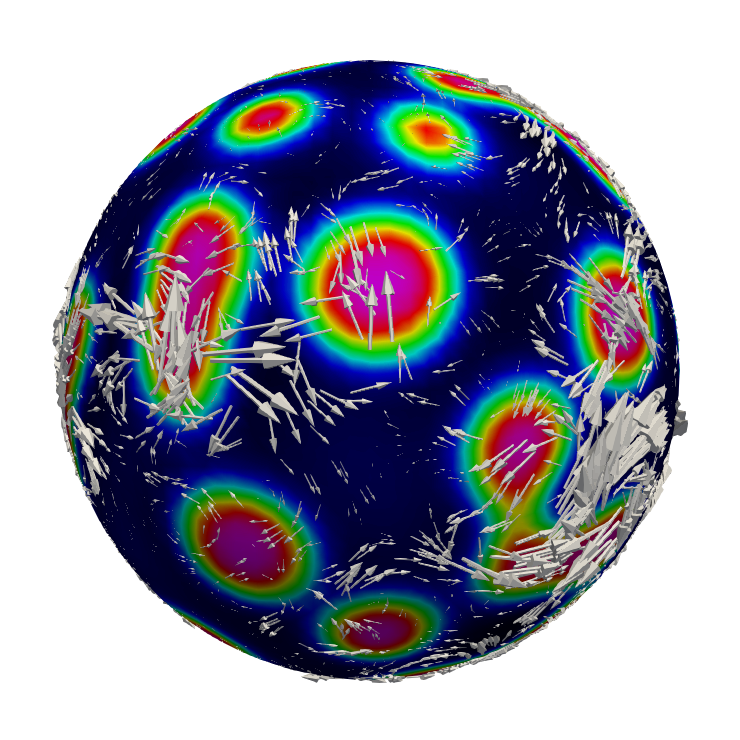}
\put(-50,45){\small{low $\eta$}}
\end{overpic}
\begin{overpic}[width=.15\textwidth,grid=false]{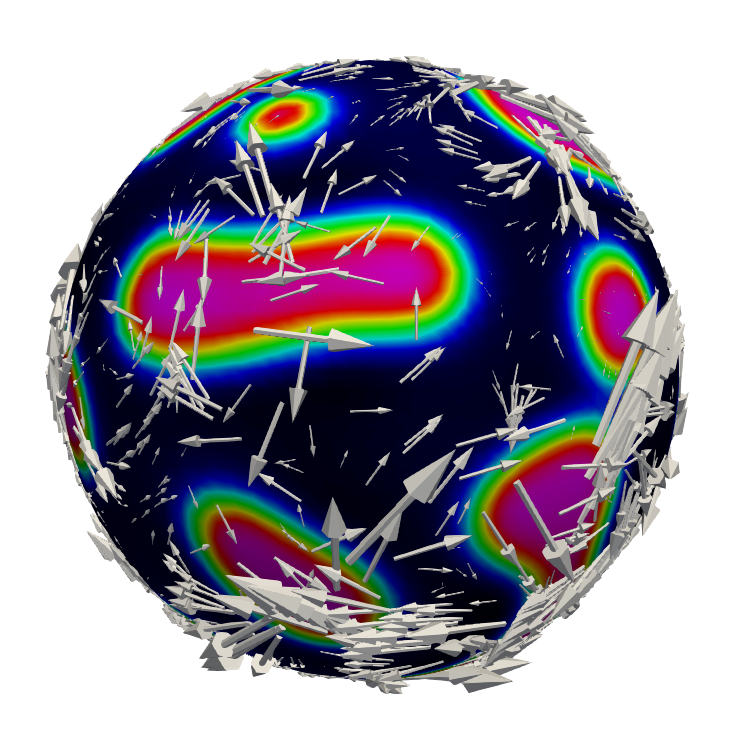}
\end{overpic}
\begin{overpic}[width=.15\textwidth,grid=false]{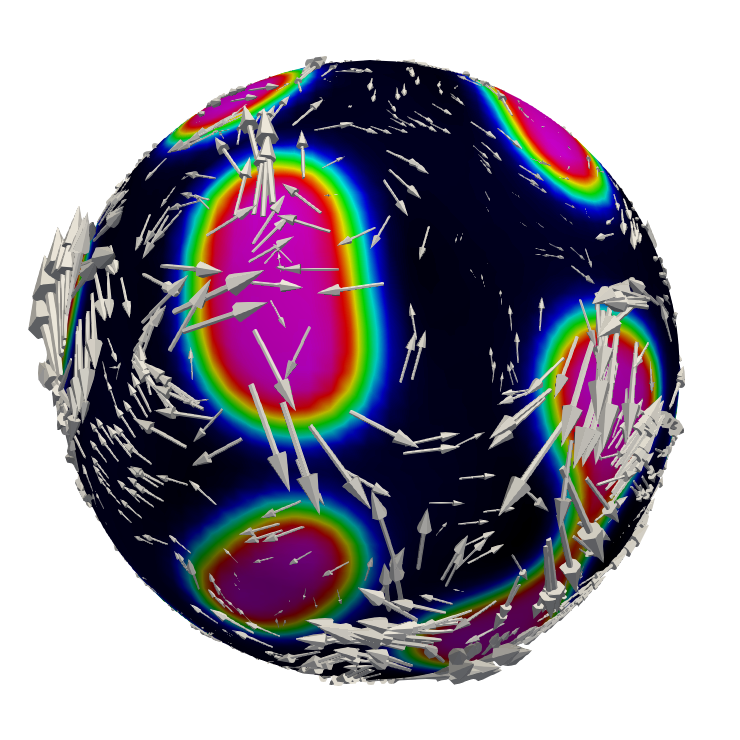}
\end{overpic}
\begin{overpic}[width=.15\textwidth,grid=false]{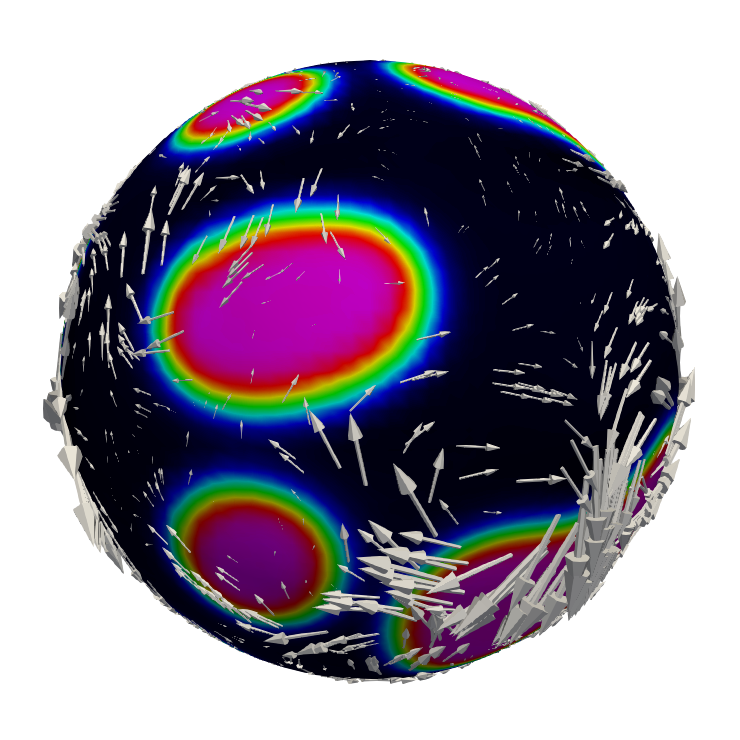}
\end{overpic}
\begin{overpic}[width=.15\textwidth,grid=false]{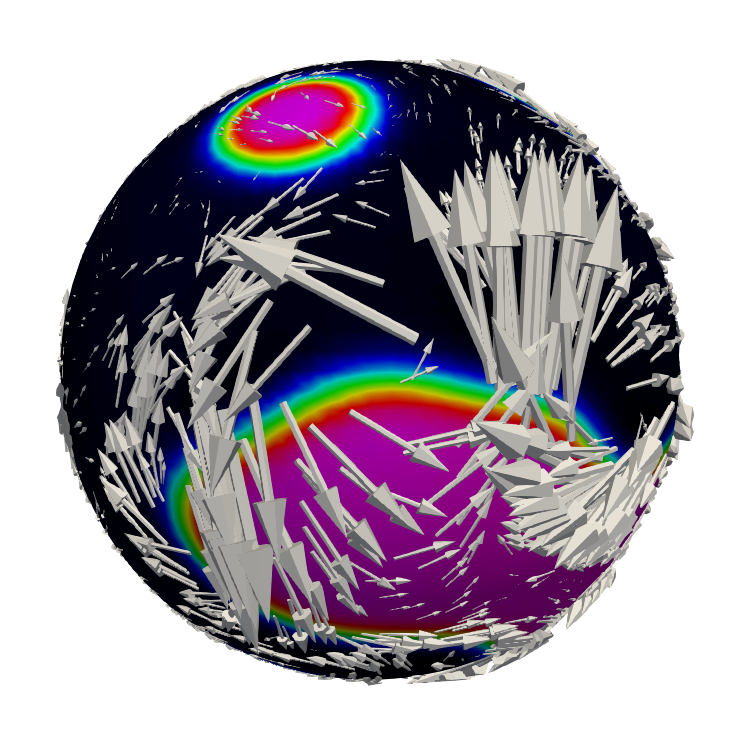}
\end{overpic}
\begin{overpic}[width=.15\textwidth,grid=false]{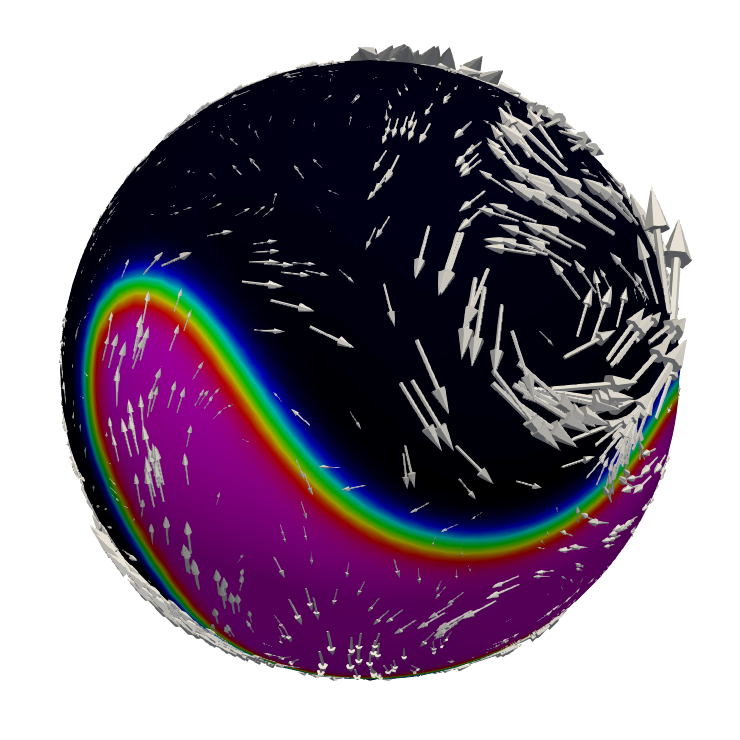}
\end{overpic}
\end{center}
\caption{Velocity vectors superimposed to the surface fraction for $\eta_1 = 10^{-2}, \eta_2 = 8 \cdot 10^{-4}$ 
(top row) and $\eta_1 = 10^{-4}, \eta_2 = 8 \cdot 10^{-6}$ (bottom row) for 
composition 30\%-70\%. For visualization purposes, the  
velocity vectors are magnified by a factor 5 in both rows. 
}\label{fig:flow_Re_30}
\end{figure}

Finally, let us take a look at the fluid flow in Fig.~\ref{fig:flow_Re_50} and Fig.~\ref{fig:flow_Re_30}
for compositions 50\%-50\% and 30\%-70\%, respectively. In both figures, the velocity vectors have been 
magnified by a factor 5. We see that the velocity magnitude in the bottom row of both figures is larger than in the top row
for every time under consideration, as one would expect
when inertial forces become more dominant over viscous forces. 

\subsection{Phase separation on a torus} \label{sec:torus}

We consider a more complex surface than the sphere used so far. We choose
an asymmetric torus with constant distance from the center of the tube to the origin
$R=1$ and variable radius of the tube: $r_{min}=0.3\leq r(x,y) \leq r_{max}=0.6$, with
$r(x,y) =r_{min} + 0.5 (r_{max} - r_{min})  (1 - \frac{x}{\sqrt{x^2 + y^2}})$.
We characterize the torus surface as the zero level set of function $\phi = (x^2 + y^2 + z^2 + R^2 - r(x,y)^2)^2 - 4 R^2 (x^2 + y^2)$.
The torus is embedded in an outer domain $\Omega=[-5/3,5/3]^3$, just like the sphere.
We also selected same mesh level, i.e. $l=5$. 

Like in Sec.~\ref{sec:variable_gamma}, we focus on composition 50\%-50\% and select
$\rho_1 = 3$, $\rho_2 = 1$. Line tension is set to $\sigma_\gamma = 0.04$ and we consider that same
high viscosity and low viscosity cases as in Sec.~\ref{sec:variable_visc_sphere}. 

\begin{figure}[htb]
\centering
\includegraphics[width=0.48\textwidth]{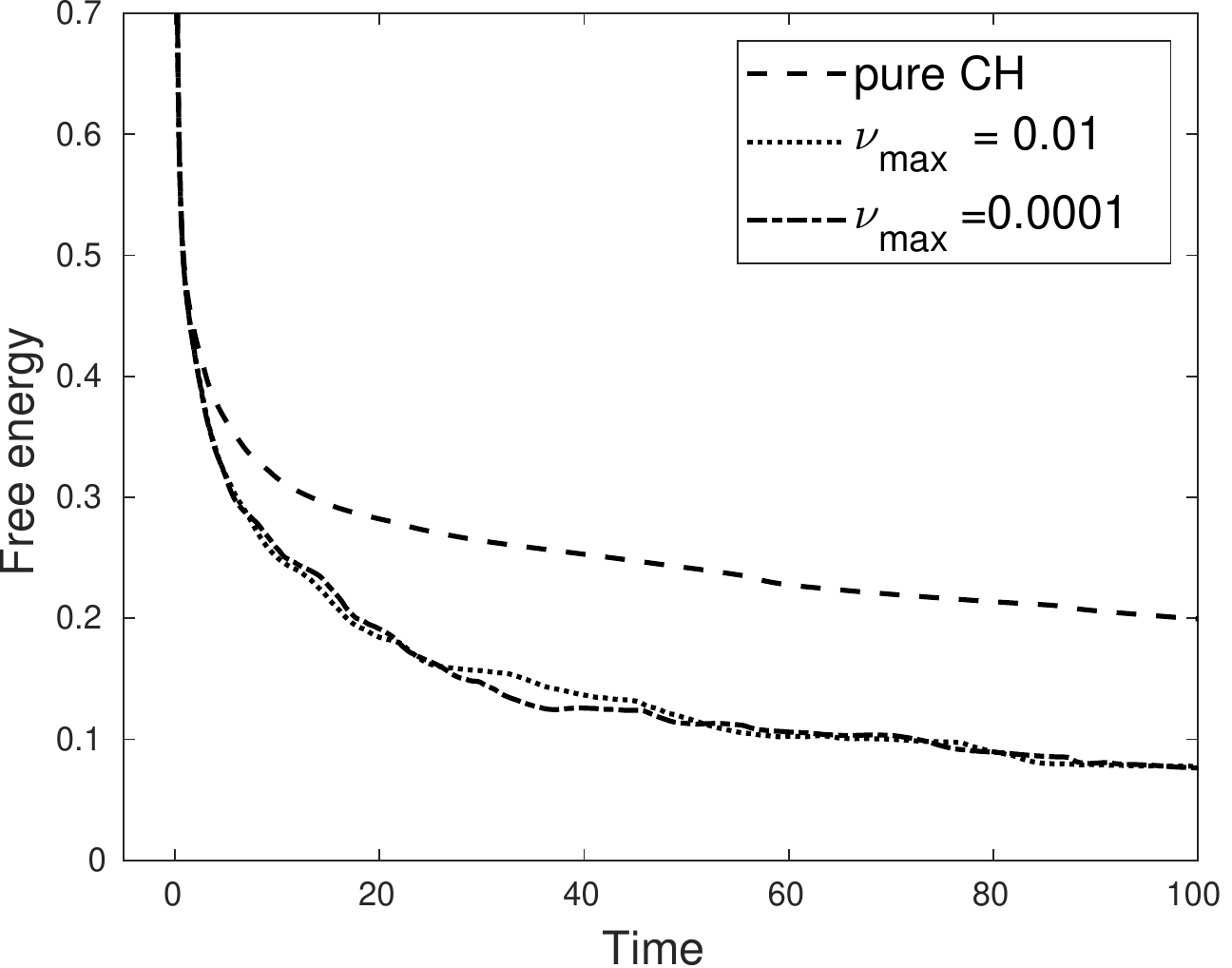}
\includegraphics[width=0.48\textwidth]{images/1.pdf}
	\caption{Discrete Lyapunov energy \eqref{eq:Lyapunov_E} given by the CH model, NSCH model with high values of viscosities
	($\eta_1 = 0.01$, $\eta_2 = 0.0008$), and NSCH model with low values of viscosities ($\eta_1 = 0.0001$, $\eta_2 = 0.000008$)
	on the torus (left) and sphere (right).}
	\label{fig:free_energy_torus}		
\end{figure}

Fig.~\ref{fig:free_energy_torus} reports the discrete Lyapunov energy \eqref{eq:Lyapunov_E}
over time computed by the CH model, NSCH model with low and high viscosities on the torus (left)
 and sphere (right). Fig.~\ref{fig:free_energy_torus} (right) is the same as Fig.~\ref{fig:free_energy_eta} (left); it is reported 
 again to facilitate the comparison. On the torus, just like on the sphere, 
the presence of surface flow leads to a faster Lyapunov energy decay. Moreover, on both surfaces
switching from high to low values of the viscosity does not produce a significant
change in the Lyapunov energy decay. On the sphere though, the energy drops to a lower value after the initial 
fast phase of phase separation and flattens faster in the subsequent slower phase. This suggests an effect of the
surface geometry on the evolution of phases.

\begin{figure}
\begin{center}
\hskip .7cm
\begin{overpic}[width=.15\textwidth,grid=false]{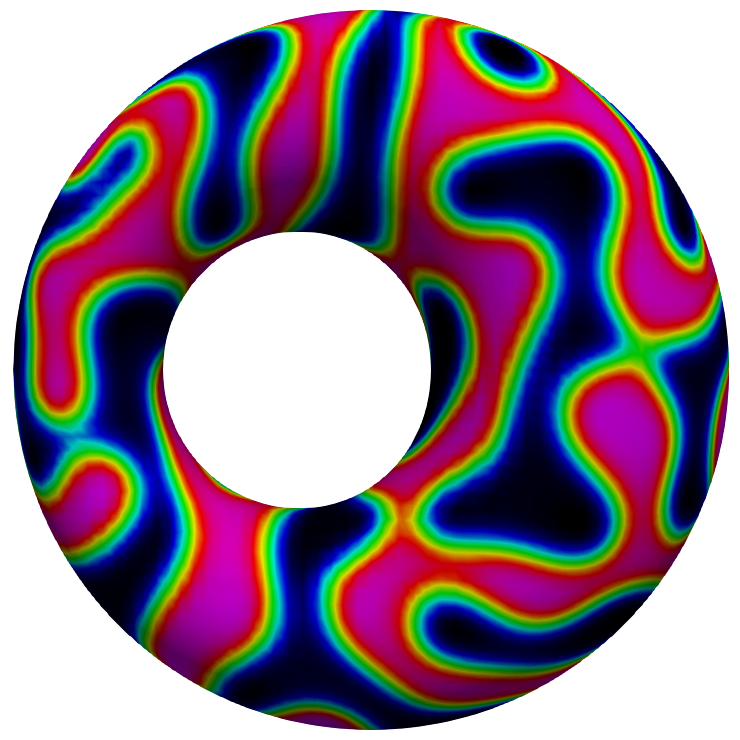}
\put(30,100){\small{$t = 2$}}
\put(-30,50){\small{CH}}
\end{overpic}
\begin{overpic}[width=.15\textwidth,grid=false]{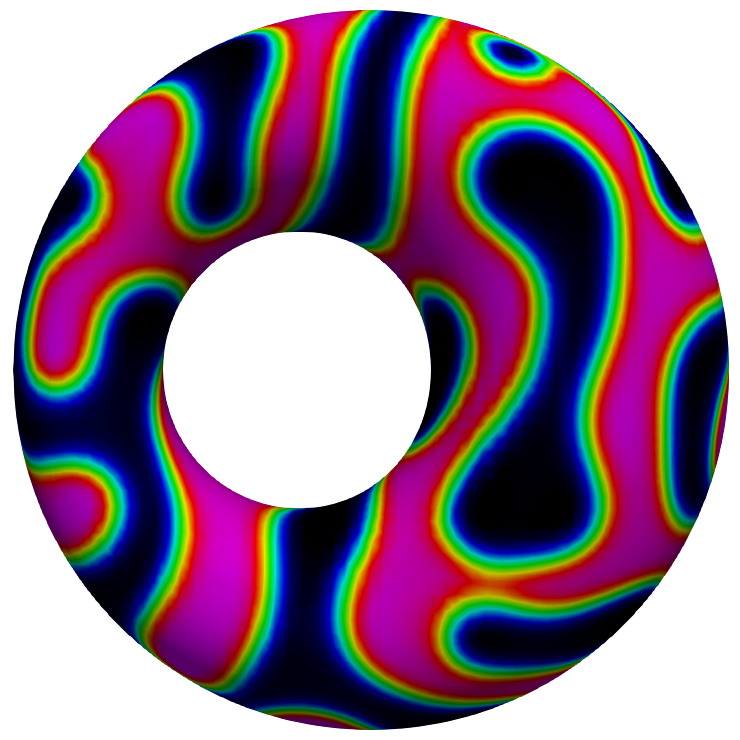}
\put(30,100){\small{$t = 5$}}
\end{overpic}
\begin{overpic}[width=.15\textwidth,grid=false]{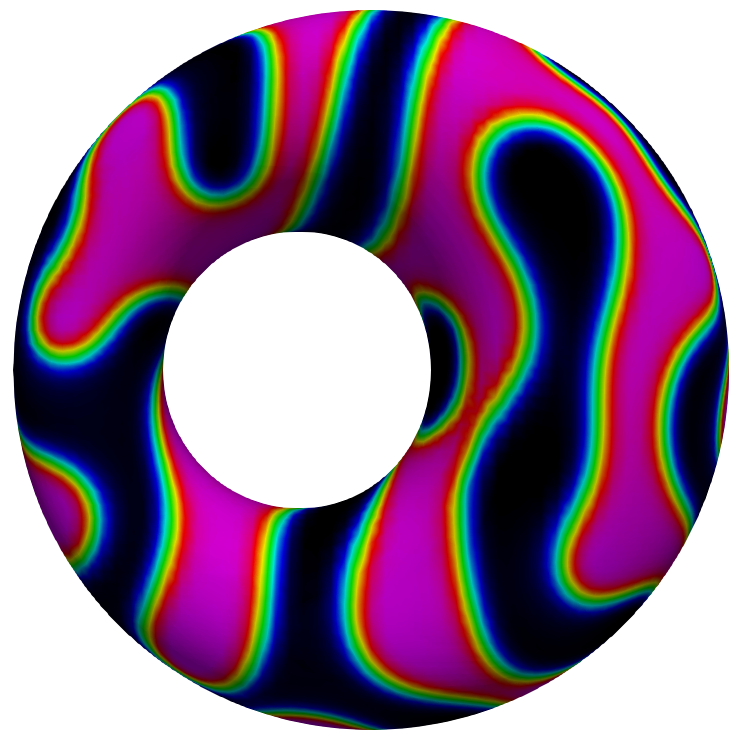}
\put(28,100){\small{$t = 10$}} 
\end{overpic}
\begin{overpic}[width=.15\textwidth,grid=false]{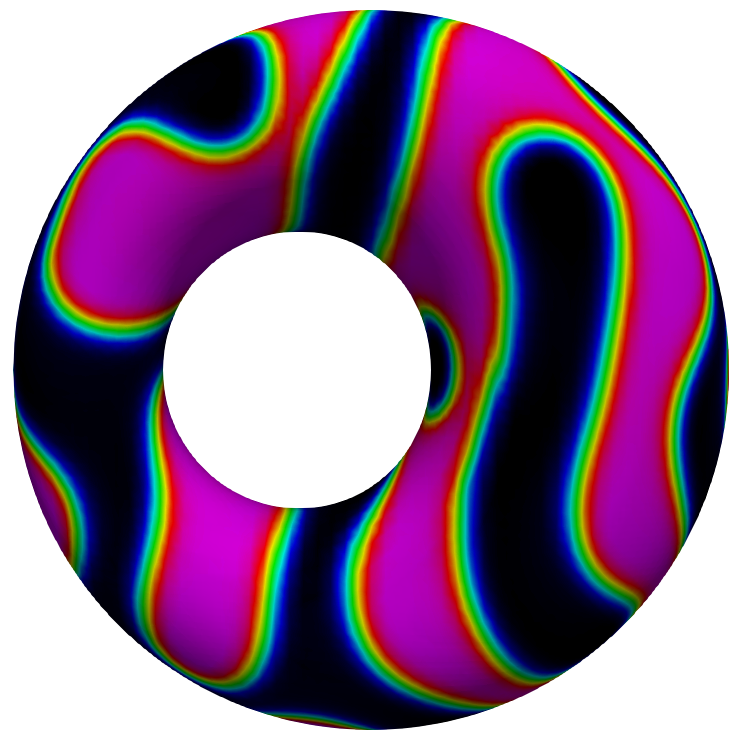}
\put(28,100){\small{$t = 20$}}
\end{overpic}
\begin{overpic}[width=.15\textwidth,grid=false]{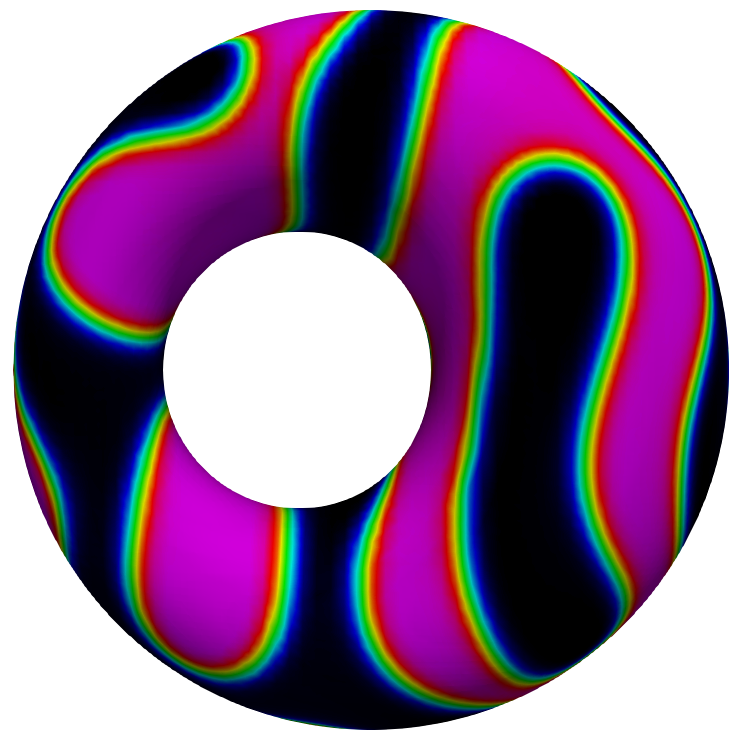}
\put(28,100){\small{$t = 40$}}
\end{overpic}
\begin{overpic}[width=.15\textwidth,grid=false]{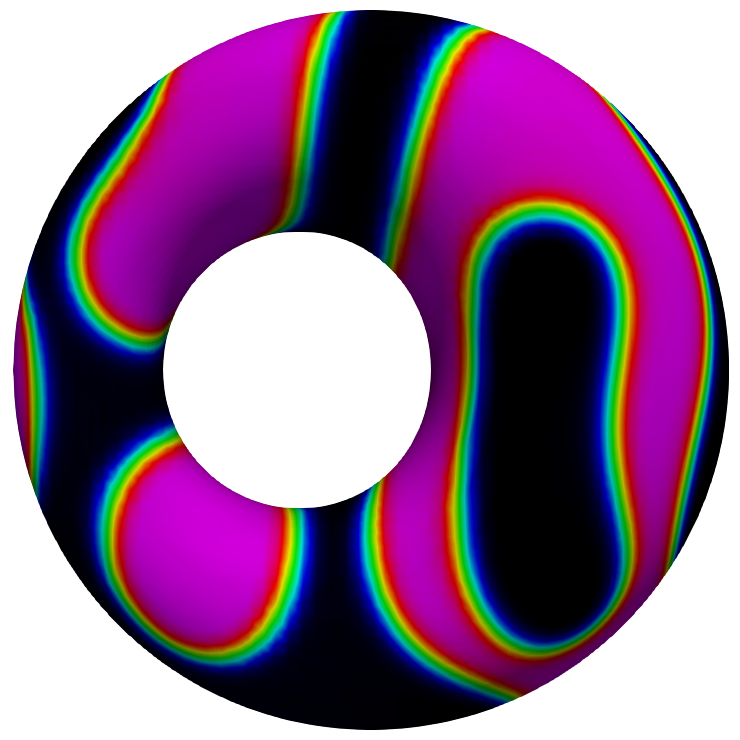}
\put(26,100){\small{$t = 100$}}
\end{overpic}
\\
\hskip .7cm
\begin{overpic}[width=.15\textwidth,grid=false]{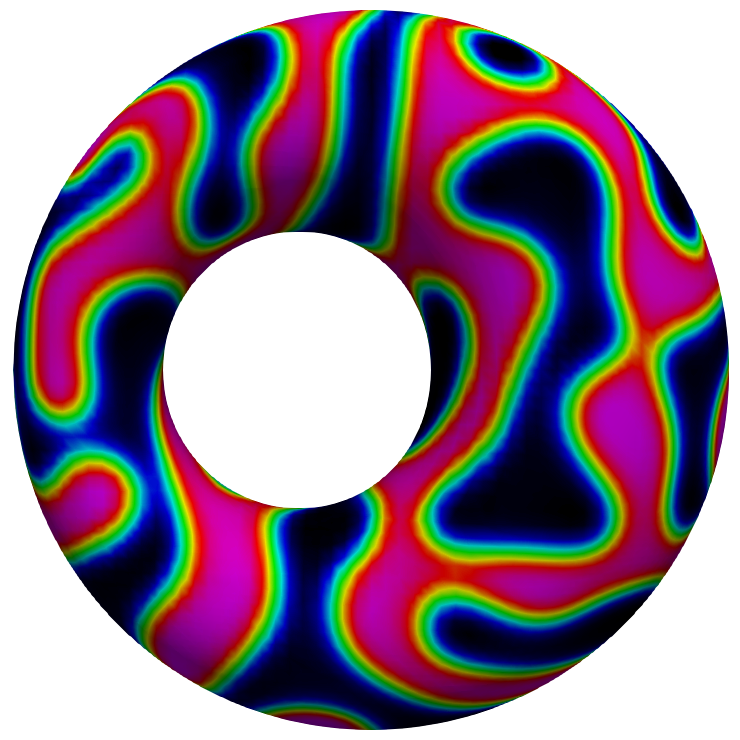}
\put(-50,50){\small{NSCH}}
\put(-50,33){\small{high $\eta$}}
\end{overpic}
\begin{overpic}[width=.15\textwidth,grid=false]{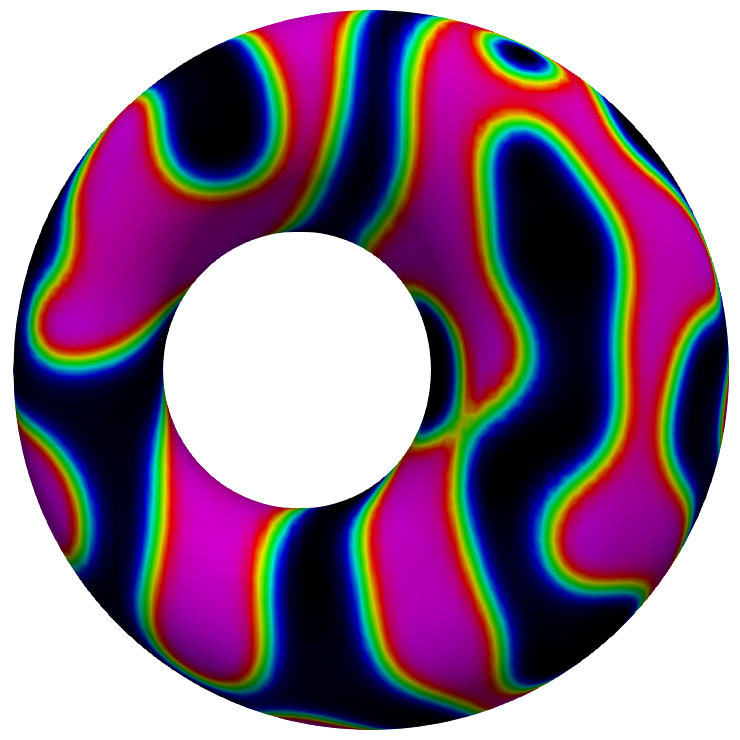}
\end{overpic}
\begin{overpic}[width=.15\textwidth,grid=false]{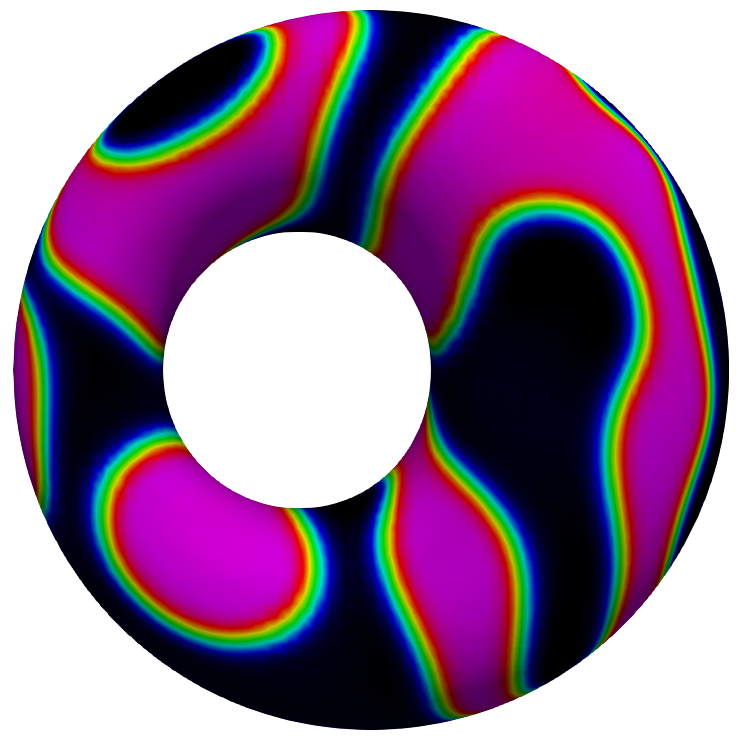}
\end{overpic}
\begin{overpic}[width=.15\textwidth,grid=false]{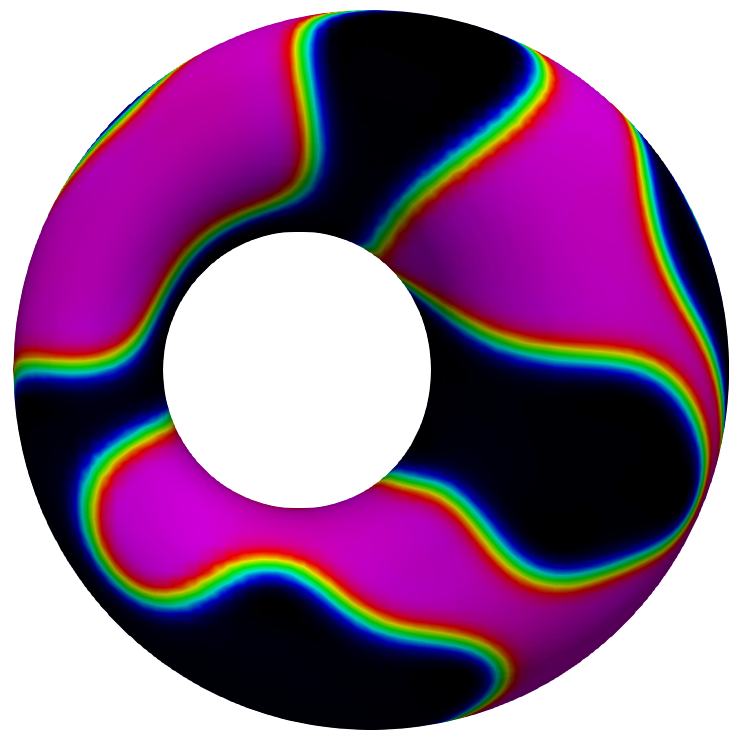}
\end{overpic}
\begin{overpic}[width=.15\textwidth,grid=false]{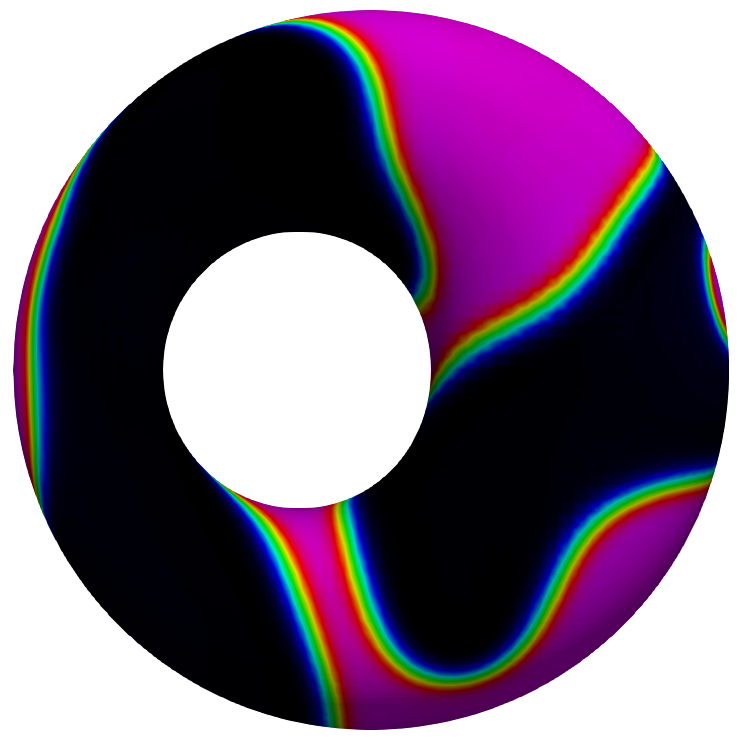}
\end{overpic}
\begin{overpic}[width=.15\textwidth,grid=false]{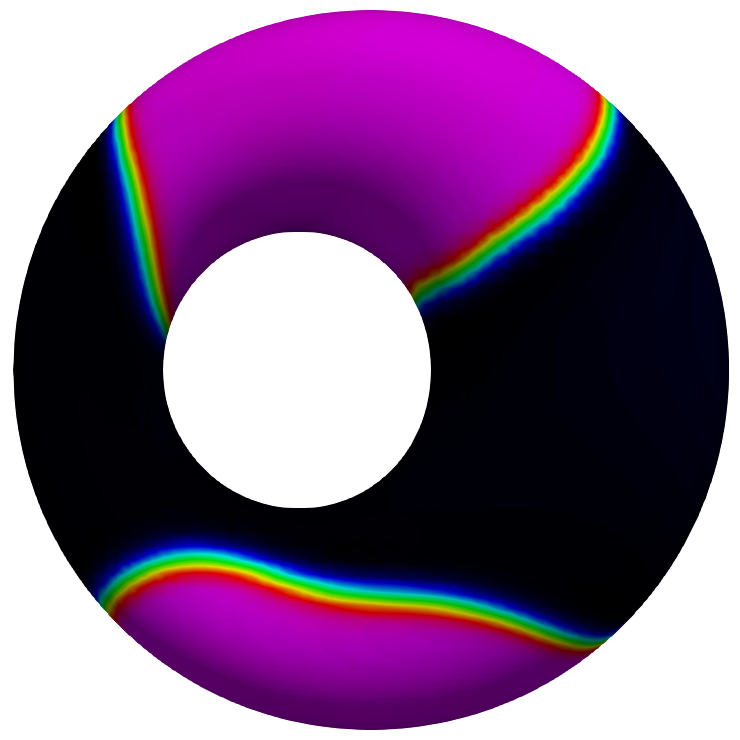}
\end{overpic}
\\
\hskip .7cm
\begin{overpic}[width=.15\textwidth,grid=false]{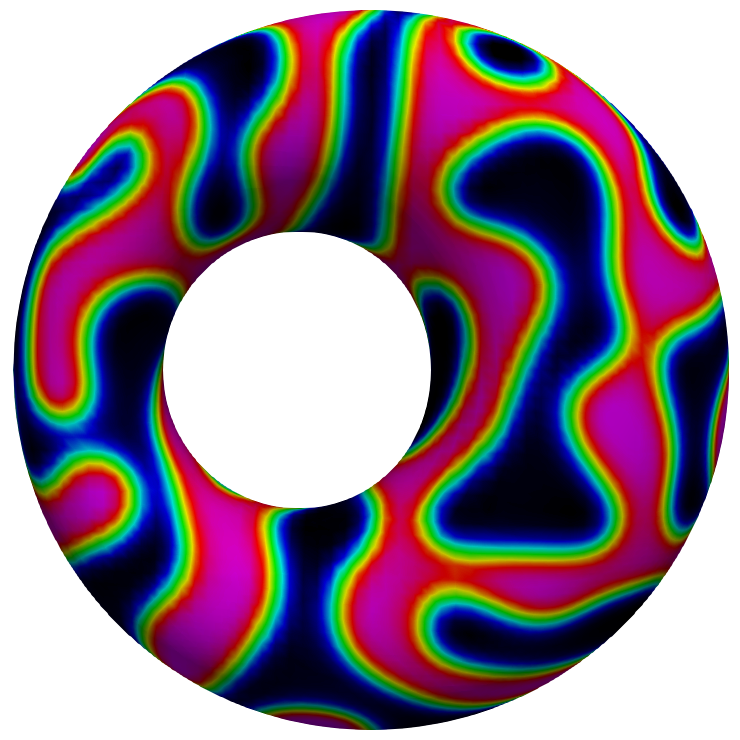}
\put(-50,50){\small{NSCH}}
\put(-50,33){\small{low $\eta$}}
\end{overpic}
\begin{overpic}[width=.15\textwidth,grid=false]{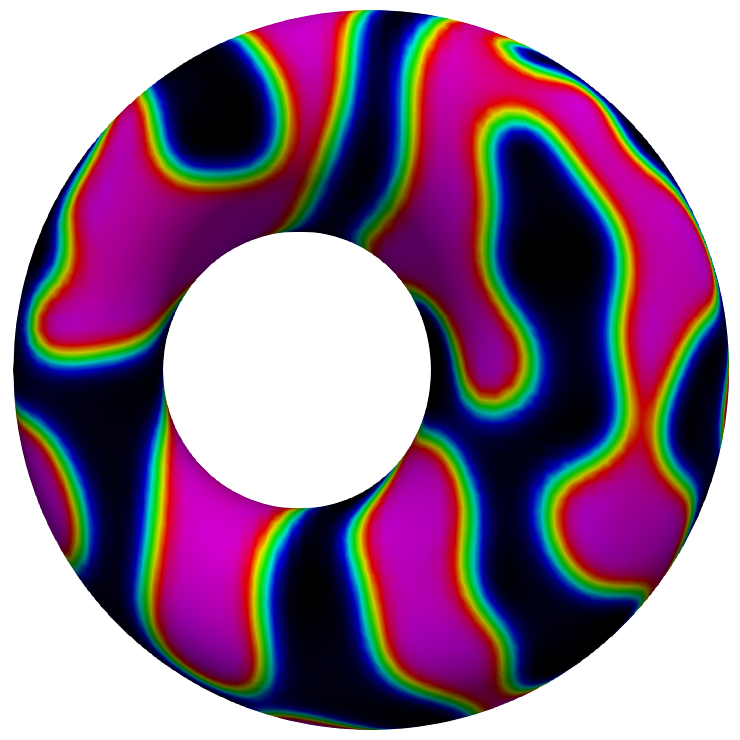}
\end{overpic}
\begin{overpic}[width=.15\textwidth,grid=false]{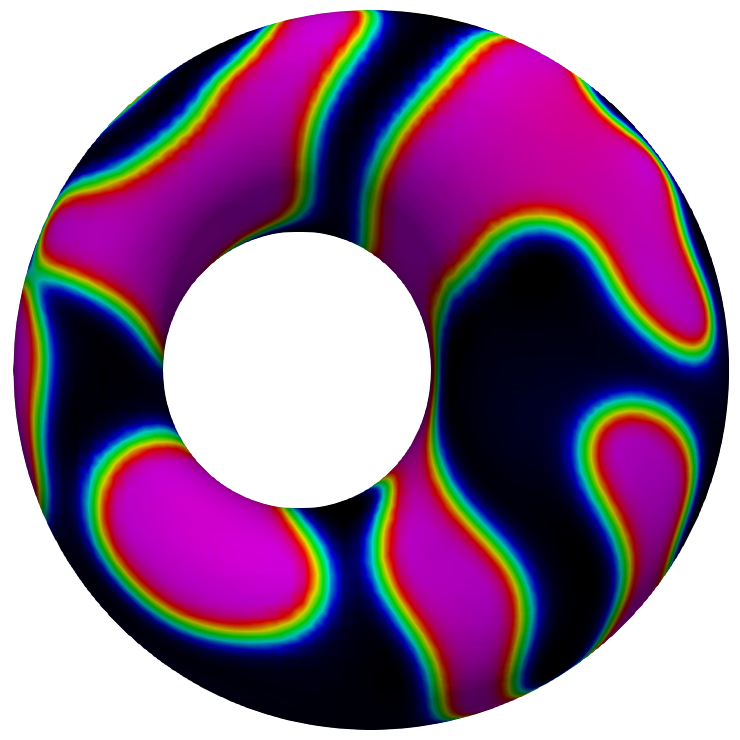}
\end{overpic}
\begin{overpic}[width=.15\textwidth,grid=false]{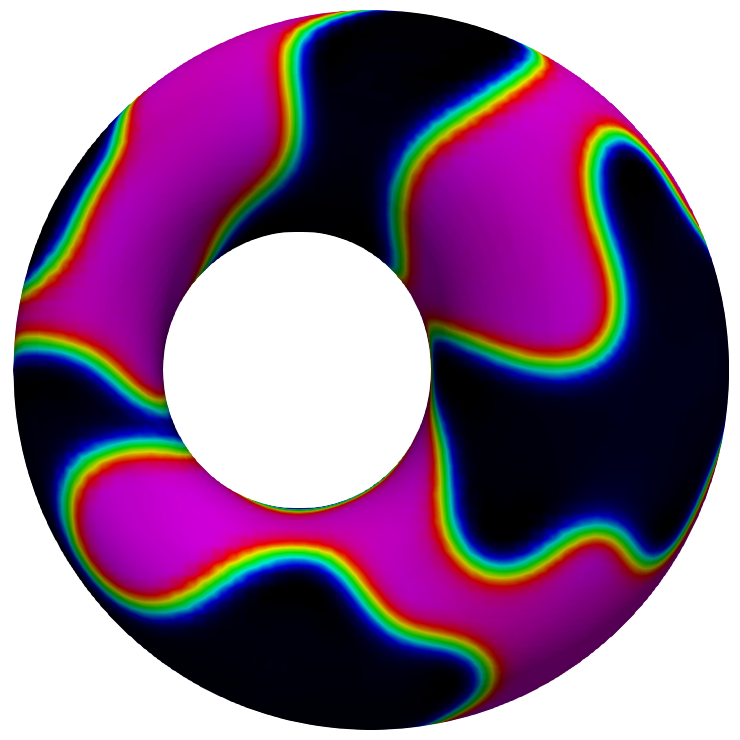}
\end{overpic}
\begin{overpic}[width=.15\textwidth,grid=false]{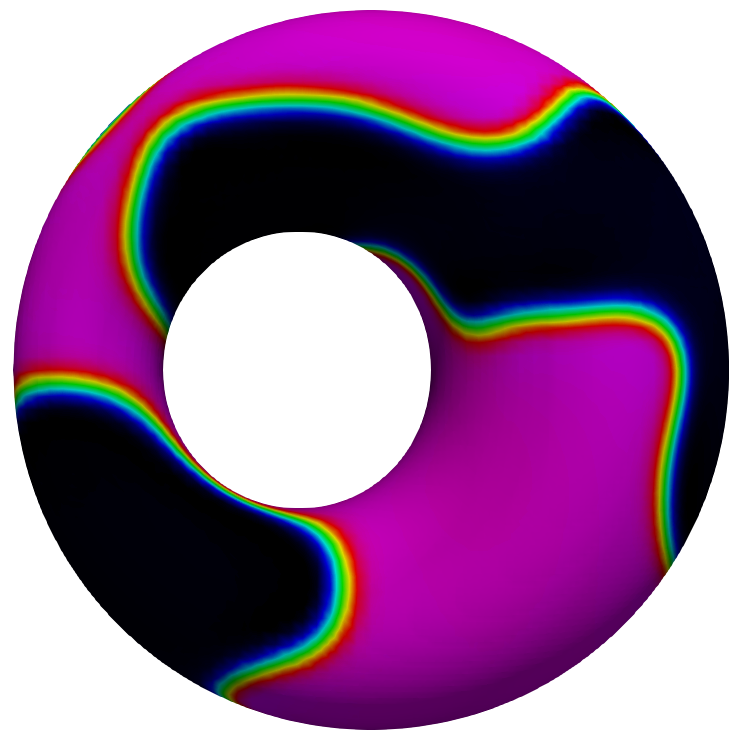}
\end{overpic}
\begin{overpic}[width=.15\textwidth,grid=false]{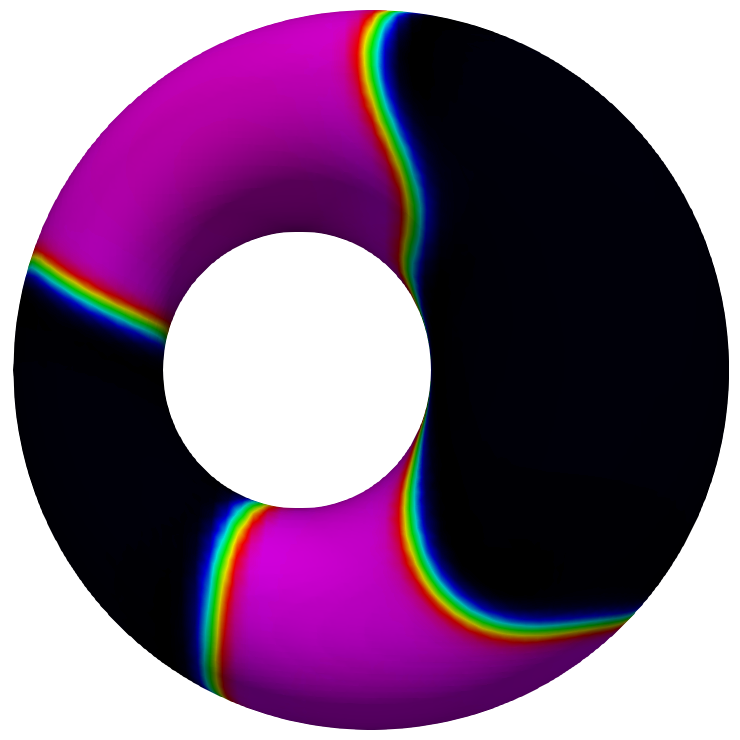}
\end{overpic}
\\
\vskip .2cm
\begin{overpic}[width=0.5\textwidth,grid=false,tics=10]{images/legend.png}
\end{overpic}
\end{center}
\caption{Phase separation given by the CH model (top), NSCH model with viscosities
$\eta_1 = 0.01$, $\eta_2 = 0.0008$ (center), and NSCH model with viscosities $\eta_1 = 0.0001$, $\eta_2 = 0.000008$ (bottom).}\label{fig:torus_Re}
\end{figure}

The effect of the geometry can bee seen also when comparing Fig.~\ref{fig:sphere_Re_50} with Fig.~\ref{fig:torus_Re}, 
which shows the evolution of phases delivered by the CH model and NSCH model for the high viscosity and low viscosity cases
on the torus.
We see that the interface separating the two phases remains tortuous for a longer period of time on the torus. 
As  within the torus itself, we do not observe a particular difference in pattern between ``skinny'' and ``fat''
side of the torus.

\begin{figure}
\begin{center}
\hskip .7cm
\begin{overpic}[width=.15\textwidth,grid=false]{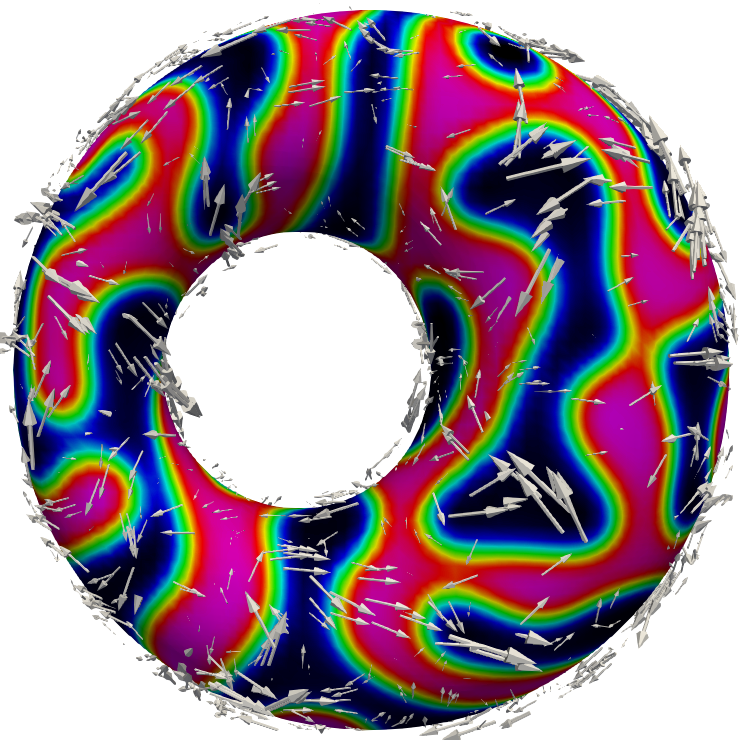}
\put(30,100){\small{$t = 2$}}
\put(-50,45){\small{high $\eta$}}
\end{overpic}
\begin{overpic}[width=.15\textwidth,grid=false]{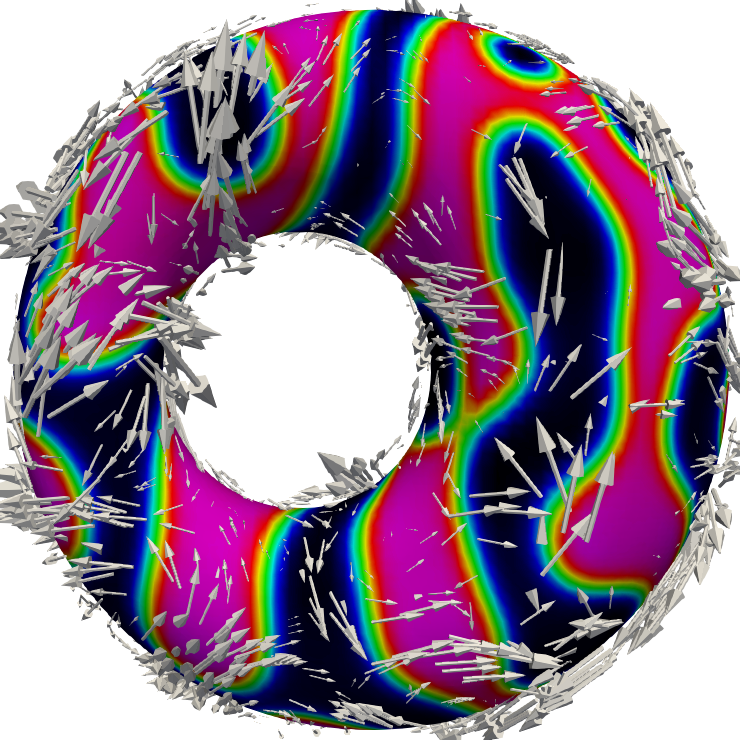}
\put(30,100){\small{$t = 5$}}
\end{overpic}
\begin{overpic}[width=.15\textwidth,grid=false]{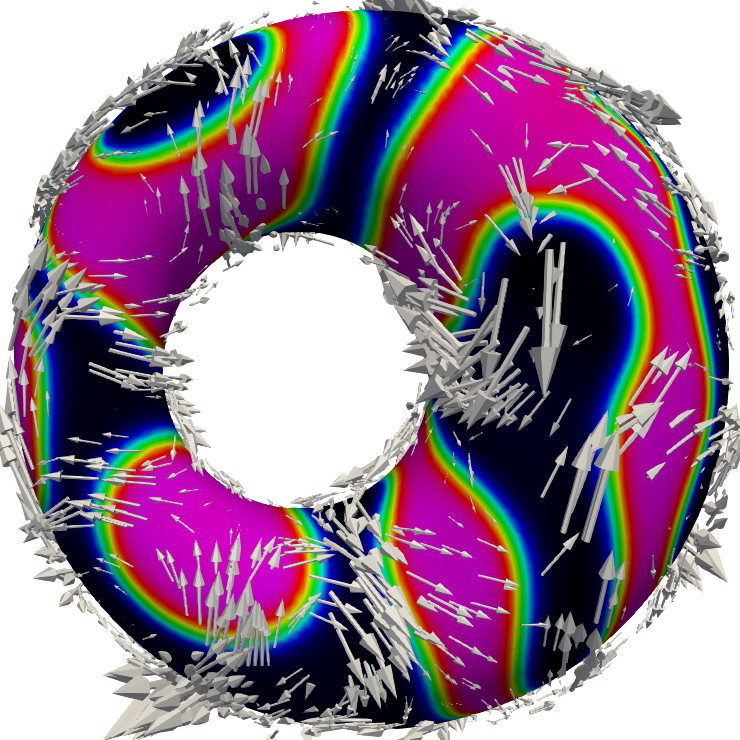}
\put(28,100){\small{$t = 10$}} 
\end{overpic}
\begin{overpic}[width=.15\textwidth,grid=false]{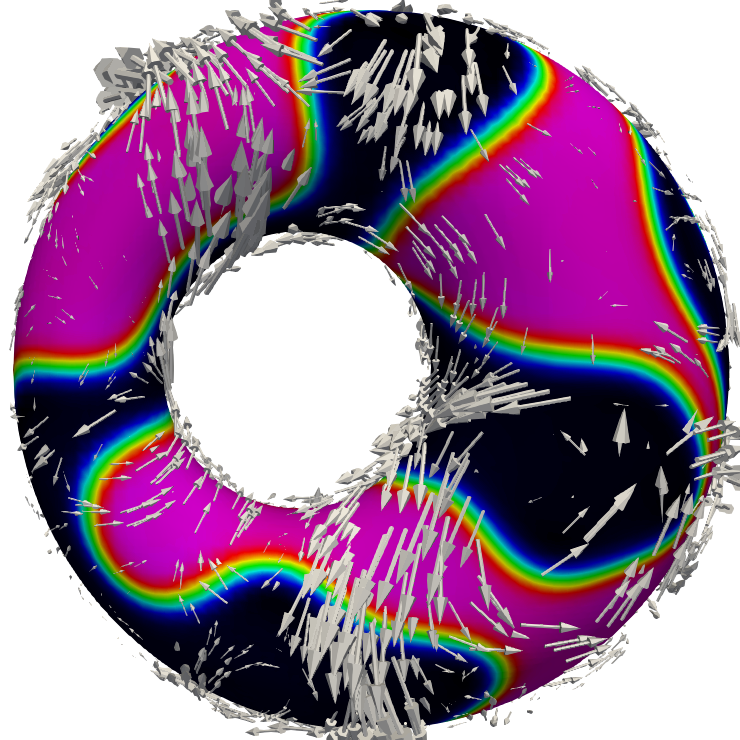}
\put(28,100){\small{$t = 20$}}
\end{overpic}
\begin{overpic}[width=.15\textwidth,grid=false]{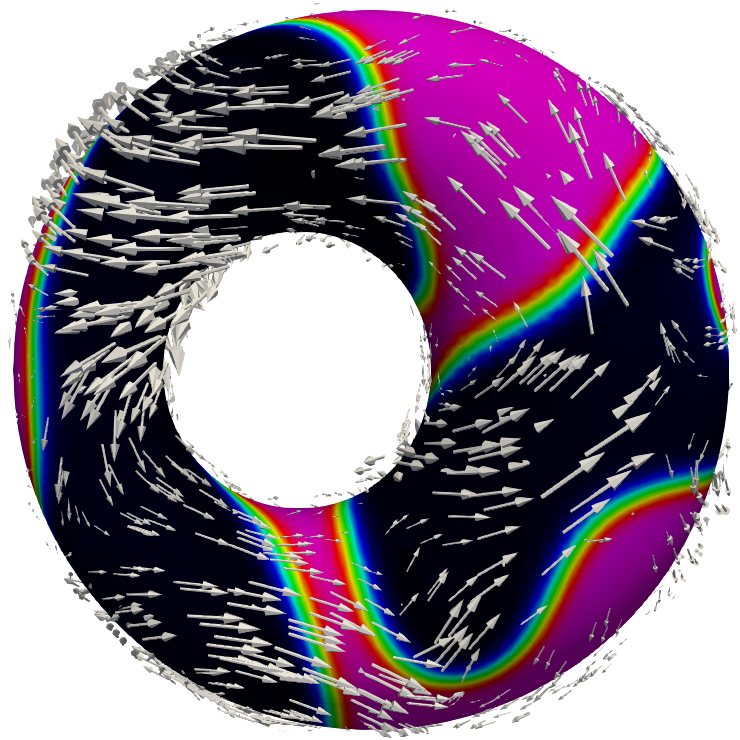}
\put(28,100){\small{$t = 40$}}
\end{overpic}
\begin{overpic}[width=.15\textwidth,grid=false]{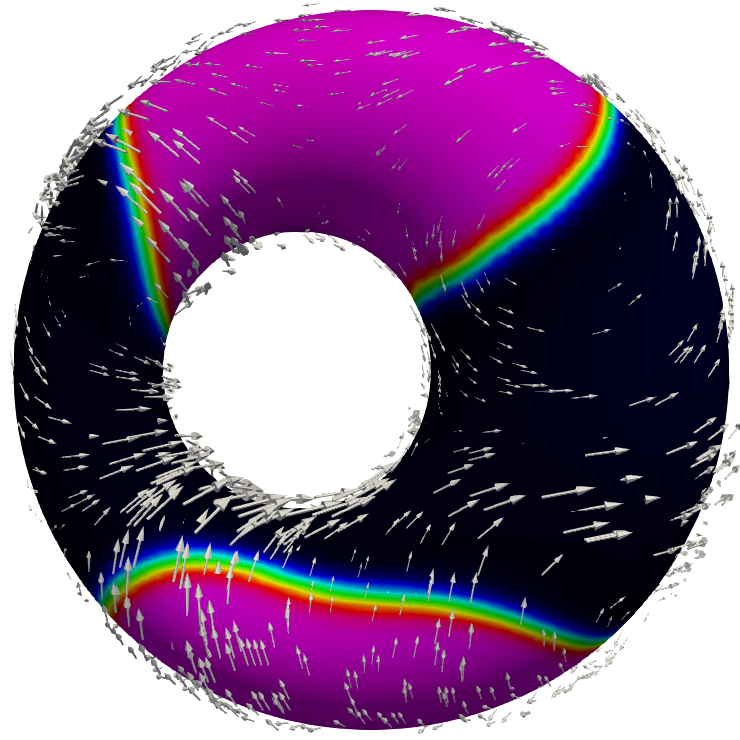}
\put(26,100){\small{$t = 100$}}
\end{overpic}
\\
\hskip .7cm
\begin{overpic}[width=.15\textwidth,grid=false]{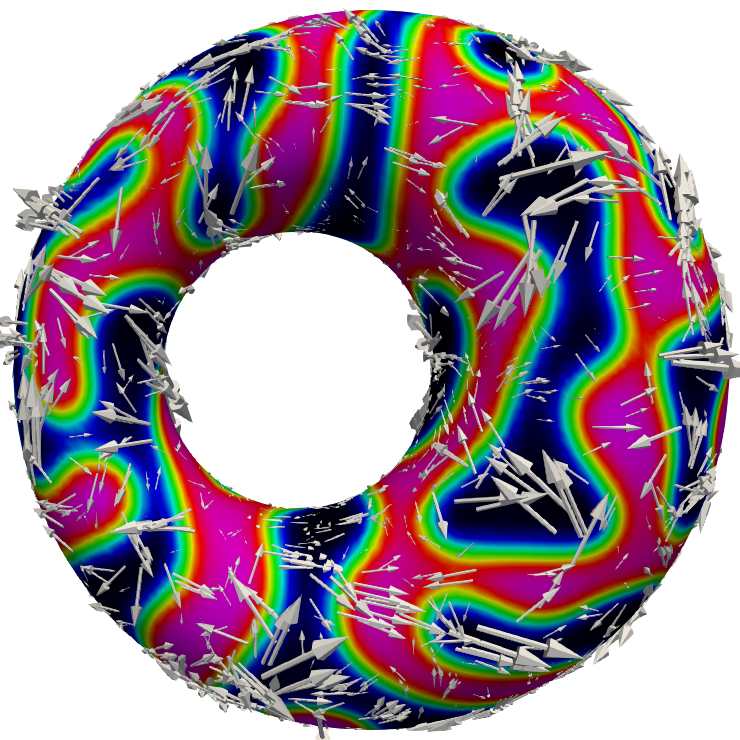}
\put(-50,45){\small{low $\eta$}}
\end{overpic}
\begin{overpic}[width=.15\textwidth,grid=false]{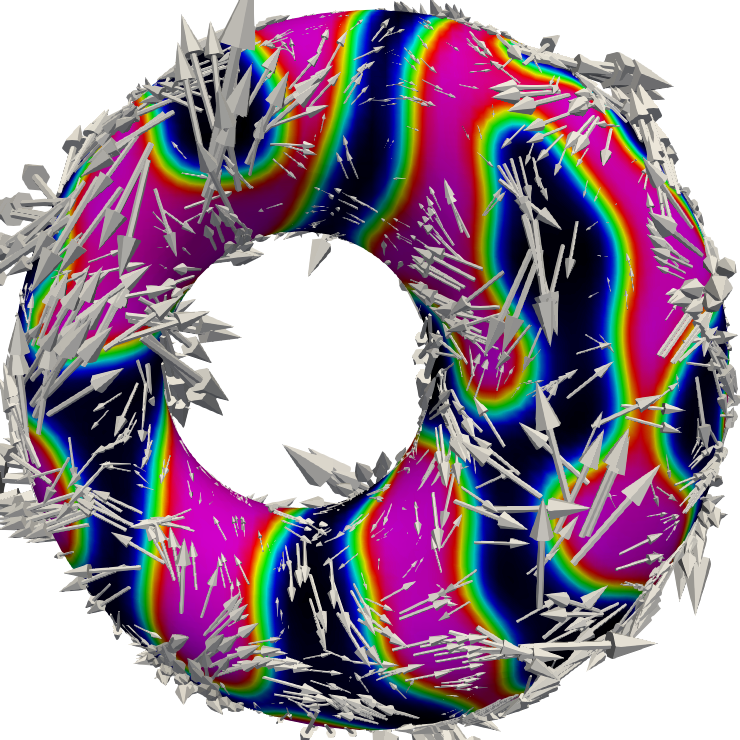}
\end{overpic}
\begin{overpic}[width=.15\textwidth,grid=false]{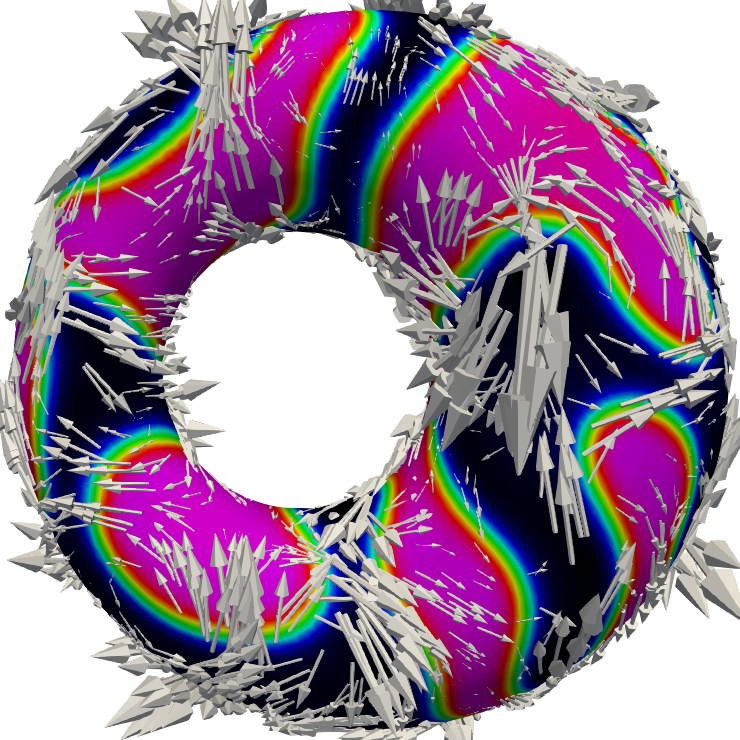}
\end{overpic}
\begin{overpic}[width=.15\textwidth,grid=false]{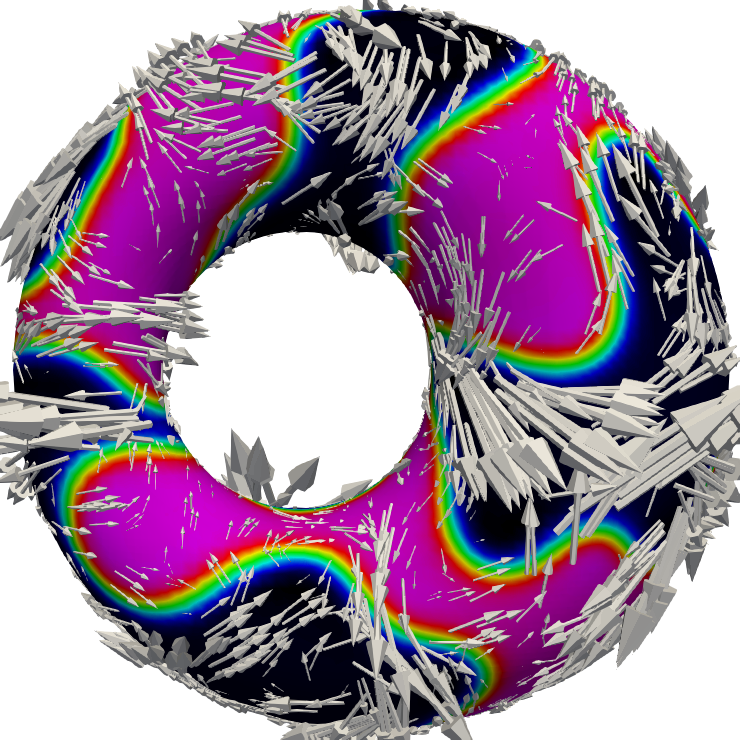}
\end{overpic}
\begin{overpic}[width=.15\textwidth,grid=false]{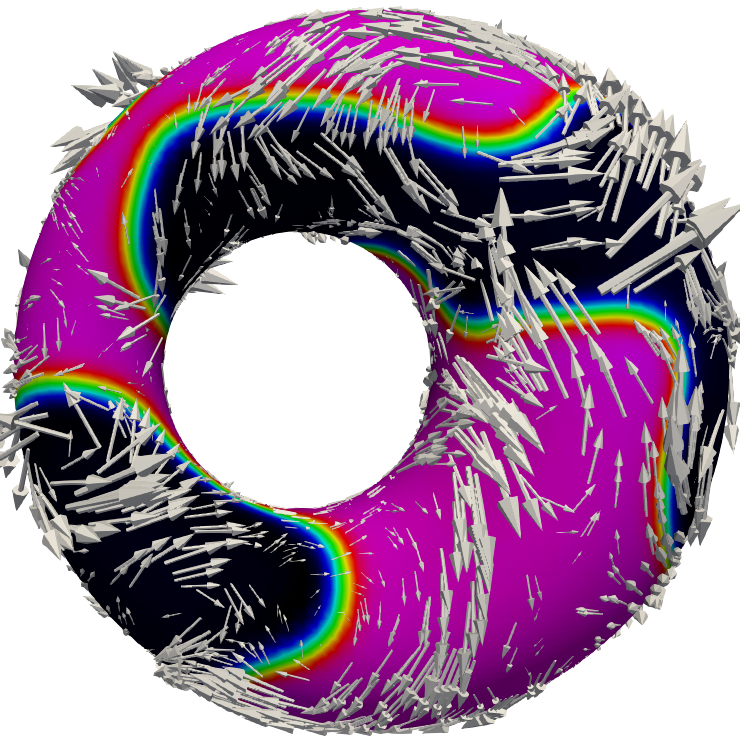}
\end{overpic}
\begin{overpic}[width=.15\textwidth,grid=false]{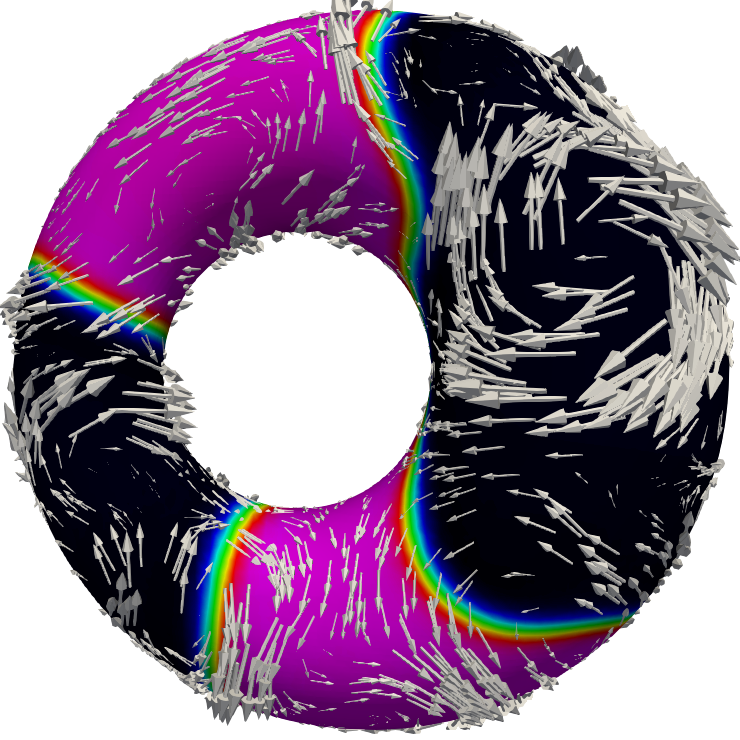}
\end{overpic}
\end{center}
\caption{Velocity vectors superimposed to the surface fraction for $\eta_1 = 10^{-2}, \eta_2 = 8 \cdot 10^{-4}$ 
(top row) and $\eta_1 = 10^{-4}, \eta_2 = 8 \cdot 10^{-6}$ (bottom row). For visualization purposes, the  
velocity vectors are magnified by a factor 5 in both rows.}\label{fig:flow_torus_Re}
\end{figure}

Next, in Fig.~\ref{fig:flow_torus_Re} we report the the velocity vectors superimposed to the surface fraction
for the bottom two cases in Fig.~\ref{fig:torus_Re}. The velocity vectors have been 
magnified by a factor 5. This allows us to compare the fluid flows on the sphere 
(in Fig.~\ref{fig:flow_Re_50}) and the torus. We observe more intricate flow patters on the torus
due to both the more complex shape and the persistence of the tortuosity in the interface
separating the phsaes.

One final note is in order.
Although the analysis of the decoupled scheme for the NSCH problem in Sec.~\ref{sec:method} 
presented in \cite{Palzhanov2021} points to some restrictions on the discretization parameters,
in practice we do not observe any restriction on  the time step, while for the mesh size 
we observed that it needs to resolve the transition layer between phases (as one would easily guess).

\section{Conclusions}\label{sec:concl}

We performed a computational study of lateral phase separation and coarsening on surfaces. 
To model these processes, we considered both the Cahn--Hilliard (phase separation alone) and the
Navier--Stokes--Cahn--Hilliard (phase separation coupled to lateral flow) equations posed on manifolds.
Both models were solved numerically using an unfitted finite element method called TraceFEM, which
allows for a flexible treatment of complex and evolving surfaces. This choice is motivated by our
interest in the computational design of lipid membranes used as drug carriers.

Through a series of numerical tests on the surface of a sphere and an asymmetric torus, we investigated
how the evolution of phases changes when switching from the Cahn--Hilliard (CH) model to the 
Navier--Stokes--Cahn--Hilliard (NSCH) model with variable line tension, viscosity, and membrane composition.
We observed that the discrete Lypunov energy decays faster when using the NSCH model. In particular, such faster
decay is more significant when the line tension is increased and the viscosity is lowered. The latter is more
evident in some membrane compositions (i.e., 30\%-70\%) than others (i.e., 50\%-50\%). Finally, by comparing
the evolution of phases on the sphere and on the torus we do observe differences that indicate an effect of the surface geometry.

\bmhead{Acknowledgments}

This work was partially supported by US National Science Foundation (NSF) through grant DMS-1953535.
 M.O.~also acknowledge the support from NSF through DMS-2011444.
A.Q.~also acknowledges support from the Radcliffe Institute for Advanced Study at Harvard University where
she has been a 2021-2022 William and Flora Hewlett Foundation Fellow.

\bibliography{literatur}


\end{document}